\documentclass[a4paper,11pt]{article}

\usepackage[T1]{fontenc}
\usepackage{geometry,	lipsum,	fancyhdr}
\usepackage{pifont, amssymb,amsmath, amsthm, latexsym,esint,thmtools, apptools}
\usepackage{epsfig}
\usepackage{mathrsfs}
\usepackage{epic,eepic,wrapfig,color, ifthen, authblk}
\usepackage{hyperref}

\usepackage{url}
\usepackage{dsfont}

\theoremstyle{plain}
\declaretheorem[title=Theorem, parent=section]{theo}
\declaretheorem[title=Lemma,sibling=theo]{lemma}
\declaretheorem[title=Proposition,sibling=theo]{prop}
\declaretheorem[title=Corollary,sibling=theo]{cor}

\theoremstyle{definition}
\declaretheorem[title=Definition,sibling=theo]{defn}
\declaretheorem[title=Remark,sibling=theo]{remark}

\declaretheorem[title=Remark, numbered=no]{remark*}
\declaretheorem[title=Example, sibling=theo]{example}

\numberwithin{equation}{section}

\renewcommand{\appendix}{\par
	\setcounter{section}{0}
	\setcounter{subsection}{0}
	\gdef\thesection{\Alph{section}}
}

\def\1{{\bf 1}}

\def\nn{\nonumber}

\def\pf{{\medskip\noindent {\bf Proof. }}}
\def\qed{{\hfill $\Box$ \bigskip}}

\def\sA {{\mathcal A}} \def\sB {{\mathcal B}} \def\sC {{\mathcal C}}
\def\sD {{\mathcal D}}
\def\sE {{\mathcal E}} 
\def\sF {{\mathcal F}}
  \def\sI {{\mathcal I}}
\def\sJ {{\mathcal J}} 
\def\sL {{\mathcal L}}

\def\sN {{\mathcal N}}

\def\sT {{\mathcal T}} 
\def\sV {{\mathcal V}}  
 
\def\R {{\mathbb R}} 
\def\N{{\mathbb N}} 
\def\E {{\mathbb E}}  \def \P{{\mathbb P}}

\def \diam{{\textrm{\rm diam}}}

\def \NT{\sT}

\def \EC{${\textrm{\rm rEC}}_\beta(s)$}
\def \TJ{${\textrm{\rm rTJ}}_\beta(s)$}
\def \PI{${\textrm{\rm rPI}}_\beta(s)$}
\def \FK{${\textrm{\rm rFK}}_\beta(s)$}
\def \CS{${\textrm{\rm rCS}}_\beta(s)$}
\def \WEHI{${\textrm{\rm rWEHI}}_\beta(s)$}
\def \EHR{${\textrm{\rm rEHR}}_\beta(s)$}

\def\la{{\langle}}
\def\ra{{\rangle}}

\def\qed{{\hfill $\Box$ \bigskip}}
\def\eps{\varepsilon}

\def\wt{\widetilde}
\def\pf{\noindent{\bf Proof. }}

\def\vp{{\varphi}}

\DeclareMathOperator*{\esssup}{ess\,sup}
\DeclareMathOperator*{\essinf}{ess\,inf}
\DeclareMathOperator*{\essosc}{ess\,osc}
\makeatletter
\@addtoreset{equation}{section}

\makeatother

	\title{Robust  estimates for elliptic nonlocal operators on doubling spaces}
	
	\author[1]{Soobin Cho}
	\affil[1]{Department of Mathematics, University of Illinois Urbana-Champaign, \protect\\ 
		Urbana, IL 61801, USA. 
		Email: soobinc@illinois.edu}
	
\begin{document}
		
		\date{}
		\maketitle

		\begin{abstract}
			We study weak Harnack inequality and a priori H\"older regularity of harmonic functions for symmetric nonlocal Dirichlet forms on  metric measure spaces with volume doubling condition. Our analysis relies on three main assumptions: the existence of a strongly local Dirichlet form with sub-Gaussian heat kernel estimates, a tail estimate of the jump measure outside balls and a local energy comparability condition. We establish the robustness of our results, ensuring that the constants in our estimates remain bounded, provided that the order of the scale function appearing in  the tail estimate and local energy comparability condition, maintains a certain distance from zero. Additionally, we establish a sufficient condition for the local energy comparability condition.
		\end{abstract}

		\medskip

		\normalsize
		
		\tableofcontents

	\section{Introduction and Main results}\label{s:intro}

\subsection{Introduction}

Regularity results for nonlocal  Dirichlet forms  are of significant importance in both analysis and probability theory.  Extensive studies have been dedicated to this area. See, for example, \cite{BaCh10,BaKa05,BaLe02, BGHH21,BSS03,ChKu03,CKW-elp,DKP14,Ka09,RoSe16,Si06}. A typical example of a nonlocal regular Dirichlet form $(\sE,\sF)$ on a metric measure $(M,d,\mu)$ is defined by
\begin{align*}
	\sE(f,g)= \int_{M\times M} (f(x)-f(y))(g(x)-g(y))J(x,y)\mu(dx)\mu(dy), \quad f,g \in \sF,
\end{align*}
where $J(x,y)$ is a symmetric measurable function on $M\times M$ satisfying
\begin{align}\label{e:jumpdensity-pointwise}
 \frac{C^{-1}}{\mu(B(x,d(x,y))) d(x,y)^{\beta s}}\le 	J(x,y) \le \frac{C}{\mu(B(x,d(x,y))) d(x,y)^{\beta s}} \quad \text{for all $x,y \in M$,}
\end{align}
with  $C>1$ and $s\in (0,1)$. Here  $\beta\ge 2$ represents the \textit{walk dimension} of $M$. This includes the Dirichlet form associated with the fractional Laplacian $-(-\Delta)^{s}$ on $\R^d$. Notably, while the walk dimension of the Euclidean space is $2$, many fractals are known to have the walk dimension $\beta>2$, potentially resulting in $\beta s$ exceeding $2$.

When $\beta s > 2$, a significant challenge in obtaining regularity estimates arises due to the fact that the domain of the Dirichlet form may not include all Lipschitz functions. This necessitates careful construction of cutoff test functions a priori to obtain the desired regularity estimates.  This issue was independently addressed by introducing cutoff Sobolev inequality  and  generalized capacity condition in  \cite{CKW-memo,CKW-jems,CKW-elp} and \cite{GHH18},  respectively. However, verifying these conditions  is challenging without the two-sided pointwise bounds \eqref{e:jumpdensity-pointwise} for the density of the jump kernel, except in cases such as when $(M,d)$ is an ultra-metric space (see \cite{BGHH21}), where the domain of the Dirichlet form includes all characteristic functions on open balls.

In this work, we deal with metric measure spaces with the walk dimension $\beta \ge 2$. The paper is divided into two parts. In the first part, we introduce the  local energy comparability condition \EC \ (see Definition \ref{def:main-condition}(ii) below) and show that it  implies robust versions of  fractional Poincar\'e inequality, Faber-Krahn inequality and  cutoff Sobolev inequality on metric measure spaces, under the  tail estimate of the jump measure outside balls  \TJ \ (see Definition \ref{def:main-condition}(i)).  The condition 
 \EC \  is easier to verify compared to the cutoff Sobolev inequality or the generalized capacity condition,   especially when $\beta s\ge 2$. We  provide a sufficient condition for \EC \  that applies to degenerate jump kernels, extending the work  of  \cite{CS20}. To the best of our knowledge, for $\beta s \ge 2$, this result gives the first general verification of the cutoff Sobolev inequality for nonlocal operators with degenerate jump kernels, except when $(M,d)$ is an ultrametric space.
 
  In the second part, by adapting recent techniques from \cite{CKW-elp, CKW-jems,CKW-memo, GHH18}, we establish robust weak Harnack inequality and a priori Hölder regularity for harmonic functions of nonlocal Dirichlet forms, ensuring that the constants in these estimates remain bounded as $s$ approaches 1. While robust results have been established for nonlocal operators on Euclidean spaces (e.g., \cite{CaSi11, FK13, DK20, GuSc12}), we extend this approach to general metric measure spaces.

\subsection{Main result}
Let $(M,d)$ be a locally compact separable metric space  and  $\mu$ be a  Radon measure on $M$ with full support. The  triple $(M,d,\mu)$ is referred to as a \textit{metric measure space}. For $p \in [1,\infty]$, denote by $L^p(M)$ the $L^p$-space $L^p(M;\mu)$, and by $\lVert \cdot \rVert_p$ the norm in $L^p(M)$. The inner product in $L^2(M)$ will be denoted by $\la \cdot, \cdot \ra$.	 
For $x \in M$ and $r>0$, we use $B(x,r)$ to denote the open ball of radius $r$ centered at $x$, and  $V(x,r)$ for  $\mu(B(x,r))$. Set 
$$	\overline R:=\diam(M)\in (0,\infty].$$ Throughout this paper,  we  assume that all open balls in $M$ are relatively compact and \textit{there exist constants $ \alpha>0$  and $C>1$ such that}
\begin{equation}\tag{VD}\label{e:VD}
	\frac{V(x,R)}{V(x,r)} \le  C \left( \frac{R}{r} \right)^{ \alpha} \quad \text{for all} \;\, x \in M \text{ and } 0<r\le R.
\end{equation}
Condition 	\ref{e:VD} is known as the \textit{volume doubling property}. Note that \ref{e:VD} is equivalent to the following condition:
\begin{equation}\tag{VD*}\label{e:VD2}
	\frac{V(y,R)}{V(x,r)} \le  C \left(\frac{R+ d(x,y)}{r} \right)^{ \alpha} \quad \text{for all} \;\, x,y \in M \text{ and } 0<r\le R.
\end{equation}

Let  	 $\sC: \sD[\sC]\times \sD[\sC] \to \R$ be a non-negative definite symmetric bilinear form where  $\sD[\sC]$ is a dense subspace of $L^2(M)$. For $\lambda>0$, define  $\sC_\lambda(f,g)=\sC(f,g)+\lambda \la f, g\ra$. The symmetric form $(\sC, \sD[\sC])$ is called a \textit{Dirichlet form} on $L^2(M)$, if (i)  $\sD[\sC]$ is a Hilbert space with respect to the inner product $\sC_1$ and  (ii) for any $f \in \sD[\sC]$,  $f_+ \wedge 1 \in \sD[\sC]$ and $\sC(f_+ \wedge 1,f_+ \wedge 1)\le \sC(f,f)$. The Dirichlet form $(\sC,\sD[\sC])$ is called \textit{regular} if $\sD[\sC] \cap C_c(M)$ is dense in $\sD[\sC]$ with $\sC_1$-norm and dense in $C_c(M)$ with the sup-norm, where  $C_c(M)$ denotes the space of all continuous functions on $M$ with compact supports.
For a given regular Dirichlet form $(\sC,\sD[\sC])$, it is known that  each function $f\in \sD[\sC]$ admits a quasi-continuous version $\wt f$ on $M$ (see  \cite[Theorem 2.1.3]{FOT}). In this paper, whenever referring to a regular Dirichlet form $(\sC, \sD[\sC])$, we always consider a quasi-continuous version of $f \in \sD[\sC]$, which we denote as $f$ without further explicit indication. 

Let $(\sC,\sD[\sC])$ be a regular Dirichlet form. The generator $\sL^\sC$ of $(\sC,\sD[\sC])$ is defined by  a self-adjoint  non-positive definite operator in $L^2(M)$ whose $\sD(\sL^C)$ consists of exactly those $f \in \sD[\sC]$ that there is unique $u \in L^2(M)$ so that
$$
\sC(f,g) =- \la u, g\ra \quad \text{for all} \;\,  g \in \sD[\sC],
$$
and $\sL^\sC f=u$.
Let $(P^\sC_t)_{t> 0}$ be the associated semigroup for $\sL^\sC$. Then $(P^\sC_t)_{t\ge 0}$ defines a contraction semigroup on $L^p(M)$ for every $p\in [1,\infty]$.  By  general theory (see \cite{FOT}), there exists a $\mu$-symmetric Hunt process $Z^\sC=\{Z^\sC_t, t \ge 0; \P^x, x \in M \setminus \sN\}$ associated with $(\sC, \sD[\sC])$, where $\sN$ is a properly exceptional set, in the following sense: For any $p \in [1,\infty]$ and any Borel function $f \in L^p(M)$,
$$
P^\sC_tf(x) = \E^x [ f(Z^\sC_t)] \quad \text{for $\mu$-a.e. } x \in M.
$$

A Dirichlet form  $(\sC,\sD[\sC])$  is called \textit{strongly local}, if $\sC(f,g)=0$ for all $f,g \in \sD[\sC]$ such that $f$ has a compact support and $g$ is constant in a neighborhood of supp$[f]$.
The form  $(\sC,\sD[\sC])$ is called \textit{conservative} if  $P^\sC_t \1_M = \1_M$ for all $t>0$.

\medskip

Let $\beta>1$.	We say that $(M,d,\mu)$ satisfies condition Exi$(\beta)$ with the walk dimension $\beta$, if there exists a conservative strongly local Dirichlet form $(\sE^L, \sF^L)$ on $L^2(M)$ 
that has a jointly measurable heat kernel $q$  with the following estimate: There exist constants  $\eta\in (0,1/2)$,  $c_i>0$, $1\le i \le 4$, such that 
for all $t>0$ and $\mu$-a.e. $x,y \in M$,
\begin{align}\label{e:HKE-diffusion}
	\frac{c_1}{V(x, t^{1/\beta})} \1_{\{d(x,y) \le \eta t^{1/\beta} \}}\le q(t,x,y) \le \frac{c_2}{V(x, t^{1/\beta})} \exp \bigg(-c_3 \bigg(\frac{d(x,y)^\beta}{t}\bigg)^{1/(\beta-1)}\bigg).
\end{align}

\begin{example}
	\rm	The Euclidean space $\R^d$  satisfies \text{\rm Exi}$(2)$ with $\sE^L(f,f)=\frac12\int_{\R^d} |\nabla f|^2dx$. 
	If $M$ is a geodesically complete Riemannian manifold and  the
	Ricci curvature of $M$ is non-negative, then \text{\rm Exi}$(2)$ is satisfied. We refer to	 \cite{Gr91, LY86, Sa92}. Moreover, many fractals including the Vicsek set,  the $d$-dimensional Sierpi\'nski gasket and Sierpi\'nski traingle satisfy Exi$(\beta)$ for some $\beta> 2$. See, for example,  \cite{Ba98, Ba03, FHK94, Ki01, Ki09} and references therein.
\end{example}

Many studies have explored the equivalent characterization of \eqref{e:HKE-diffusion} using Harnack inequality and functional inequalities such as  Poicar\'e and  cutoff Sobolev inequalities,  and the generalized capacity condition for strongly local Dirichlet forms. See, for example, \cite{AB15,BBK06,BGK12, GH14, GHL15}.

 We recall the following consequences of condition  Exi$(\beta)$ from \cite[Theorem 4.5 and Lemma 4.7]{KM23}
  and \cite[Theorem 3.1]{BGK12}.

\begin{prop}\label{p:Exi-basic}
	Suppose that \text{\rm Exi}$(\beta)$ holds for $\beta>1$. Then $\beta \ge 2$ and the following reverse volume doubling property holds: There exist constants $C>1$ and $\alpha_0>0$ such that
	\begin{equation}\tag{RVD}\label{e:RVD}
		\frac{V(x,R)}{V(x,r)} \ge 	C^{-1} \left( \frac{R}{r} \right)^{ \alpha_0} \quad \text{for all} \;\, x \in M \text{ and } 0<r\le R<	\overline R.
	\end{equation}
Moreover,		the heat kernel $q(t,x,y)$ can be chosen to be jointly  continuous on $(0,\infty) \times M \times M$ and the  $\mu$-symmetric Hunt process $Z=\{Z_t, t \ge 0; \P^x, x \in M\}$ associated with $(\sE^L,\sF^L)$ can be modified to start from every point in $M$.
\end{prop} 

\medskip

\begin{center}
	\textbf{Condition \text{\rm Exi}$(\beta)$ will be in force throughout this paper.}
\end{center} 

\medskip
This assumption allows us to construct a specific class of nonlocal Dirichlet forms through subordination, which in turn facilitates the establishment of a priori robust function inequalities. 
By Proposition \ref{p:Exi-basic}, we have $\beta \ge 2$ and $(M,d,\mu)$ satisfies \ref{e:RVD} with $\alpha_0 \in (0,\alpha]$, where $\alpha$ is the constant in \ref{e:VD}. Further,  the  process $Z$ associated with $(\sE^L,\sF^L)$ has a jointly continuous heat kernel $q(t,x,y)$ and can start from any point in $M$.

\medskip

We consider a nonlocal bilinear form $(\sE,\sF)$ on $L^2(M)$ defined as
\begin{equation}\label{e:def-sE}
	\begin{split} 
		\sE(f,g)&= \int_{M \times M} (f(x)-f(y))(g(x)-g(y)) J(dx,dy),\\
		\sF&= \sE_1\text{-closure of } \left\{ f \in C_c(M): \sE(f,f)<\infty\right\},
	\end{split}
\end{equation}
where $J(dx,dy)$ is a symmetric Borel measure  on $M \times M$.   Denote by $\sB(M)$ the $\sigma$-algebra of all Borel sets of $M$.  We assume that there exists  $J: M \times \sB(M) \to [0,\infty)$ satisfying the following properties:

\smallskip

$(1)$ For every $x \in M$, $J(x,\cdot)$ is a Borel measure.

$(2)$ For every  $E \in \sB(M)$, the map $x\mapsto J(x,E)$ is a Borel function on $M$.

$(3)$  $J(dx,dy) = J(x,dy)\mu(dx)$ in $M\times M$.

\medskip

To state our results precisely, we introduce several conditions for $(\sE,\sF)$. From now on, we fix constants $R_0 \in (0,	\overline R)$ and $s_0 \in (0,1)$, and denote a constant as `independent of $s\in [s_0,1)$', if it may rely on $s_0$ but remains unaffected by the particular value of $s$.  

\begin{defn}\label{def:main-condition}
	 Let $s\in [s_0,1)$.

	\smallskip
	
	\noindent (i) We say that $\sE$ satisfies the \textit{robust tail estimate condition}  \TJ=rTJ$_\beta(R_0,s_0,s)$ if there exists a constant $\Lambda>0$ independent of $s$ such that for all  $x\in M$ and $r\in (0,R_0)$,
	\begin{align*}
		J(x,B(x,r)^c)\le \frac{\Lambda (1-s)}{r^{\beta s}}.
	\end{align*}

	\noindent	(ii) We say that $\sE$ satisfies the \textit{robust energy comparability condition}  \EC=rEC$_\beta(R_0,s_0,s)$  if there exist constants $C,K_0\ge1$ independent of $s$ such that
	for all $x_0 \in M$,  $r\in (0,R_0)$ and $f\in L^2(B(x_0,r))$,
	\begin{align}\label{e:EC}
	&\int_{B(x_0,r)\times B(x_0,r)} (f(x)-f(y))^2	J(dx,dy) \nn\\
	&\ge C^{-1} (1-s)	 \int_{B(x_0,r/K_0)\times B(x_0,r/K_0)} \frac{(f(x)-f(y))^2}{V(x,d(x,y)) d(x,y)^{\beta s}} \mu(dx)\mu(dy).
\end{align}
\end{defn}

\begin{remark}
\noindent (i)	If $J(dx,dy)$ has a density $J(x,y)$ with respect to $\mu\times \mu$ and  there exists $C \ge 1$ independent of $s$ such that 
	\begin{equation*}
		\frac{C^{-1} (1-s)}{V(x,d(x,y)) d(x,y)^{\beta s}}\le J(x,y) \le \frac{C (1-s)}{V(x,d(x,y)) d(x,y)^{\beta s}} \quad \text{for all} \;\, x,y \in M,
	\end{equation*}
	then  \TJ \ and \EC \  hold.

	\noindent (ii) The opposite direction for \eqref{e:EC} holds  under \TJ. See Proposition \ref{p:EC-upper} below.

\noindent (iii) Condition \EC \ is motivated by condition (A) in \cite{DK20}, where robust estimates in the Euclidean space were established. For general criteria and examples of nonlocal Dirichlet forms on $L^2(\mathbb{R}^d)$ satisfying \EC \ (with $K_0=1$), we refer to \cite{BKS19,  DK20}.  Our novelty   is allowing the constant $K_0$ to be strictly larger than 1. 
 This makes it possible to verify  
\EC \ using a more accessible non-degeneracy condition.  Refer to  \cite[(1.4)]{CS20} for the Euclidean space  case and Theorem \ref{t:coercivity} for its generalization to general metric measure spaces.
\end{remark}

For an open set $U \subset M$ with positive measure and  $f\in L^1(U)$, we let
$$\overline f_U:=\frac{1}{\mu(U)} \int_U f d\mu.$$

\begin{defn}
	Let  $s\in [s_0,1)$. We say that $(\sE,\sF)$ satisfies the \textit{robust Poincar\'e inequality}  \PI=rPI$_\beta(R_0,s_0,s)$,  if there exist constants $C>0$ and  $K_1\ge 1$ independent of $s$ such that  for any  $x_0 \in M$, $0<r<R_0/K_1$ and  $f \in \sF$, 
	\begin{equation*}
		\int_{B(x_0,r)} (f- \overline{f}_{B(x_0,r)})^2 d\mu \le Cr^{\beta s}\int_{B(x_0, K_1r) \times B(x_0,K_1r)} (f(x) - f(y))^2J(dx,dy).
	\end{equation*}
\end{defn}

For a non-empty open set $D \subset M$, 
let $\sF_D$ be the closure of $\sF\cap C_c(D)$ in $\sF$ with respect to the $\sE_1$-norm.  The form $(\sE,\sF_D)$ is called the \textit{part of $(\sE,\sF)$ on $D$}.   It is known  that if  $(\sE,\sF)$ is a regular Dirichlet form on $L^2(M)$, then $(\sE,\sF_D)$ is  a regular Dirichlet form on $L^2(D)$ (see \cite[Theorem 4.4.3]{FOT}).

\begin{defn}
	Let  $s\in [s_0,1)$. We say that $(\sE,\sF)$ satisfies the \textit{robust Faber-Krahn inequality}  \FK=rFK$_\beta(R_0,s_0,s)$, if	there exist constants  $C>0$ and $K_2\ge 1$ independent of $s$ such that for any  $x_0 \in M$, $0<r<R_0/K_2$,   non-empty open set $D\subset B(x_0,r)$ and any $f \in \sF_D$ with $\lVert f\rVert_2=1$,
		\begin{align}\label{e:FK}	
			\left( \frac{V(x_0,r)}{\mu(D)}\right)^{\beta s/ \alpha}	 \le 	  C\bigg(  r^{\beta s}\int_{B(x_0,K_2r)\times B(x_0,K_2r)} (f(x)-f(y))^2J(dx,dy) + 1\bigg).	\end{align}
\end{defn}

Note that we use a local energy term incremented by 1  on the right-hand side of \eqref{e:FK}, differing from the typical use of $r^{\beta s}\sE(f,f)$ in previous works.

For  open subsets $U$ and $V$ of $M$, the notation $U \Subset V$ means that $\overline U \subset V$.

\begin{defn}
	Let $U$ and $V$ be open sets of $M$ with $U \Subset V$ and $\kappa \ge 1$. We say that a  measurable function $\vp$ on $M$ is a \textit{$\kappa$-cutoff function} for $U \Subset V$, if $0\le \varphi \le \kappa$ in $M$, $\varphi \ge 1$ in $U$ and $\varphi=0$ in $V^c$. Any $1$-cutoff function is simply referred to as a \textit{cutoff function}.
\end{defn}

Denote $\wt \sF:=\{ u + a : u \in \sF, \, a \in \R\}$ and $\wt \sF_b:= \wt \sF \cap L^\infty(M)$.   Since $\sE$ has no killing part, we can extend the form $\sE$ from $\sF$ to $\wt \sF$  by letting
$$
\sE(u+a, v+b) :=\sE (u, v), \quad \text{for all}\;\, u,v \in \sF, \; a,b \in \R.
$$

\begin{defn} 
	Let $s\in [s_0,1)$. We say that $(\sE,\sF)$ satisfies the \textit{robust cutoff Sobolev inequality}  \CS=rCS$_\beta(R_0,s_0,s)$, if	for any $\eps\in (0,1)$, there exists  $C=C(\eps)>0$ independent of $s$ such that the following  holds:  For any $x_0 \in M$ and $R, r > 0$ with $R + 2r < R_0$,  there exists a cutoff function $\phi \in \sF$ for $B(x_0,R) \Subset B(x_0,R+r)$ such that  for all $f \in \wt \sF_b$,
	\begin{align}\label{e:rCS}
		&\int_{B(x_0,R+2r) \times M}  f(x)^2(\phi(x)-\phi(y))^2 J(dx,dy) \nn\\
		&\le  \eps  \int_{B(x_0,R+r) \times B(x_0,R+2r)}  \phi(x)^2 (f(x) - f(y))^2 J(dx,dy)   + \frac{C}{r^{\beta s}}\int_{B(x_0,R+2r)} f^2 d\mu.
	\end{align}	
\end{defn}

\bigskip

For $s\in (0,1)$, define  a Sobolev-Slobodeckij space $W^{\beta s/2,2}(M)$ by
\begin{equation}\label{e:def-Sobolev-space}
	W^{\beta s/2,2}(M):=\left\{ u \in L^2(M): \int_{M\times M} \frac{(u(x)-u(y))^2}{V(x,d(x,y))d(x,y)^{\beta s}} \mu(dx)\mu(dy)<\infty \right\}.
\end{equation}

Our first main result is the following theorem. The proof of this theorem will be given at the end of Section \ref{section:robust-function-inequalities}.
\begin{theo}\label{t:main-1}
	Suppose that {\rm \ref{e:VD}} and \text{\rm Exi}$(\beta)$ hold. Let  $R_0  \in (0,	\overline R)$ and $s_0 \in (0,1)$ be given constants, and $(\sE,\sF)$ be  a bilinear form  on $L^2(M)$ given by \eqref{e:def-sE}.
	If $(\sE,\sF)$ satisfies  \text{\rm rTJ}$_\beta(R_0,s_0,s)$ and \text{\rm rEC}$_\beta(R_0,s_0,s)$  for some  $s \in [s_0,1)$, then $(\sE,\sF)$ is a regular Dirichlet form on $L^2(M)$ and $\sF=W^{\beta s/2,2}(M)$. Moreover, \text{\rm rFK}$_\beta(R_0,s_0,s)$, \text{\rm rCS}$_\beta(R_0,s_0,s)$ and   \text{\rm rPI}$_\beta(R_0,s_0,s)$  hold for $(\sE,\sF)$.
\end{theo}

For an open set $D\subset M$, we say that a Borel function $u$ is locally in $\sF_D$, denoted as $u \in \sF^{\rm loc}_D$, if for any relatively compact open set $U\Subset D$, there exists $f \in \sF_D$ such that $u=f$ a.e. on $U$.

Let $\sL$ be the generator of $(\sE,\sF)$.

\begin{defn}\label{def:harmonic}
	\rm Let $D\subset M$ be a non-empty open subset and $f \in L^1(D)$. We say that $u \in \sF^{\rm loc}_D$ is a weak  solution (resp. subsolution, supersolution) of 
	\begin{align}\label{e:Poisson}
		-	\sL  u  = f \;\; \text{ in} \;\, D,
	\end{align}
	if $u$ is locally bounded on $D$ and satisfies the following two properties:

	(1) For any relatively compact open subset $U$ and $V$ of $D$ with $U \Subset V \Subset D$,
	\begin{align}\label{e:harmonicity-1}
		\int_{U \times V^c} |u(y)| \,J(dx,dy)<\infty.
	\end{align}
	
	(2) For any non-negative $\phi \in \sF \cap C_c(D)$, it holds that
	$$
	\sE(u,\phi)=\la f, \phi \ra \quad\;\;(\text{resp. } \sE(u,\phi)\le \la f, \phi \ra , \;\, \sE(u,\phi)\ge \la f, \phi \ra ).
	$$
	
	\medskip

\noindent 	We use the notation ``$-\sL u = f$ in $D$" (resp.\,``$-\sL u \le f$ in $D$", ``$-\sL u \ge f$ in $D$") to denote that $u$ is a weak solution (resp. subsolution, supersolution) to \eqref{e:Poisson}. We say that  $u$ is  $\sE$-\textit{harmonic} in $D$ if $-\sL u = 0$ in $D$.
\end{defn}

For Borel subsets $D_1$ and $D_2$  of $M$ with $\text{dist}(D_1,D_2)>0$ and a Borel function $u$ on $D_2$, we define a nonlocal tail $\sT(u,D_1,D_2)$  of $u$ by 
\begin{align}\label{e:def-tail}
	\sT(u,D_1,D_2):=\sup_{x \in D_1} \int_{D_2} |u(y)| J(x,dy).
\end{align}

\begin{defn}\label{def:WEHI}
	\rm 
	Let  $s\in [s_0,1)$.

	\smallskip

	\noindent (i) We say that \textit{robust weak elliptic Harnack inequality}  \WEHI=rWEHI$_\beta(R_0,s_0,s)$ holds for  $(\sE,\sF)$, if	there exist constants $\delta,C >0$ and $K \ge 1$  independent of $s$ such that for any $x_0 \in M$, $R \in (0,R_0)$, $r \in (0,R/(K+2))$, and any Borel function $u$  that is bounded, non-negative and $-\sL u \ge f$ in $B(x_0,R)$ for  $f \in L^\infty(B(x_0,R))$,  
	\begin{align*}
		&\bigg(\frac{1}{V(x_0,r)}\int_{B(x_0,r)} u^\delta d\mu \bigg)^{1/\delta} \\
		&\le C \left[\, \essinf_{B(x_0,r)} u +  r^{\beta s}\Big(   \NT\left(u_-, B(x_0,Kr), B(x_0,R-2r)^c \right) + \lVert f_-\rVert_{L^\infty(B(x_0,R))}\Big) \,\right].
	\end{align*}

	\noindent (ii) We say that \textit{elliptic H\"older regularity}  \EHR=rEHR$_\beta(R_0,s_0,s)$ holds for  $(\sE,\sF)$, if
	there exist constants $\gamma \in (0,1]$ and $C>0$ independent of $s$ such that for any $x_0 \in M$,  $ R\in (0,R_0)$ and any Borel function $u$ that is bounded in $M$ and $\sE$-harmonic in $B(x_0,R)$,
	\begin{align*}
		|u(x)-u(y)| \le C  \bigg( \frac{d(x,y)}{R}\bigg)^\gamma  \lVert u \rVert_{\infty} \quad \text{for $\mu$-a.e.} \;\, x,y \in B(x_0,R/4).
	\end{align*} 
\end{defn}

Our second main result is the following theorem. Using arguments developed in \cite{CKW-memo, CKW-jems, CKW-elp} and \cite{GHH18}, we obtain regularity results for  $(\sE,\sF)$ under \TJ \ and \EC. Notably, our result  allows the constants to remain independent of  $s$ as $s\to 1$, which was not addressed in the aforementioned works. The proof of this theorem will be presented at the end of Section \ref{section:WEHI-EHR}.
\begin{theo}\label{t:main-2}
Suppose that {\rm \ref{e:VD}} and \text{\rm Exi}$(\beta)$ hold. Let $(\sE,\sF)$ be a  bilinear form on $L^2(M)$ given by \eqref{e:def-sE}. For any fixed $R_0  \in (0,	\overline R)$ and $s_0 \in (0,1)$, the following implications hold:
\begin{align*}
	\text{\TJ \ $+$  \EC} &\  \Rightarrow \ 	\text{\TJ \ $+$  \FK \ $+$  \CS \ $+$ \PI \ $+$ $(\sE,\sF)$ is regular} \\
	&\ \Rightarrow \ 	\text{\TJ \ $+$  \WEHI} \\
	&\ \Rightarrow \ 	\text{\EHR}.
\end{align*}
\end{theo} 

\begin{remark}
	For $s \in (0,1)$, consider the  quadratic form 
	\begin{equation*} 
		\sC^{(s)}(f,g)= (1-s)\int_{\R^d \times \R^d} \frac{(f(x)-f(y))(g(x)-g(y))}{|x-y|^{d+2s}}dxdy, \quad f,g \in W^{s,2}(\R^d),
	\end{equation*}
	where $W^{s,2}(\R^d)$ is the  Sobolev space.  The  form  $(\sC^{(s)},W^{s,2}(\R^d))$ satisfies \TJ \ and \EC. Consequently, by Theorem \ref{t:main-2},   $(\sC^{(s)},W^{s,2}(\R^d))$ satisfies \PI, \WEHI \ and \EHR \ (with $R_0=\infty$). This conclusion reaffirms  the robust fractional  Poincar\'e inequality on $\R^d$, which follows from the results of \cite{BBM02,MS02,MS02-era}, and the robust weak Harnack inequality and elliptic H\"older regularity estimates for the fractional Laplacian, which were established in \cite{DK20}. 
\end{remark}

Lastly, we provide a sufficient condition for \EC, generalizing the main result of \cite{CS20} for the Euclidean case. 
\begin{theo}\label{t:coercivity}
	Let $s \in [s_0,1)$. Suppose that {\rm \ref{e:VD}} holds and
	\begin{align}\label{e:surface}
		\mu\left(\{y\in M: d(x,y)=r\}\right)=0 \quad \text{for all $x\in M$ and $r \in (0,R_0)$}.
	\end{align}
	If	 $J(dx,dy)$ has a density $J(x,y)$  with respect to $\mu\times \mu$ and there exist constants  $\delta_0,\sigma \in (0,1)$ and $\theta>0$ such that  for all $x \in M$, $r\in (0,R_0)$ and $y \in B(x, (1+\delta_0)r)$,
	\begin{align}\label{e:ass-coercivity}
		\mu\left( \left\{ z \in B(y,r) : J(x,z) \ge \frac{\theta(1-s)}{ V(x,d(x,z)) d(x,z)^{\beta s}}  \right\} \right) \ge \sigma V(y,r),
	\end{align}
	then \EC \ holds with  $C,K_0$ depending  on $s_0,\delta_0,\sigma, \theta$ and the constants $C,\alpha$ in {\rm \ref{e:VD}}  only. 
\end{theo}

We say that $(M,d,\mu)$ satisfies \textit{the annular decay property}  if there exist constants $\gamma \in (0,1]$ and $C>1$ such that for all $x \in M$, $r>0$ and $\eps \in (0,1)$,
\begin{align}\label{e:ADP}
	\mu\left( B(x,r) \setminus  B(x, (1-\eps)r) \right) \le C\eps^\gamma V(x,r).
\end{align}
\begin{remark}
	\noindent (i) The result of \cite{CS20} is also robust with respect to the parameter  $s$ since the dependency on $s$ occurs only through the constant  $c_b$	in \cite[Definition 3.7]{CS20}, which has an upper bound in terms of  $s_0$.
	
	\noindent (ii) If $(M,d,\mu)$ satisfies the annular decay property, then \eqref{e:surface} follows directly. Moreover, by the proof of \cite[Lemma 2.1]{CS20},  \eqref{e:ass-coercivity} is equivalent to that \eqref{e:ass-coercivity} holds for $y \in B(x,r)$ only, instead of $y \in B(x,(1+\delta_0)r)$.
	
	\noindent (iii) If $(M,d)$ is a length space, then by \cite[Corollary 2.2]{Bu99}, $(M,d,\mu)$ satisfies the annular decay property.

	\noindent (iv) Suppose that $(M,d,\mu)$ is complete and the  heat kernel $q(t,x,y)$ in Exi$(\beta)$ satisfies the following  lower bound in addition to \eqref{e:HKE-diffusion}: There exist $c_4,c_5>0$ such that  
	\begin{align*}
	 q(t,x,y) \ge \frac{c_4}{V(x, t^{1/\beta})} \exp \bigg(-c_5 \bigg(\frac{d(x,y)^\beta}{t}\bigg)^{1/(\beta-1)}\bigg) \quad \text{for all $t>0$ and $\mu$-a.e. $x,y \in M$.} 
	\end{align*}
	Then, by \cite[Corollary 1.8 and Remark 1.9(a)]{Mu20}, there exists a geodesic metric $\wt d$ such that $d$ and $\wt d$ are bi-Lipschitz equivalent. By replacing  $d$ with $\wt d$, we have that  $(M,\wt d,\mu)$ is a length space and thus satisfies the annular decay property.

	\noindent (v) In our setting, $(M,d,\mu)$ may not satisfy the annular decay property. For instance, let $M:= \{ (x,0): x<-1\} \cup  \{ (x,\sqrt{1-x^2}): -1\le x \le 0\}\cup \{ (0,y): -1\le y\le 1\} \cup  \{ (x,-\sqrt{1-x^2}):   0\le x\le 1\}  \cup \{ (x,0): x\ge 1\} $, $d$ be the metric inherited from $\R^2$, and $\mu$ be the one-dimensional Hausdorff measure.  Since $d$ is comparable with the intrinsic metric on $M$,  \ref{e:VD} holds with $\alpha=1$ and, by \cite[Corollary 3.1 and Theorem 4.1]{LY86}, Exi$(2)$ holds. However, since $\mu(B(0,1+\eps)-B(0,1)) \ge \pi$ for all $\eps>0$, \eqref{e:ADP} fails for any $\gamma \in (0,1]$.
\end{remark}

For the proof of Theorem \ref{t:coercivity},  we mainly follow the strategy of \cite{CS20}. However, significant non-trivial modifications are required since some of their arguments rely on the geodesic and annular decay properties of Euclidean space, which do not apply in our context. The proof of Theorem \ref{t:coercivity} will be provided in Section \ref{section:Thm1}.

\medskip

The paper is organized as follows. Section \ref{section:subordinate} presents robust estimates for subordinators and the Dirichlet forms associated with subordinate processes. In Sections \ref{section:robust-function-inequalities} and \ref{section:mean-value-inequality},  we establish \FK, \CS, \PI \  and the $L^2$-mean value inequality  for $(\sE,\sF)$ under \EC \ and \TJ, by comparing $(\sE,\sF)$ with the one constructed in Section \ref{section:subordinate}.  In Section \ref{section:WEHI-EHR}, we establish the  weak Harnack inequality and H\"older regularity for $(\sE,\sF)$,  thereby concluding the proof of Theorem \ref{t:main-2}. Section \ref{section:Thm1} contains the proof of Theorem \ref{t:coercivity}.

\smallskip

 {\it Notation:} We use $\alpha,\alpha_0$ for the constants in \ref{e:VD} and \ref{e:RVD}, and $\beta,\eta$ for the constants in \eqref{e:HKE-diffusion}.  
We use the  notations $a\wedge b:=\min\{a,b\}$ and  $a\vee b:=\max\{a,b\}$ for $a,b \in \R$. For a subset $D$ of $M$, the set $D^c$ denotes its complement $M\setminus D$. Values of lower case letters with subscripts $c_i$, $i=0,1,2,...$ are fixed in each statement and proof, and the labeling of these constants starts anew in each proof.  
The notation $f(x) \asymp g(x)$ means that there exist constants $c_1,c_2>0$ such that $c_1g(x)\leq f (x)\leq c_2 g(x)$ in the common domain of the definition of $f$ and $g$.

\section{Analysis of subordinate processes}\label{section:subordinate}

A $C^\infty$ function $\phi$ on $(0,\infty)$ is called a \textit{Bernstein function} if $(-1)^{n-1} \phi^{(n)}(\lambda) \ge 0$ for all $n \ge 0$ and $\lambda>0$. It is well known that every Bernstein function $\phi$ can be expressed as
\begin{align*}
	\phi(\lambda) = a + b\lambda + \int_0^\infty (1-e^{-\lambda t}) \Pi(dt),
\end{align*}
where $a,b \ge 0$ are constants and $\Pi$ is a Borel measure  on $(0,\infty)$ satisfying $\int_0^\infty (1 \wedge t) \Pi(dt)<\infty$.
The triplet $(a,b,\Pi)$ is called the \textit{L\'evy triplet} of $\phi$.  A   process $\xi=(\xi_t)_{t \ge 0}$ is called a \textit{subordinator}, if it is a real-valued non-decreasing L\'evy process. For every   subordinator $\xi$, there exists a unique Bernstein function $\phi$ such that 
$$
\E[e^{-\lambda \xi_t}]=e^{-t\phi(\lambda)} \quad \text{for all} \;\, \lambda>0, \, t\ge0.
$$ The function $\phi$ is called the \textit{Laplace exponent} of $\xi$.  We refer to  \cite{SSV12} for fundamental results on Bernstein functions and their connections to subordinators.

For $a \in [0,1]$, we define a constant $m_a$ by
\begin{equation}\label{e:def-kappa}
	m_a:= \left(\frac{1+ \alpha/\beta}{1-2e^{-1}}\right)^{1-a}.
\end{equation}
Note that 
\begin{align}\label{e:kappa-bound}
	1\le 	m_a \le  m_0 \quad \text{for all} \;\, a \in [0,1].
\end{align}
Fix $s_0 \in (0,1)$ and let $s \in [s_0,1)$. Define a  measure $\Pi_{s}$ on $(0,\infty)$ by
\begin{align*}
	\Pi_{s}(dt) :=m_s(1-s)  t^{-1-s}  dt.
\end{align*}
Observe that $\int_0^\infty (1\wedge t) \Pi_{s}(dt)<\infty$. Let  $\phi_{s}$ be the  Bernstein function with  the L\'evy triplet $(0,0, \Pi_{s})$ and  $\xi^{s}$ be a subordinator correponding to the Laplace exponent $\phi_{s}$,  independent of the process $Z$.  Define  a time-changed process $Y^{s}=(Y^s_t)_{t\ge 0}$ by
\begin{equation}\label{e:def-subordinate-process}
	Y^{s}_t:=Z_{\xi^{s}_t}, \quad t\ge 0.
\end{equation}
According to \cite[Theorem 2.1(ii)]{Oku02}, $Y^{s}$ is a  Hunt process associated with  a regular Dirichlet form $(\sE^{s},\sF^{s})$. Moreover, by  \cite[Theorem 2.1(v)]{Oku02}, we have
\begin{align*}
		\sE^{s}(f,g)= \int_{M \times M} (f(x)-f(y))  (g(x)-g(y)) J_{s}(x,y)\mu(dx)\mu(dy), \quad \; f,g \in \sF^s,
\end{align*}
where
\begin{align}\label{e:def-subordinate-jump}
	J_{s}(x,y):=\frac12\int_0^\infty q(t,x,y) \Pi_{s}(dt).
\end{align}
Define a function $q_s$ on $(0,\infty) \times M \times M$ by
\begin{align}\label{e:def-qs}
	q_s(t,x,y) := \int_0^\infty q(a,x,y) \,\P(\xi^{s}_t \in da), \quad t>0, \, x,y \in M.
\end{align} 
From the definition \eqref{e:def-subordinate-process} of $Y^s$, one can deduce  that  $q_s(t,x,y)$ is a heat kernel of $(\sE^{s},\sF^{s})$.

Recall that  the Sobolev-Slobodeckij space    $W^{\beta s/2,2}(M)$ is defined as \eqref{e:def-Sobolev-space}.

\begin{prop}\label{e:subordinate-domain}
	  $\sF^s=W^{\beta s/2,2}(M)$.
\end{prop}
\pf  Let $(Q_t^{s})_{t \ge 0}$ be the  semigroup of $(\sE^s,\sF^s)$. By the general theory, we have 
\begin{align*}
	\sF^s = \Big\{ f \in L^2(M): \sup_{t\in (0,1)} t^{-1}\la f- Q_t^s f, f\ra <\infty\Big\}. 
\end{align*}
See, for example, \cite[Lemma 1.3.4]{FOT}. By the  conservativeness  and the symmetry of $(Q_t^{s})_{t \ge 0}$, we have for all $f \in L^2(M)$ and  $t>0$,
\begin{align}\label{e:Dirichlet-form-qs}
	\frac1t\la f- Q_t^s f, f\ra 	&= \frac{1}{2t}\int_{M\times M} (f(x)-f(y))^2 q_s(t,x,y) \mu(dx)\mu(dy).
\end{align}
 Following the arguments in  the proof of \cite[Theorem 3.1]{BSS03} (where $M$ is assumed to be a unbounded $d$-set), one can deduce that   $q_s(t,x,y)$ satisfies the following estimate: 
\begin{align*}
q_s(t,x,y) \asymp   \frac{1}{V(x, t^{1/(\beta s)})}  \wedge \frac{t}{V(x,d(x,y))d(x,y)^{\beta s}} \quad \text{for $0<t\le 1$ and $x,y\in M$}.
\end{align*}
Combining this with \eqref{e:Dirichlet-form-qs}, we get that for all $f \in L^2(M)$,
\begin{align*}
	\sup_{t \in (0,1)}	\frac1t\la f- Q_t^s f, f\ra  \le \frac{c_1}{2}\int_{M\times M} \frac{(f(x)-f(y))^2}{V(x,d(x,y))d(x,y)^{\beta s}} \mu(dx)\mu(dy)
\end{align*}
and
\begin{align*}
	\sup_{t \in (0,1)}	\frac1t\la f- Q_t^s f, f\ra  &\ge \frac{1}{2c_1}	\sup_{t \in (0,1)}\int_{M\times M: d(x,y) \ge t^{1/(\beta s)}} \frac{(f(x)-f(y))^2}{V(x,d(x,y))d(x,y)^{\beta s}} \mu(dx)\mu(dy)\\
	&=\frac{1}{2c_1}\int_{M\times M} \frac{(f(x)-f(y))^2}{V(x,d(x,y))d(x,y)^{\beta s}} \mu(dx)\mu(dy),
\end{align*}
and the result follows. \qed

\subsection{Analysis of the subordinator $\xi^{s}$}

In this subsection, we establish some  estimates on the distribution of $\xi^{s}$ which do not depend on  $s \in [s_0,1)$.  Define  $H_{s}:(0,\infty) \to (0,\infty)$ by 
\begin{align*}
	H_{s}(\lambda):=\phi_{s}(\lambda) - \lambda \phi_{s}'(\lambda) = \int_0^\infty (1-e^{-\lambda t}- \lambda t e^{-\lambda t}) \Pi_{s}(dt).
\end{align*}
Note that  $\phi_s$ and $H_s$ are increasing continuous functions with $\phi_s(0+)=H_s(0+)=0$ and 
 $\lim_{\lambda \to \infty}\phi_s(\lambda) =\lim_{\lambda \to \infty}H_s(\lambda) =\Pi_s((0,\infty))=\infty$, and  $\phi_s'$ is a decreasing continuous function with $\lim_{\lambda \to \infty}\phi_s'(\lambda) =0$. Let $\phi_s^{-1}$ and $(\phi_s')^{-1}$ be the inverse functions of $\phi_s$ and $\phi_s'$ respectively.  
For each fixed $t>0$,  define a function $g_{s,t}:(0, t\phi_s'(0+)) \to (0,\infty)$ by
\begin{align*}
	g_{s,t}(a):=(\phi_s')^{-1}(a/t).
\end{align*}
By \cite[Lemma 5.2]{JP87}, we get the following left tail probability estimates  for  $\xi^s$.
\begin{prop}\label{p:JP87}
	(i)	For all $t>0$ and $a \in (0, t\phi_s'(0+))$,
	\begin{align*}
		\P( \xi^s_t \le a) \le \exp \big( - t \,(H_{s} \circ g_{s,t})(a) \big).
	\end{align*}
	(ii)	For all $k>0$, $t>0$ and $a \in (0, t\phi_s'(0+))$,
	\begin{align*}
		\P( \xi^s_t \le a) \ge \bigg( 1-\frac{(1+k ) c_0}{k^2 t (H_s \circ g_{s,t})(a)} \bigg) \exp \big( - (1+2k)t (H_s\circ g_{s,t})(a)\big),
	\end{align*}
	where  $c_0:=\sup_{\lambda>0} (-\lambda^2\phi_{s}''(\lambda))/H_{s}(\lambda)$.
\end{prop}

We will  use the following elementary inequality several times in this paper.
\begin{align}\label{e:exp-poly}
	e^{-r} \le (er/a)^{-a} \quad \text{for all} \;\, r,a>0,
\end{align}

\begin{lemma}\label{l:robust-aux-functions}	
  (i) $-\lambda^2 \phi_s''(\lambda) \le 2H_s(\lambda)$ for all $\lambda>0$.

	\noindent (ii) 	$ \phi_s^{-1}(7/t)^{-1}\le 	t \phi_s'(H_s^{-1}(1/t)) \le \phi_s^{-1}(1/t)^{-1}$ for all $t>0$.
\end{lemma}
\pf (i) Using $1-e^{-r} -re^{-r} \ge r^2 e^{-r}/2$ for all $r \ge 1$, we get that for all $\lambda>0$,
\begin{align*}
	&	-\lambda^2\phi_s''(\lambda) =m_s (1-s) \int_0^{\infty}  (\lambda t)^2 e^{-\lambda  t}t^{-1-s}dt \le  2H_s(\lambda).
\end{align*}

\noindent (ii) The result follows from \cite[Lemma 2.4(i)]{CK20} and  the inequality $(e^2-e)/(e-2)<7$.
\qed

\begin{lemma}\label{l:robust-aux-functions-2}	
  (i)  $\phi_s'(\lambda) \ge e^{-1}m_s \lambda^{s-1}$ for all $\lambda>0$.
 
 \noindent (ii) For all $\lambda>0$,
	$$
	H_s(\lambda ) \ge 	\frac{(1-2e^{-1}) m_s (1-s)}{s} \lambda^s.
	$$ 
	
	\noindent (iii) For all $\lambda>0$,
	\begin{align}\label{e:phi-s}
		\frac{m_s}{es}\lambda^s \le 	\phi_{s}(\lambda )  \le 		\frac{m_s}{s}\lambda^s.
	\end{align}
	Consequently,  it holds that
	\begin{align*}
		\bigg(\frac{m_s}{es}\bigg)^{1/s} t^{1/s}  \le 	\phi_s^{-1}(1/t)^{-1} \le \bigg(\frac{m_s}{s}\bigg)^{1/s} t^{1/s} \;\; \text{ for all} \;\, t>0.
	\end{align*}
\end{lemma}
\pf  
\noindent (i) For all $\lambda>0$, we get
\begin{align*}
	\phi_s'(\lambda)
	&\ge m_s (1-s)\int_0^{1/\lambda} t^{-s} e^{-\lambda t} dt  \ge e^{-1}m_s(1-s) \int_0^{1/\lambda} t^{-s}dt = e^{-1}m_s \lambda^{s-1}. 
\end{align*}
(ii) Using the inequality  $1-e^{-r} -re^{-r} \ge 1-2e^{-1}$ for  $r \ge 1$, we get that for all $\lambda>0$,
\begin{align*}
	H_s(\lambda) 	&\ge  (1-2e^{-1}) m_s(1-s)\int_{1/\lambda}^{\infty} t^{-1-s} dt=\frac{(1-2e^{-1}) m_s (1-s)}{s} \lambda^s.
\end{align*}

\noindent (iii)    Since $1-e^{-r} \ge (r\wedge 1)/e$ for all $r>0$, it holds  that for all $\lambda>0$,
\begin{align*}
	\phi_s(\lambda) 	&\ge  \frac{m_s(1-s)}{e} \bigg( \lambda  \int_0^{1/\lambda} t^{-s} dt+  \int_{1/\lambda}^{\infty } t^{-1-s} dt\bigg)= \frac{m_s}{es}\lambda^s.
\end{align*}
Moreover, since $1-e^{-r} \le r \wedge 1$ for all $r>0$, we have 
\begin{equation*}
	\phi_s(\lambda)\le m_s (1-s)  \int_0^{1/\lambda}  \lambda t^{-s} dt + (1-s) \int_{1/\lambda}^{\infty} t^{-1-s}dt= \frac{m_s}{s}\lambda^s.
\end{equation*}
Hence, \eqref{e:phi-s} holds.  Now,  by \eqref{e:phi-s}, we get that  for all $t>0$, 
$$
\frac{m_s}{es}\phi_s^{-1}(1/t)^s\le \frac{1}{t} \le \frac{m_s}{s}\phi_s^{-1}(1/t)^s.
$$ \qed

\begin{lemma}\label{l:subordinator}   (i)  For all $t>0$ and  $a>0$,
	\begin{align*}
		\P( \xi^{s}_t \le a) 
		&\le   \exp \bigg( -\frac{(1+ \alpha/\beta)(1-s)}{s} \bigg(\frac{t^{1/s}}{ea}\bigg)^{s/(1-s)}\bigg).\end{align*}

	\noindent (ii) There exists a  constant $K>1$   independent of $s$ such that 
	\begin{align}\label{e:subordinator-left}
		\P\left(K^{-1}t^{1/s}< \xi^{s}_t\le K t^{1/s}\right)  \ge  e^{-7}/18\quad \text{for all} \;\, t>0.
	\end{align}
\end{lemma}
\pf (i) Let $t,a>0$. By  Lemma \ref{l:robust-aux-functions-2}(i), since
$\phi_s'(g_{s,t}(a))=a/t $, we have 
\begin{align*}
	g_{s,t}(a) \ge (e^{-1}m_s t/a)^{1/(1-s)}.
\end{align*}
Combining this with  Lemma \ref{l:robust-aux-functions-2}(ii), we obtain
\begin{equation}\label{e:H-lambda-estimate}
	(H_s \circ g_{s,t})(a) \ge 	\frac{(1-2e^{-1})m_s(1-s)}{s} \bigg(\frac{m_s t}{ea}\bigg)^{s/(1-s)}.
\end{equation}
Using Proposition \ref{p:JP87}(i)  and \eqref{e:H-lambda-estimate} in the first line below, and the  definition  \eqref{e:def-kappa} of $m_s$  in the second, we deduce that
\begin{align*}
	\P( \xi^{s}_t \le a)& \le \exp \bigg( -  		\frac{(1-2e^{-1})m_s^{1/(1-s)}(1-s)}{s} \bigg( \frac{ t^{1/s}}{ea} \bigg)^{s/(1-s)}  \bigg)\\
	&=  \exp \bigg( -\frac{(1+ \alpha/\beta)(1-s)}{s} \bigg(\frac{t^{1/s}}{ea}\bigg)^{s/(1-s)}\bigg).
\end{align*}
(ii) Let $t>0$ and set $a_0:=t \phi_s'(H_s^{-1}(1/t))$.  By  Proposition \ref{p:JP87}(ii) (with $k=3$) and  Lemma \ref{l:robust-aux-functions}(i),  we have
\begin{equation*}
	\P( \xi^{s}_t \le a_0) \ge \bigg( 1-\frac{8}{9t(H_{s} \circ g_{s,t}) (a_0)} \bigg)  \exp \big( - 7t (H_{s} \circ g_{s,t})(a_0))\big)=  e^{-7}/9.
\end{equation*}Using this, the monotonicity of $\xi^s$, Lemma \ref{l:robust-aux-functions}(ii),  Lemma  \ref{l:robust-aux-functions-2}(iii) and \eqref{e:kappa-bound}, we get
\begin{equation}\label{e:subordinator-left-1}
	e^{-7}/9 \le  \P(\xi^s_t \le \phi_s^{-1}(1/t)^{-1})\le \P\left( \xi^{s}_t \le (m_s/s)^{1/s} t^{1/s}	 \right) \le \P\left( \xi^{s}_t \le (m_0/s_0)^{1/s_0} t^{1/s}	 \right).
\end{equation}
Let $c_1:=\log(18e^7)$ and $a_1:= t \phi_s'(H_s^{-1}(c_1/t))$. By Proposition \ref{p:JP87}(i), 
$$
\P\left(\xi^s_t\le a_1 \right) \le e^{-c_1}=e^{-7}/18.$$ Hence, using the monotonicity of $\xi^s$,  Lemma  \ref{l:robust-aux-functions}(ii),  Lemma \ref{l:robust-aux-functions-2}(iii) and \eqref{e:kappa-bound},  we obtain
\begin{align*}
	& e^{-7}/18\ge 		\P\left(\xi^s_t\le c_1 \phi_s^{-1}(7c_1/t)^{-1} \right) \\
	& \ge  	\P\left(\xi^s_t\le c_1  (m_s/(es))^{1/s} (t/(7c_1))^{1/s} \right) \ge \P\left( \xi^{s}_t \le c_1 (7ec_1)^{-1/s_0}  t^{1/s} \right).
\end{align*}
Combining this with \eqref{e:subordinator-left-1}, we arrive at \eqref{e:subordinator-left} (with $K= (m_0/s_0)^{1/s_0}\vee ((7ec_1)^{1/s_0}/c_1)$). \qed

\subsection{Analysis of the Dirichlet form $(\sE^{s},\sF^{s})$}

Recall that a constant is considered independent of $s\in [s_0,1)$, if it may depend on $s_0$ but remains unaffected by the specific value of $s$. We also recall that $\overline R$ denotes $\diam(M)$.

We begin with the following elementary lemma.

\begin{lemma}\label{l:robust-jumptail} There exists $C>0$  such that for all   $b >0$, $x\in M$ and $r>0$,
	\begin{align*}
		\int_{B(x,r)^c} \frac{1}{V(x,d(x,y)) d(x,y)^{b}}\,\mu(dy) \le   \frac{	C e^b}{b r^b}.
	\end{align*}
\end{lemma}
\pf Using  \ref{e:VD}  and the inequality $1-e^{-b} \ge be^{-b}$, we get
\begin{align*}
	&\int_{ B(x,r)^c} \frac{1}{V(x,d(x,y)) d(x,y)^{b}}\,\mu(dy)\le  \sum_{n=1}^{\infty} \int_{B(x,e^n r)\setminus B(x, e^{n-1}r)}  \frac{	1}{V(x,e^{n-1} r)(e^{n-1}r)^{b}}\,\mu(dy) \\
	&\le \frac{1}{r^b}\sum_{n=1}^{\infty} \frac{V(x,e^{n} r)}{V(x,e^{n-1} r)e^{(n-1)b}}\le \frac{c_1}{r^{b}}\sum_{n=1}^{\infty}  e^{-(n-1)b}  \le \frac{c_1e^b}{br^b}.
\end{align*}
\qed

\begin{prop}\label{p:robust-jumpkernel}
	There exists $C>1$ independent of $s$ such that 
	\begin{align*}
		\frac{	C^{-1}(1-s) }{V(x,d(x,y)) d(x,y)^{\beta s}}\le 	J_{s}(x,y) \le \frac{	C(1-s) }{V(x,d(x,y)) d(x,y)^{\beta s}}  \quad \text{for all} \;\, x,y \in M.
	\end{align*}
\end{prop}
\pf  Let $x,y \in M$. Then $(2d(x,y)/\eta)^\beta \le (2\overline R/\eta)^\beta$.  Using  \eqref{e:def-subordinate-jump} and  \eqref{e:HKE-diffusion}  in the first inequality below, 
\ref{e:VD} in the second, and \ref{e:RVD} in the third,  we get
\begin{align*}
	J_s(x,y) &\ge  c_1(1-s)\int_{(d(x,y)/\eta)^\beta }^{(2d(x,y)/\eta)^\beta } \frac{t^{-1-s} }{V(x,t^{1/\beta})} dt\ge  c_2(1-s)\int_{(d(x,y)/\eta)^\beta }^{(2d(x,y)/\eta)^\beta } \frac{t^{-1-s} }{V(x,\eta t^{1/\beta}/2)} dt\\
	&\ge  \frac{c_3(1-s) }{V(x,d(x,y))}\int_{(d(x,y)/\eta)^\beta }^{(2d(x,y)/\eta)^\beta  } t^{-1-s }  \bigg(\frac{d(x,y)}{\eta t^{1/\beta}/2}\bigg)^{ \alpha_0}  dt\\
	&=  \frac{c_3\eta^{\beta s + \alpha_0}(1- 2^{- \alpha_0-\beta s} ) (1-s)  }{(\eta/2)^{\alpha_0}(s+\alpha_0/\beta)V(x,d(x,y)) d(x,y)^{\beta s}}\ge \frac{2^{\alpha_0}c_3\eta^{\beta}(1- 2^{- \alpha_0} )(1-s)}{ (1+ \alpha_0/\beta)V(x,d(x,y)) d(x,y)^{\beta s}}. 
\end{align*}
On the other hand, by \eqref{e:def-subordinate-jump} and  \eqref{e:HKE-diffusion}, we have
\begin{align*}
	J_s(x,y) & \le  c_4(1-s)\int_0^{d(x,y)^\beta } \frac{t^{-1-s}}{V(x, t^{1/\beta})} \exp \bigg(-c_5 \bigg(\frac{d(x,y)^\beta}{t}\bigg)^{1/(\beta-1)}\bigg)  dt\\
	&\quad +c_4(1-s)\int_{d(x,y)^\beta }^{\overline R^\beta}  \frac{t^{-1-s}}{V(x, t^{1/\beta})}  dt +c_4(1-s)\int_{\overline R^\beta}^\infty  \frac{t^{-1-s}}{V(x, \overline R)}  dt \\
	&=: c_4(1-s)(I_1+I_2+I_3).
\end{align*}
For $I_1$, using  \eqref{e:exp-poly} in the first inequality below and  \ref{e:VD} in the second, we get
\begin{align*}
	I_1 &\le \frac{c_6}{V(x, d(x,y))}\int_0^{d(x,y)^\beta } t^{-1-s} \frac{V(x, d(x,y))}{V(x, t^{1/\beta})} \bigg(\frac{t}{d(x,y)^\beta}\bigg)^{2+ \alpha/\beta}  dt \\
	&\le\frac{c_7}{V(x, d(x,y)) d(x,y)^{2\beta}}\int_0^{d(x,y)^\beta} t^{1-s} dt\\
	&=\frac{c_7 }{(2-s)V(x, d(x,y)) d(x,y)^{\beta s}}\le \frac{c_7}{V(x, d(x,y)) d(x,y)^{\beta s}}.
\end{align*}
For $I_2$, using \ref{e:RVD}, we obtain
\begin{align*}
	I_2&\le \frac{c_8 d(x,y)^{ \alpha_0}}{V(x,d(x,y))}\int_{d(x,y)^\beta }^{\overline R^\beta}   t^{-1-s- \alpha_0/\beta}   dt\le \frac{c_{8}}{( \alpha_0/\beta)V(x, d(x,y)) d(x,y)^{\beta s}}.
\end{align*}
When $\overline R<\infty$, using \ref{e:VD}, we also get
\begin{align*}
	I_3= \frac{1}{s V(x,\overline R) \overline R^{\beta s}} \le  \frac{1}{s_0 V(x,d(x,y)) d(x,y)^{\beta s}} .
\end{align*}
The proof is complete.\qed

\begin{cor}\label{c:robust-jumptail}
	There exists $C>0$ independent of $s$  such that
	\begin{align*}
		\int_{B(x,r)^c} J_{s}(x,y)\mu(dy) \le  \frac{C(1-s)}{r^{\beta s}} \quad \text{for all} \;\, x \in M \text{ and } \,r>0.
	\end{align*}
\end{cor}
\pf Using  Proposition \ref{p:robust-jumpkernel} in the first inequality below and Lemma \ref{l:robust-jumptail} in the second, we get that for all  $x \in M$ and $r>0$,
\begin{align*}
	&	\int_{B(x,r)^c}  J_s(x,y)\mu(dy)  \\
	&\le c_1(1-s)	\int_{ B(x,r)^c} \frac{\mu(dy)}{V(x,d(x,y)) d(x,y)^{\beta s}} \le  \frac{c_2e^{\beta s}(1-s)}{\beta s r^{\beta s}} \le  \frac{c_2e^{\beta }(1-s)}{\beta s_0 r^{\beta s}}.
\end{align*} \qed

Recall that $q_s(t,x,y)$ is the heat kernel of $(\sE^s,\sF^s)$ defined as \eqref{e:def-qs}. We establish a robust near-diagonal estimate for $q_s(t,x,y)$.
\begin{prop}\label{p:NDU}
	There exists $C>0$ independent of $s$  such that 
	\begin{align}\label{e:NDU}
		q_s(t,x,y) \le \frac{C}{V(x,t^{1/(\beta s)})} \quad \text{for all} \;\, t>0\, \text{ and } \, x,y \in M.
	\end{align} 
\end{prop}
\pf For all $t>0$ and $x,y \in M$,  using   \eqref{e:HKE-diffusion} in the first line below, \ref{e:VD} in the third,
Lemma \ref{l:subordinator}(i) in the fourth, and the inequality $e^r \ge r$ for all $r>0$ in the fifth, we obtain
\begin{align*}
	q_s(t,x,y) 	&\le c_1\int_0^{\infty} \frac{\P(\xi^{s}_t \in da)	}{V(x, a^{1/\beta})}\\
	&\le    \frac{ c_1 \P(\xi^{s}_t \ge e^{-1}t^{1/s})}{V(x, (e^{-1}t^{1/s})^{1/\beta})} + c_1 \sum_{n=0}^\infty \frac{  \P(e^{-n-2}t^{1/s}<\xi^{s}_t \le e^{-n-1}t^{1/s})}{V(x, (e^{-n-2}t^{1/s})^{1/\beta})} \nn\\
	&\le \frac{c_2 e^{ \alpha/\beta}}{V(x,t^{1/(\beta s)})} +  \frac{c_2e^{2 \alpha/\beta}}{V(x,t^{1/(\beta s)})}\sum_{n=0}^\infty e^{n \alpha/\beta}   \, \P(\xi^{s}_t \le e^{-n-1}t^{1/s}) \nn\\
	&\le \frac{c_2 e^{ \alpha/\beta}}{V(x,t^{1/(\beta s)})} +  \frac{c_2e^{2 \alpha/\beta}}{V(x,t^{1/(\beta s)})}\sum_{n=0}^\infty  \exp \left( \frac{n \alpha}{\beta} -\frac{(1+ \alpha/\beta)(1-s)}{s} e^{ns/(1-s)}\right)\nn\\
	&\le \frac{c_2 e^{ \alpha/\beta}}{V(x,t^{1/(\beta s)})} +  \frac{c_2e^{2 \alpha/\beta}}{V(x,t^{1/(\beta s)})}\sum_{n=0}^\infty  \exp(-n) = \frac{c_3}{V(x,t^{1/(\beta s)})}.
\end{align*}  \qed

For a non-empty open set $D \subset M$, let $(\sE^{s}, \sF^{s}_D)$ be the part of $(\sE^{s},\sF^{s})$ on $D$ and   $(Q_t^{s,D})_{t \ge 0}$ be its semigroup.

\begin{prop}\label{p:existence-Dirichlet}
	For any bounded non-empty quasi-open set $D\subset M$,  $(\sE^s, \sF^s_D)$ has a  heat kernel $q^D_s(t,x,y)$ defined on $(0,\infty) \times M \times M$.
\end{prop}
\pf Let $x_0\in M$ and $R>0$ be such that $D\subset B(x_0,R)$. For every $t>0$ and $f \in L^1(D)$ with $\lVert f \rVert_1=1$, by the $L^\infty$-contractivity of $(Q^{s,D}_t)_{t\ge 0}$,   Proposition \ref{p:NDU} and \ref{e:VD2}, we have
\begin{align*}
	&	\lVert Q_t^{s,D} f \rVert_\infty \le \lVert Q_{t\wedge 1}^{s,D} f \rVert_\infty \\
	&\le \inf_{x \in D}  \frac{c_1}{V(x, (t\wedge 1)^{1/(\beta s)})}  \le \frac{c_2}{V(x_0,R+ (t\wedge 1)^{1/(\beta s)})} \bigg(\frac{R+(t\wedge 1)^{1/(\beta s)}}{ (t\wedge 1)^{1/(\beta s)}} \bigg)^{ \alpha}<\infty.
\end{align*}
Now,  the result follows from \cite[Theorem 2.1]{GHH21}. \qed

\begin{prop}\label{p:FK}
 $(\sE^s,\sF^s)$ satisfies \FK \  with $R_0=\overline R$ and $K_2=2$.
\end{prop}
\pf Let  $x_0 \in M$, $0<r<\overline R/2$ and  $D\subset B(x_0,r)$ be a non-empty open set.  Set $U:=B(x_0,r)$ and $V:=B(x_0,2r)$. Let $f \in \sF^s_D$ be such that $\lVert f \rVert_2=1$. Since $f=0$ in $D^c$,  we have
\begin{align*}
	\sE^{s}(f,f) &\le 2 \bigg( \int_{D \times V}  +  \int_{D \times V^c} \bigg)  (f(x)-f(y))^2 J_s(x,y)\mu(dx)\mu(dy) =:2(I_1+I_2).
\end{align*}
By Proposition \ref{p:robust-jumpkernel}, since $D\subset V$, we have 
\begin{align}\label{e:FK-1-1}
	I_1&\le c_1(1-s) \int_{V \times V} \frac{(f(x)-f(y))^2}{V(x,d(x,y))d(x,y)^{\beta s}} \mu(dx)\mu(dy) .
\end{align}
Further, using  Corollary \ref{c:robust-jumptail}, since $\lVert f\rVert_2=1$, we get
\begin{align}\label{e:FK-1-2}
	I_2 &\le 2 \int_{D}  f(x)^2 \int_{B(x,r)^c} J_s(x,y)\mu(dy)\mu(dx)\le \frac{c_2(1-s)}{r^{\beta s}}.
\end{align}
On the other hand, by \cite[Lemma 1.3.4(i)]{FOT},  it holds that 
\begin{align}\label{e:FK-2}
	\sE^{s}(f,f) \ge \sup_{t>0}\left[ \frac1t\langle f-Q^{s, D}_t f , f \rangle \right] = \sup_{t>0}\left[ \frac1t \big( 1- \langle Q^{s,D}_t f , f \rangle \big)  \right].
\end{align}
For all $t>0$, using the symmetry of $q_s^D$, the AM-GM inequality, and  Fubini's theorem in the first inequality below, and  Proposition \ref{p:NDU} in the second,  we obtain
\begin{align*}
	&	\langle Q^{s,D}_t f , f \rangle= \int_{D \times D} q_s^D(t,x,y) f(x)f(y) \, \mu(dx)\mu(dy)\\
	&\le \int_D  f(x)^2 \int_{D} q_s^D(t,x,y)\mu(dy)\, \mu(dx)\le c_3\mu(D)\int_D  \frac{f(x)^2}{V(x, t^{1/(\beta s)} )}  \mu(dx).
\end{align*}
Consequently,  using \ref{e:VD2}, we get that for all $0<t\le  r^{\beta s} $,
\begin{align}\label{e:FK-3}
	\langle Q^{s,D}_t f , f \rangle	\le \frac{c_4\mu(D)}{V(x_0,r)}\bigg( \frac{2r }{t^{1/(\beta s)}}\bigg)^{ \alpha} \int_D  f^2 d\mu = \frac{2^{ \alpha}c_4 \mu(D)}{V(x_0,r)}\bigg( \frac{r }{t^{1/(\beta s)}}\bigg)^{ \alpha}.
\end{align}

We assume that,  without loss of generality, the constant $c_4$ in \eqref{e:FK-3} is greater than $1$.   If  $2^{ \alpha+1}c_4\mu(D)\ge V(x_0,r)$,  then by choosing $C$ to be greater than $ (2^{ \alpha+1}c_4)^{\beta/ \alpha}$, we get \eqref{e:FK}. Assume that  $2^{ \alpha+1}c_4\mu(D)< V(x_0,r)$. By taking 
$$t=\bigg( \frac{2^{ \alpha+1}c_4\mu(D)}{V(x_0,r)} \bigg)^{\beta s/  \alpha} r^{\beta s},$$
we get from \eqref{e:FK-2} and \eqref{e:FK-3} that
\begin{align*}
	\sE^{s}(f,f)& \ge \frac{1}{2r^{\beta s}} \bigg( \frac{V(x_0,r)}{2^{ \alpha+1}c_4\mu(D)} \bigg)^{\beta s/  \alpha}   \ge \frac{1}{2^{ 1+( \alpha+1)\beta/ \alpha} c_4^{\beta/ \alpha} r^{\beta s}} \bigg( \frac{V(x_0,r)}{\mu(D)} \bigg)^{\beta s/  \alpha} .
\end{align*}
Combining this with \eqref{e:FK-1-1} and \eqref{e:FK-1-2}, the  result  follows. The proof is complete.  \qed

For a non-empty open set $D \subset M$, let $(\sE^{L},\sF^{L}_D)$ be the part of $(\sE^L,\sF^L)$ on $D$ and  $(Q^{D}_t)_{t \ge 0}$ be its  semigroup. Since  $Z$ has a jointly continuous transition density $q$, and $q(t,x,y)=q(t,y,x)\le c/V(x,t^{1/\beta})<\infty$ for all $t>0$ and $x,y \in M$, from the strong Markov property, we deduce that  $(Q^{D}_t)_{t \ge 0}$ has a heat kernel $q^D$ given by the  Dynkin-Hunt formula:
\begin{align}\label{e:Dynkin-Hunt}
	q^D(t,x,y)=q(t,x,y)-\E^x\left[ q(t-\tau^Z_D, Z_{\tau^Z_D},y): \tau^Z_D<t\right], \quad\; t>0, \, x,y \in D,\end{align}
where  $\tau^Z_D:=\inf\{t>0:Z_t \notin D\}$ denotes the first exit time of $Z$ from $D$.

\begin{prop}\label{p:NDL-diffusion}
	There exist constants $\eps_1, \eps_2 \in (0,1/2)$ and  $C>0$ such that for any $x_0 \in M$ and  $r \in (0, \overline R)$, 
	\begin{align*}
		q^{B(x_0,r)}(t,x,y) \ge \frac{C}{V(x,t^{1/\beta})} \quad \text{for all } \; t \in (0, \eps_1 r^\beta] \text{ and } x,y \in B(x_0, \eps_2 t^{1/\beta}).
	\end{align*}
\end{prop}
\pf  Let   $\eps_1 \in (0,1/2)$ and $\eps_2\in (0,\eta/2)$ be constants determined later, where $\eta \in (0,1)$ is the constant in \eqref{e:HKE-diffusion}. Fix $x_0 \in M$ and $r \in (0,\overline R)$. Write $B:=B(x_0,r)$. By \eqref{e:Dynkin-Hunt} and \eqref{e:HKE-diffusion},  we have for all $t \in (0, \eps_1 r^\beta]$ and  $x,y \in B(x_0, \eps_2 t^{1/\beta})$,
\begin{align}\label{e:NDL-diffusion-1}
	&q^{B}(t,x,y)  \nn\\
	&\ge \frac{c_1}{V(x, t^{1/\beta})} -\E^x\left[ \frac{c_2}{V(Z_{\tau^Z_{B}}, (t-\tau^Z_B)^{1/\beta})} \exp \bigg(-c_3 \bigg(\frac{d(Z_{\tau^Z_{B}},y)^\beta}{t-\tau^Z_B}\bigg)^{1/(\beta-1)}\bigg): \tau^Z_B<t\right].
\end{align}
Observe that 
\begin{equation}\label{e:NDL-diffusion-2}
	d(Z_{\tau^Z_B},x) \wedge	d(Z_{\tau^Z_B},y)  \ge r-  \eps_1^{1/\beta}\eps_2 r>r/2 \quad \text{and} \quad t- \tau^Z_B \le t \le \eps_1 r^\beta.
\end{equation}
In particular, since  $d(x,y) \le 2\eps_2 t^{1/\beta} < r$, we have
\begin{equation}\label{e:NDL-diffusion-3}
	d(Z_{\tau^Z_B},x) \le 	d(Z_{\tau^Z_B},y) + r < 3d(Z_{\tau^Z_B},y).
\end{equation}
Using \eqref{e:exp-poly} in the first inequality below, \eqref{e:NDL-diffusion-3} in the second,   \ref{e:VD2} in the third, and  \eqref{e:NDL-diffusion-2} in the fourth, we get
\begin{align*}
	&	\E^x\left[ \frac{c_2}{V(Z_{\tau^Z_{B}}, (t-\tau^Z_B)^{1/\beta})} \exp \bigg(-c_3 \bigg(\frac{d(Z_{\tau^Z_{B}},y)^\beta}{t-\tau^Z_B}\bigg)^{1/(\beta-1)}\bigg): \tau^Z_B<t\right]\\
	&\le 	\E^x\left[ \frac{c_4}{V(Z_{\tau^Z_{B}}, (t-\tau^Z_B)^{1/\beta})}  \bigg(\frac{(t-\tau^Z_B)^{1/\beta}}{d(Z_{\tau^Z_{B}},y)}\bigg)^{ \alpha+1}: \tau^Z_B<t\right]\\
	&\le 	\E^x\left[ \frac{6^{d_2+1}c_4}{V(Z_{\tau^Z_{B}}, (t-\tau^Z_B)^{1/\beta})}  \bigg(\frac{(t-\tau^Z_B)^{1/\beta}}{2d(Z_{\tau^Z_{B}},x)}\bigg)^{ \alpha+1}: \tau^Z_B<t\right]\\
	&\le 	\E^x\left[ \frac{c_5 (t-\tau^Z_B)^{1/\beta}}{2d(Z_{\tau^Z_{B}},x) \, V(x, 2 d(Z_{\tau^Z_B},x))} : \tau^Z_B<t\right]\le  \frac{c_5 \eps_1^{1/\beta}}{V(x,t^{1/\beta})}.
\end{align*}
Combining this with \eqref{e:NDL-diffusion-1} and taking $\eps_1$ smaller than $(c_1/(2c_5))^\beta$, we get the result. \qed

\begin{prop}\label{p:NDL-subordinate}
	There exist constants $\delta_1,\delta_2\in (0,1/2)$ and $C>0$ independent of $s$   such that for any $x_0 \in M$ and   $r \in (0, \overline R)$,
	\begin{align*}
		q^{B(x_0,r)}_s(t,x,y) \ge  \frac{C}{V(x,t^{1/(\beta s)})} \quad \text{for all } \; t \in (0,  (\delta_1 r)^{\beta s}] \text{ and $\mu$-a.e. } x,y \in B(x_0, \delta_2 t^{1/(\beta s)}).
	\end{align*}
\end{prop}
\pf  Let $x_0 \in M$ and $r\in (0,\overline R)$. Write $B:=B(x_0,r)$. Define for $t>0$ and $x,y \in B$,
\begin{align*}
	r_{s}^B (t,x,y) := \int_0^\infty q^B(a,x,y) \P(\xi^{s}_t \in da).
\end{align*}
By \cite[Proposition 3.1]{SV08},    we have  $q_s^B(t,x,y) \ge 	r_s^B(t,x,y)$ for all $t>0$ and $\mu$-a.e. $x,y\in B$. Thus, it suffices to show that
there exist constants $\delta_1,\delta_2\in (0,1/2)$ and $c_1>0$ independent of $s$, $x_0$ and $r$ such that
\begin{equation*}
	r^B_s(t,x,y) \ge \frac{c_1}{V(x,t^{1/(\beta s)})} \quad\text{ for all  $t \in (0, (\delta_1 r)^{\beta s}] $ and $x,y \in B(x_0,  \delta_2 t^{1/(\beta s)})$} .
\end{equation*}
Let $K>1$ be the constant in Lemma \ref{l:subordinator}(ii), and $\eps_1, \eps_2 \in (0,1/2)$ be the constants in Proposition \ref{p:NDL-diffusion}. Set
$\delta_1:=(\eps_1/K)^{1/\beta}$  and $\delta_2:=\eps_2/K^{1/\beta} .$  
For all $t \in (0, (\delta_1 r)^{\beta s}] $,  we have $	Kt^{1/s}\le   K \delta_1^{\beta} r^\beta = \eps_1r^\beta$ and  $\delta_2t^{1/(\beta s)} = \eps_2 (K^{-1} t^{1/s})^{1/\beta}.$ Using Proposition \ref{p:NDL-diffusion} in the first inequality below, and  \ref{e:VD}  and Lemma \ref{l:subordinator}(ii) in the third, we get that for all $t \in (0, \delta_1r^{\beta s}] $ and $x,y \in B(x_0,  \delta_2 t^{1/(\beta s)})$,
\begin{align*}
	&r_{s}^B(t,x,y) \ge c_2\int_{K^{-1}t^{1/s}}^{ K t^{1/s}} \frac{\P(\xi^{s}_t \in da)}{V(x, a^{1/\beta})} \ge\frac{c_2\,\P( K^{-1}t^{1/s} \le  \xi^{s}_t \le  Kt^{1/s} ) }{V(x,K^{1/\beta} t^{1/(\beta s)})}\ge \frac{c_3}{V(x,t^{1/(\beta s)})}.
\end{align*} \qed

\begin{prop}\label{p:PI}
 $(\sE^s,\sF^s)$ satisfies \PI \  with $R_0=\overline R$.
\end{prop}
\pf  We follow the proof of \cite[Proposition 3.5(i)]{CKW-jems}, which was motivated by the argument in \cite[Theorem 5.1]{KS87}. Let $\delta_1,\delta_2\in (0,1/2)$ be the constants in Proposition \ref{p:NDL-subordinate} and set $K_1:=1/(\delta_1\delta_2)$. Let $x_0 \in M$ and $0<r<\overline R/K_1$. Write $B:=B(x_0,r)$ and $V:=B(x_0,K_1r)$.  

Consider a bilinear form
\begin{align*}
	\sC^{s,V}(f,f)&:= \int_{V \times V} (f(y) - f(x)) (g(y) - g(x)) J_s(x,y)  \mu(dx)\mu(dy),\\
	\sD[\sC^{s,V}]&:=\big\{ f \in L^2(V): \sC^{s,V} (f,f)<\infty\big\}.
\end{align*}
By Proposition \ref{p:robust-jumpkernel}(i), we have
\begin{align}\label{e:PI-con-1}
	\sC^{s,V}(f,f)  \le 	c_1(1-s)	\int_{V \times V} \frac{(f(x) - f(y))^2}{V(x,d(x,y))d(x,y)^{\beta s}} \mu(dx)\mu(dy) \quad \text{for all} \;\, f \in \sF^s.
\end{align}
Hence, $\sF^s|_V \subset \sD[\sC^{s,V}]$.
Further, by Fatou's lemma, we see that $(\sC^{s,V}, \sD[\sC^{s,V}])$ is closable and is a Dirichlet form on $L^2(V)$. Let $(\overline Q^{s,V}_t)_{t> 0}$ be the  semigroup associated with $(\sC^{s,V}, \sD[\sC^{s,V}])$. Since  $\sC^{s,V}(\1_V,\1_V)=0$ and $\1_V \in \sD[\sC^{s,V}]$, the semigroup $(\overline Q^{s,V}_t)_{t> 0}$ is conservative.  Let $\sD[\sC^{s,V}]_V$ be the closure of $\sD[\sC^{s,V}] \cap C_c(V)$ in $L^2(V)$ and let  $(\overline R^{s,V}_t)_{t> 0}$ be the  semigroup associated with $(\sC^{s,V}, \sD[\sC^{s,V}]_V)$.  According to \cite[Theorem 5.2.17]{CF12}, $(\sC^{s,V}, \sD[\sC^{s,V}]_V)$ is the resurrected Dirichlet form of $(\sE^s, \sF^s_V)$. 
Since $(\sC^{s,V}, \sD[\sC^{s,V}]_V)$ is a part of  $(\sC^{s,V}, \sD[\sC^{s,V}])$, it follows that 
\begin{align}\label{e:PI-1}
	\overline Q^{s,V}_tf \ge	\overline R^{s,V}_tf \ge Q^{s, V}_t f \quad \text{for all $t>0$ and  $0\le f \in L^2(V)$}.
\end{align}

Set $t_0:=(\delta_1K_1r)^{\beta s}=(r/\delta_2)^{\beta s}$. For all $f \in \sF^{s}$,  using \cite[Lemma 1.3.4(i)]{FOT} in the first line below, the   conservativeness and the symmetry   of  $(\overline Q^{s,V}_t)_{t \ge 0}$ in the second,  \eqref{e:PI-1} in the third,  Proposition \ref{p:NDL-subordinate} in the fourth and \ref{e:VD} in the last, we obtain
\begin{align}\label{e:PI-con-2}
	\sC^{s,V}(f,f)&\ge \frac{1}{t_0} \int_{V} f(x) (f(x) - \overline Q^{s,V}_{t_0}f(x)) \,\mu(dx)\nn\\
	& =  \frac{1}{t_0} \int_{V} \overline  Q^{s,V}_{t_0}  \left( \frac12f(z)^2 -f(z) f(\cdot)  + \frac12 f(\cdot)^2\right) (x)\big|_{z=x} \,\mu(dx)  \nn\\
	&\ge \frac{1}{2t_0} \int_{V}  Q^{s,V}_{t_0} \left( f(z)-f(\cdot)\right)^2 (x)\big|_{z=x}\mu(dx)\nn\\
	&\ge \frac{c_2}{2t_0} \int_{B} \int_B \frac{(f(x)-f(y))^2}{V(x,t_0^{1/(\beta s)})} \,\mu(dy)\mu(dx)\nn\\
	&\ge \frac{c_3\delta_2^{\beta }}{2 r^{\beta s}} \int_{B} \int_B \frac{(f(x)-f(y))^2}{V(x,r)} \,\mu(dy)\mu(dx) .
\end{align}
On the other hand,  using \ref{e:VD2} in the third inequality below, we see that
\begin{align*}
	&\int_{B} \int_B \frac{(f(x)-f(y))^2}{V(x,r)} \,\mu(dy)\mu(dx) \ge \inf_{a \in \R}\int_{B}  \frac{(f(x)-a)^2}{V(x,r)} \mu(dx)\, \int_B \mu(dy)\\
	& \ge \inf_{x \in B} \frac{V(x_0,r)}{V(x,r)}\inf_{a \in \R}\int_{B}(f(x)-a)^2  \mu(dx) \ge  c_4\inf_{a \in \R}\int_{B}(f(x)-a)^2  \mu(dx)  =c_4 \int_{B} (f- \overline{f}_{B})^2 d\mu.
\end{align*}
Combining this with \eqref{e:PI-con-1} and  \eqref{e:PI-con-2},   we get the desired result.  \qed

For a non-empty open set $D \subset M$, denote by $(G^{s, D}_{\lambda})_{\lambda>0}$  the resolvent corresponding to the semigroup $(Q^{s,D}_t)_{t>0}$ defined by
$$	G^{s,D}_\lambda f:= \int_0^\infty e^{-\lambda t} Q^{s, D}_tf dt \quad \text{for $\lambda>0$ and $f \in L^2(D)$.}$$ 
By \cite[Theorem 4.4.1(i)]{FOT}, for any $\lambda>0$ and $ f \in L^2(D)$, we have $G^{s, D}_\lambda f \in \sF^{s}_D$ and 
\begin{align}\label{e:resolvent}
	\sE^{s}_\lambda(G_\lambda^{s,D}f, v ) = \la f, v \ra  \quad \text{for all} \;\, v \in \sF^{s}_{D}.
\end{align}

\begin{lemma}\label{l:resolvent}
	There exists a constant $\kappa \ge 1$ independent of $s$ such that the following  holds:	Let $x_0 \in M$, $R>0$ and $r\in (0,\overline R)$. For   $U_0:=B(x_0,R)$, $U_1:=B(x_0,R+r)$ and $\lambda:=r^{-\beta s}$, we have
	\begin{align}\label{e:GU-1}
		\lambda  G_\lambda^{s,U_1} \1_{U_1} \le 1 \quad \text{$\mu$-a.e. on $U_1$} \quad \text{and} \quad 
		\lambda	G_\lambda^{s,U_1} \1_{U_1} \ge 1/\kappa \quad \text{$\mu$-a.e. on $U_0$}.
	\end{align}
\end{lemma}
\pf For any $0 \le f \in L^2(U_1)$, we have
\begin{align*}
	&	\la \lambda	G_\lambda^{s,U_1} \1_{U_1}, f \ra=  \int_0^{\infty} \lambda  e^{-\lambda t} \la Q^{s,U_1}_t \1_{U_1}, f \ra \le \lVert f \rVert_{L^1(U_1)} \int_0^{\infty}\lambda e^{-\lambda t}  dt = \lVert f \rVert_{L^1(U_1)}.
\end{align*}
Hence, $	\lambda  G_\lambda^{s,U_1} \1_{U_1} \le 1 $ $\mu$-a.e. on $U_1$.   On the other hand, by  Proposition \ref{p:NDL-subordinate} and \ref{e:VD}, for  any $t \in (0, (\delta_1r)^{\beta s}]$ and  a.e. $x \in U_0$,
\begin{align*} 
	Q^{s, U_1}_t \1_{U_1} (x) \ge 	Q^{s, B(x,r)}_t \1_{B(x,r)} (x)\ge  \frac{c_1}{V(x,t^{1/(\beta s)})} \int_{B(x, \delta_2 t^{1/(\beta s)})} \mu(dy) \ge c_2.
\end{align*}
Using this, we get that    for any $0 \le f \in L^2(U_0)$,
\begin{align*}
	&\la \lambda	G_\lambda^{s,U_1} \1_{U_1}, f \ra \ge  \int_0^{\delta_1^{\beta s}/\lambda} \lambda  e^{-\lambda t} \int_{U_0}  f(x) Q_t^{s,U_1} \1_{U_1}(x) \mu(dx) \, dt\nn\\
	&\ge  c_2\lVert f \rVert_{L^1(U_0)} \int_0^{\delta_1^{\beta s}/\lambda} \lambda  e^{-\lambda t} dt = c_2 (1- e^{-\delta_1^{\beta s}}) \lVert f \rVert_{L^1(U_0)} \ge c_2 (1- e^{-\delta_1^{\beta}})\lVert f \rVert_{L^1(U_0)}.
\end{align*}
Therefore, the second inequality in \eqref{e:GU-1} holds with $\kappa=1/ (c_2(1-e^{-\delta_1^{\beta}}))$. \qed

Denote $\wt \sF^{s}:=\{ f + a : f \in \sF^{s}, \, a \in \R\}$ and $\wt \sF^{s}_b:= \wt \sF^{s} \cap L^\infty(M)$. Since $\sE^{s}$ has no killing part, the bilinear form $\sE^{s}$ can be extended to functions from  $\wt \sF^{s}$ by letting
$$
\sE^{s}(f+a, g+b) :=\sE^{s} (f, g) \quad \text{for all}\;\, f,g \in \sF^{s}, \; a,b \in \R.
$$

We now establish a robust generalized capacity condition for $(\sE^s,\sF^s)$. The proof of the following result is originally due to  \cite[Lemma 5.4]{AB15} and \cite[Lemma 2.8]{GHH18}.

\begin{lemma}\label{l:Gcap-robust}
	For any $x_0 \in M$, $R>0$ and $r\in (0,\overline R)$, there exists a $\kappa$-cutoff function $\vp \in \sF^{s}$ for $B(x_0,R) \Subset B(x_0,R+r)$ such that 
	\begin{align}\label{e:Gcap-robust}
		\sE^{s}(f^2 \vp, \vp) \le   \frac{\kappa^2 }{r^{\beta s}} \int_{B(x_0,R+r)} f^2 d\mu \quad \text{for any} \;\, f \in \wt\sF^{s}_b,
	\end{align}	
	where $\kappa \ge 1$ is the constant in Lemma \ref{l:resolvent}.
\end{lemma}
\pf 
Let $x_0 \in M$, $R>0$ and $r\in (0,\overline R)$. Set  $\lambda:=r^{-\beta s}$ and  $\vp :=\kappa\lambda G_\lambda^{s,B(x_0,R+r)} \1_{B(x_0,R+r)}.$
Then $\vp \in  \sF^{s}_{B(x_0,R+r)}$. Further,  by Lemma \ref{l:resolvent},  $\vp$ is a $\kappa$-cutoff function for $B(x_0,R) \Subset B(x_0,R+r)$. According to  \cite[Proposition 15.1]{GHH23}, 
whenever $f \in \wt\sF_b^{s}$,  we have $f^2 \vp \in \sF_{B(x_0,R+r)}^{s}$. Hence, using \eqref{e:resolvent}, we get that   for any $f \in \wt\sF_b^{s}$,
\begin{align*}
	&\sE^{s}(f^2\vp, \vp)\le \sE^{s}_\lambda(f^2 \vp, \vp) = \kappa \lambda \la f^2 \vp, \1_{U_1} \ra \le \kappa^2 \lambda \int_{U_1} f^2 d\mu.
\end{align*} \qed

Recall  the following property of nonlocal Dirichlet forms from \cite[Lemma 3.5]{CKW-memo}. 
\begin{lemma}\label{l:l.3.5}	For any Borel subset $D \subset M$, constant $k>1$, and any  $f,g \in \wt \sF^{s}_b$,	\begin{align}\label{e:l.3.5}	&	(1-k^{-1}) \int_{D\times D} f(x)^2 (g(x)-g(y))^2J_{s}(x,y) \mu(dx)\mu(dy)\nn\\	&\le \int_{D \times D} (f(x)^2 g(x) - f(y)^2 g(y))(g(x)-g(y))J_{s}(x,y)  \mu(dx)\mu(dy) \nn\\	&\quad + k \int_{D\times D} g(x)^2 (f(x)-f(y))^2 J_{s}(x,y)  \mu(dx)\mu(dy),	\end{align}	provided that all three integrals in \eqref{e:l.3.5} are absolutely integrable.\end{lemma}

Following the proof of \cite[Lemma 2.4]{GHH18}, we deduce the next result from  Lemma \ref{l:Gcap-robust}.   We give a detailed proof for completeness.

\begin{prop}\label{p:CS-0}
 $(\sE^s,\sF^s)$ satisfies \CS \  with $R_0=\overline R$.
\end{prop}
\pf  Let $x_0 \in M$, $R>0$ and $r \in (0,\overline R)$.  Define $U_a=B(x_0,R + a r)$ for $a \ge 0$. Let $\kappa \ge 1$ be the constant in Lemma \ref{l:resolvent}. By Lemma \ref{l:Gcap-robust}, there exists a $\kappa$-cutoff function $\vp \in \sF^{s}$ for $U_0 \Subset U_{1}$  such that 
\begin{align}\label{e:Gcap-robust-2}
	\sE^{s}(f^2 \vp, \vp) \le   \frac{\kappa^2 }{r^{\beta s}} \int_{U_1} f^2 d\mu \quad \text{for all} \;\, f \in \wt\sF^{s}_b.
\end{align}	 
Set 
$\phi:=\vp \wedge 1.$ Then $\phi\in \sF^{s}$ and $\phi$ is a cutoff function for $U_0\Subset  U_1$.  Moreover, using Proposition \ref{p:robust-jumpkernel} and  the fact $|\phi(x)-\phi(y)|\le |\vp(x)-\vp(y)|$ for  $x,y \in M$ in the first inequality below, and  Lemma \ref{l:l.3.5} (with $k=2$) in the second, we get that  for all $f \in \wt \sF^{s}_b$,
\begin{align*}
	&(1-s)\int_{U_2 \times U_2}   \frac{f(x)^2(\phi(x)-\phi(y))^2}{V(x,d(x,y))d(x,y)^{\beta s}} \,\mu(dx)\mu(dy)\nn\\
	& \le  c_1\int_{U_{2 } \times U_{2 }} f(x)^2(\vp(x)-\vp(y))^2 J_s(x,y)  \,\mu(dx)\mu(dy) \nn\\
	&\le 2c_1 \int_{U_{2 } \times U_{2 }} (f(x)^2 \vp(x)-f(y)^2 \vp(y))(\vp(x)-\vp(y)) J_s(x,y) \, \mu(dx)\mu(dy)\nn\\
	&\quad + 4c_1\int_{U_{2 } \times U_{2}} \vp(x)^2 (f(x)-f(y))^2 J_s(x,y) \,\mu(dx)\mu(dy)\nn\\
	&=:I_1+I_2.
\end{align*}
For all $x,y \in (M \times M) \setminus (U_{2}\times U_{2 })$, since either $\vp(x)=0$ or $\vp(y)=0$, we have  $(f(x)^2 \vp(x)-f(y)^2 \vp(y))(\vp(x)-\vp(y))  \ge 0.$ Hence, using  \eqref{e:Gcap-robust-2}, we get
\begin{align}\label{e:CS-0-I1}
	I_1 &\le 2c_1\sE^{s}(f^2\vp, \vp) \le  \frac{2c_1\kappa^2 }{r^{\beta s}} \int_{U_1} f^2 d\mu.
\end{align}
Moreover, by Proposition \ref{p:robust-jumpkernel}, since $\vp^2 \le \kappa^2 \phi^2$ and $\vp=0$ in $U_1^c$, we have
\begin{align}\label{e:CS-0-I2}
	I_2 &\le c_2(1-s)\int_{U_1 \times U_2} \frac{\vp(x)^2 (f(x) - f(y))^2}{ V(x,d(x,y))d(x,y)^{\beta s}}  \,\mu(dx)\mu(dy)\nn\\
	& \le \kappa^2 c_2(1-s)\int_{U_1 \times U_2} \frac{\phi(x)^2 (f(x) - f(y))^2}{ V(x,d(x,y))d(x,y)^{\beta s}}  \,\mu(dx)\mu(dy).
\end{align}
On the other hand, since $\phi$ is a cutoff function for $U_0 \Subset U_1$, by Lemma \ref{l:robust-jumptail}, it holds that 
\begin{align}\label{e:CS-0-I3}
	&(1-s)\int_{U_2 \times U_2^c}   \frac{f(x)^2(\phi(x)-\phi(y))^2}{V(x,d(x,y))d(x,y)^{\beta s}} \,\mu(dx)\mu(dy)\nn\\
	&\le(1-s)\int_{U_1}  f(x)^2 \int_{U_2^c}\frac{\mu(dy)}{V(x,d(x,y))d(x,y)^{\beta s}} \,\mu(dx)\nn\\
	&\le(1-s)\int_{U_1}  f(x)^2 \int_{B(x,r)^c}\frac{\mu(dy)}{V(x,d(x,y))d(x,y)^{\beta s}} \,\mu(dx)\nn\\
	&\le \frac{c_3e^{\beta s}(1-s)}{(\beta s)r^{\beta s}} \int_{U_1} f^2 d\mu\le \frac{c_3e^{\beta }(1-s)}{(\beta s_0)r^{\beta s}} \int_{U_1} f^2 d\mu.
\end{align}
Combining \eqref{e:CS-0-I1}, \eqref{e:CS-0-I2} and \eqref{e:CS-0-I3}, we arrive at the desired result.\qed

\section{Robust function inequalities under  \TJ \ and \EC}\label{section:robust-function-inequalities}

We begin with a standard covering lemma. 
\begin{lemma}\label{l:VD-covering}
For any $k\ge 1$,	there exists $N_0=N_0(k) \in \N$ such that for each $r>0$, there exists an open covering  $\{B(z_i,r)\}_{i=1}^\infty$ of $M$ satisfying
	\begin{align}\label{e:VD-covering}
		\sum_{i=1}^\infty \1_{B(z_i,kr)}  \le N_0\;\; \text{on} \;\,  M.
	\end{align}
\end{lemma}
\pf By the Vitali covering lemma,  there exists a collection $\{B(z_i,r/5)\}_{i=1}^\infty$ of pairwise disjoint open balls in $M$ such that $M =\cup_{i=1}^\infty B(z_i,r)$. Suppose that $y \in M$ is in $N$ of the balls $B(z_i,kr)$ ($N$ may be infinite). Then $B(y,(k+1)r)$ contains at least $N$ balls $B(z_i,r)$. By \ref{e:VD2}, when $y \in B(z_i,kr)$, we have $
V(z_i,r) \ge c_1V(y,(k+1)r)$. It follows that
\begin{align*}
	V(y,(k+1)r) \ge \sum_{i: y \in B(z_i,kr)} V(z_i,r) \ge c_1 N V(y,(k+1)r).
\end{align*}
This leads to the conclusion that $N\le 1/c_1$, establishing \eqref{e:VD-covering}. \qed

\begin{prop}\label{p:EC-upper}
Suppose that \TJ \ holds.	There exists $C\ge1$ independent of $s$ such that
	for all $x_0 \in M$,  $r\in (0,R_0)$ and $f\in L^2(B(x_0,r))$,
	\begin{align*}
		&\int_{B(x_0,r)\times B(x_0,r)} (f(x)-f(y))^2	J(dx,dy) \nn\\
		&\le C (1-s)	 \int_{B(x_0,r)\times B(x_0,r)} \frac{(f(x)-f(y))^2}{V(x,d(x,y)) d(x,y)^{\beta s}} \mu(dx)\mu(dy).
	\end{align*}
\end{prop}
\pf Let $B=B(x_0,r)$ be an open ball in $M$ with $r\in (0, R_0)$ and $f:B\to \R$.  For all $x,y \in B$ and $z \in B(y,d(x,y)/2)$, we have
\begin{align}\label{e:singular-region}
	\frac23 d(x,z)\le  d(x,y)\le 2d(x,z) \quad \text{and} \quad d(x,y) \ge 2d(y,z).
\end{align} 
Using \eqref{e:singular-region}, \ref{e:VD2}, Fubini's theorem and the symmetry of $J$ in the second inequality below,  and \TJ \  in the third,  we obtain
\begin{align*}
	&	\int_{B\times B} (f(x)-f(y))^2 J(dx,dy) \\
	&\le 2\int_{B\times B}  \int_{B(y,d(x,y)/2)} \frac{((f(x)-f(z))^2+(f(y)-f(z))^2)}{ V(y,d(x,y)/2)}\mu(dz) J(x,dy)\mu(dx) \\
	& \le c_1\int_{B\times B}   \frac{(f(x)-f(z))^2}{ V(x,d(x,z))} \int_{B \setminus B(x, 2d(x,z)/3)} J(x,dy)\mu(dx) \mu(dz) \\
	&\quad + 2\int_{B\times B}   \frac{(f(y)-f(z))^2}{ V(y,d(y,z))} \int_{B \setminus B(y, 2d(y,z))} J(y,dx)\mu(dy) \mu(dz)\\
	& \le c_2\Lambda(1-s)\int_{B\times B} \frac{(f(x)-f(z))^2}{V(x,d(x,z))d(x,z)^{\beta s}} \mu(dx)\mu(dz).
\end{align*} \qed 

\begin{prop}\label{p:domain=domain}
	Suppose that \TJ \ and \EC \ hold. There exists $C>1$ independent of $s$ such that
	\begin{align*}
		C^{-1}\sE^{s}_1(f,f) \le \sE_1(f,f) \le C\sE^{s}_1(f,f) \quad \text{for all} \;\, f \in L^2(M).
	\end{align*}
\end{prop}
\pf  Let $K_0 \ge 1$ be the constant in \EC. Set $r:= 1 \wedge (R_0/(3K_0))$. By Lemma \ref{l:VD-covering}, there exist  $N_0=N_0(K_0)\in \N$ and  an open covering $\{B(z_i,r)\}_{i=1}^\infty$ of $M$ such that
	\begin{align}\label{e:VD-covering-appl}
	\sum_{i=1}^\infty \1_{B(z_i,2K_0r)}  \le N_0\;\; \text{on} \;\,  M.
\end{align}
Let $f \in L^2(M)$. Since   $\{B(z_i,r)\}_{i=1}^\infty$  is an open covering, we have
\begin{align}\label{e:domain=domain-1}
	&\sE^{s}(f,f)\le \sum_{i=1}^\infty \int_{B(z_i,r)\times M}  (f(x)-f(y))^2 J_s(x,y)  \mu(dy)\mu(dx).
\end{align}
By Proposition \ref{p:robust-jumpkernel}(i), \EC \ and \eqref{e:VD-covering-appl}, we see that
\begin{align}\label{e:domain=domain-first}
	&\sum_{i=1}^\infty \int_{B(z_i,r)\times B(z_i,2r)}(f(x)-f(y))^2 J_s(x,y)  \mu(dy)\mu(dx)\nn\\
	&\le c_1 \sum_{i=1}^\infty  \int_{B(z_i,2K_0r) \times B(z_i,2K_0r)} (f(x)-f(y))^2 J(dx,dy) \le c_2N_0  \sE(f,f).
\end{align}
On the other hand, using the Cauchy inequality and  Fubini's theorem in the first inequality below, \eqref{e:VD-covering-appl} and the symmetry of $J_s$ in the second, and  Corollary \ref{c:robust-jumptail} in the third, we get
\begin{align}\label{e:domain=domain-second}
	&\sum_{i=1}^\infty \int_{B(z_i,r)\times B(z_i,2r)^c} (f(x)-f(y))^2 J_s(x,y)  \mu(dy)\mu(dx)\nn\\
	&\le 2\sum_{i=1}^\infty\int_{B(z_i,r)} f(x)^2\int_{B(x,r)^c} J_s(x,y)  \mu(dy)\mu(dx)\nn\\
	&\quad + 2\int_{M} f(y)^2 \sum_{i=1}^\infty\int_{B(z_i,r) \setminus B(y,r)} J_s(y,x)  \mu(dx)\mu(dy)\nn\\
	&\le 4N_0\int_{M} f(x)^2\mu(dx) \sup_{z\in M} \int_{ B(z,r)^c} J_s(z,y)  \mu(dy) \nn\\
	&\le c_2N_0(1-s)r^{-\beta s} \lVert f \rVert_2^2 \le c_2N_0(1-s)r^{-\beta} \lVert f \rVert_2^2.
\end{align}
Combining  \eqref{e:domain=domain-1}, \eqref{e:domain=domain-first} and \eqref{e:domain=domain-second}, we arrive at $\sE^s(f,f) \le c_3\sE_1(f,f)$. Similarly, using  \TJ \ and Proposition \ref{p:EC-upper} in place of Corollary \ref{c:robust-jumptail} and \EC \ respectively and repeating the preceding argument, we can establish $\sE(f,f) \le c_4\sE^s_1(f,f)$. \qed

\begin{cor}\label{c:domain=domain}
	Suppose that  \TJ \ and \EC \ hold. Then $ \sF_D = \sF_D^{s}$ for any open set $D\subset M$. In particular, $\sF=\sF^s=W^{\beta s/2,2}$ and  $(\sE,\sF)$ is a regular Dirichlet form on $L^2(M)$.
\end{cor}
\pf  Let $f \in \sF^{s}_D$. Since $(\sE^{s}, \sF^{s}_D)$ is regular, there exists a $\sE^{s}_1$-Cauchy sequence $(f_n)_{n\ge 1}$ in $ C_c(D) \cap \sF^{s}$  converging to $f$ in  $\sE^{s}_1$-norm. By Proposition \ref{p:domain=domain},  $(f_n)_{n \ge 1}$ is a $\sE_1$-Cauchy sequence in $C_c(D)$ that converges to $f$ in $\sE_1$-norm. Hence, we get $f \in \sF_D$, implying  that $\sF^{s}_D \subset \sF_D$. Similarly, we can deduce that $\sF_D \subset \sF^{s}_D$. Hence, $\sF_D=\sF^s_D$. Now by letting $D=M$, we obtain $\sF=\sF^s=W^{\beta s/2,2}$. Since $\sF$ is non-empty,  $(\sE,\sF)$ is a regular Dirichlet form on $L^2(M)$.
\qed

By Propositions \ref{p:FK} and \ref{p:PI}, and Corollary \ref{c:domain=domain}, we get the next proposition. 
\begin{prop}\label{p:FK-PI}
	If \TJ \ and \EC \ hold, then $(\sE,\sF)$ satisfies \FK \ and \PI.
\end{prop}

In the remainder of this section, we establish \CS \ for $(\sE,\sF)$ under \TJ \ and \EC. The verification of \CS \ is challenging since we cannot generally expect a pointwise bound for the integral $\int_M (\phi(x)-\phi(y))^2J(dx,dy)$ that appears in \CS. Note that
 \EC \   allows us to compare double integrals only, and  the constant  $K_0$ in \EC \ can be strictly larger than $1$. 
To address these issues, we use a covering argument to manage the constant  $K_0$ and take advantage of the fact that the cutoff $\phi \in \sF^s$ used in \CS \ for  $(\sE^s,\sF^s)$ is almost radial.

We begin with the next lemma.

\begin{lemma}\label{l:EC-modify}
	Suppose that \TJ \ and \EC \ hold. There exists $C>0$ independent of $s$  such that
	for all $x_0 \in M$,  $R>0$, $r \in (0,R_0/2)$ and $f\in L^2(B(x_0,R+r))$,
		\begin{align}\label{e:EC-modify-1}
		&(1-s)\int_{B(x_0,R)\times B(x_0,R)} \frac{(f(x)-f(y))^2}{V(x,d(x,y)) d(x,y)^{\beta s}} \mu(dx)\mu(dy) \nn\\
		&\le C 	\bigg(  \int_{B(x_0,R+r)\times B(x_0,R+r)}(f(x)-f(y))^2	J(dx,dy)  + \frac{1-s }{r^{\beta s}}  \int_{B(x_0,R)} f^2d\mu \bigg)
	\end{align}
	and
	\begin{align}\label{e:EC-modify-2}
		&\int_{B(x_0,R)\times B(x_0,R)} (f(x)-f(y))^2	J(dx,dy) \nn\\
		&\le C 	\bigg( (1-s) \int_{B(x_0,R+r)\times B(x_0,R+r)}  \frac{(f(x)-f(y))^2}{V(x,d(x,y)) d(x,y)^{\beta s}} \mu(dx)\mu(dy) + \frac{1-s }{r^{\beta s}}  \int_{B(x_0,R)} f^2d\mu \bigg).
	\end{align}
\end{lemma}
\pf Let $K_0 \ge 1$ be constant in \EC.  Set
$$
E_1:= \left\{ (x,y) \in B(x_0,R)\times B(x_0,R) : d(x,y) <r/(4K_0)\right\}
$$
and $E_2:= (B(x_0,R)\times B(x_0,R) ) \setminus E_1$.
 By Lemma \ref{l:VD-covering}, there exist  $N_0 \in \N$ and  an open  covering  $\{B(z_i,r/(4K_0))\}_{i=1}^\infty$ of $M$ such that
\begin{align}\label{e:VD-covering-appl-2}
	\sum_{i=1}^\infty \1_{B(z_i,r/2)}  \le N_0\;\; \text{on} \;\,  M.
\end{align}
Since  $\{B(z_i,r/(4K_0))\}_{i=1}^\infty$ is an open covering of $M$, we have
\begin{align*}
	&(1-s)\int_{E_1} \frac{(f(x)-f(y))^2}{V(x,d(x,y)) d(x,y)^{\beta s}} \mu(dx)\mu(dy) \\
	&\le (1-s) \sum_{i: B(z_i,r/(4K_0)) \cap B(x_0,R) \neq \emptyset} \int_{B(z_i,r/(4K_0))} \int_{B(x,r/(4K_0))}\frac{(f(x)-f(y))^2}{V(x,d(x,y)) d(x,y)^{\beta s}} \mu(dy) \mu(dx)\\
		&\le (1-s) \sum_{i: z_i \in  B(x_0,R+r/4)} \int_{B(z_i,r/(2K_0))\times B(z_i,r/(2K_0))}\frac{(f(x)-f(y))^2}{V(x,d(x,y)) d(x,y)^{\beta s}}  \mu(dx)\mu(dy).
\end{align*}
Hence, applying \EC \ and using \eqref{e:VD-covering-appl-2}, we get
\begin{align}\label{e:EC-modify-1-1}
	&(1-s)\int_{E_1} \frac{(f(x)-f(y))^2}{V(x,d(x,y)) d(x,y)^{\beta s}} \mu(dx)\mu(dy) \nn\\
	&\le c_1 \sum_{i: z_i \in  B(x_0,R+r/2)} \int_{B(z_i,r/2)\times B(z_i,r/2)}(f(x)-f(y))^2 J(dx,dy)\nn\\
		&\le c_1N_0 \int_{B(x_0,R+r)\times B(x_0,R+r)}(f(x)-f(y))^2 J(dx,dy).
\end{align}
On the other hand, since $V(x,d(x,y))\asymp V(y,d(x,y))$ for all $x,y \in M$ by \ref{e:VD2},  using Lemma \ref{l:robust-jumptail}, we obtain
\begin{align}\label{e:EC-modify-1-2}
	&(1-s)\int_{E_2} \frac{(f(x)-f(y))^2}{V(x,d(x,y)) d(x,y)^{\beta s}} \mu(dx)\mu(dy) \nn\\
		&\le c_2(1-s)\int_{B(x_0,R)}f(x)^2 \int_{B(x,r/(4K_0))^c} \frac{\mu(dy)}{V(x,d(x,y)) d(x,y)^{\beta s}} \mu(dx) \nn\\
	&\le \frac{c_3(1-s)}{\beta s(r/(4K_0))^{\beta s}}  \int_{B(x_0,R)} f^2d\mu \le \frac{c_3(1-s)(4K_0)^\beta }{\beta s_0r^{\beta s}}  \int_{B(x_0,R)} f^2d\mu.
\end{align}
Combining \eqref{e:EC-modify-1-1} with \eqref{e:EC-modify-1-2}, we deduce that \eqref{e:EC-modify-1} holds. 

 Similarly, following the arguments for \eqref{e:EC-modify-1}  and using Proposition \ref{p:EC-upper}, we obtain \eqref{e:EC-modify-2}.  \qed

Note that, by Corollary \ref{c:domain=domain}, $\wt \sF = \wt \sF^{s}$ and  $\wt \sF_b = \wt \sF^{s}_b$.  

\begin{lemma}\label{l:CS-pre}
	Suppose that \TJ \ and \EC \ hold. There exist  $C_1,C_2>0$  independent of $s$ such that the following  holds: Let $x_0 \in M$ and $R, r > 0$ be such that $R + 2r < R_0$.  There exists a cutoff function $\phi\in \sF$ for $B(x_0,R) \Subset B(x_0,R+r)$ such that  for all $f \in \wt \sF_b$,
	\begin{align*}
		&\int_{B(x_0,R+2 r) \times M}   f(x)^2(\phi(x)-\phi(y))^2 J(dx,dy) \nn\\
		&\le  C_1(1-s)\int_{B(x_0,R+ 2r) \times B(x_0,R+2r)} \frac{(f(x) - f(y))^2}{ V(x,d(x,y))d(x,y)^{\beta s}} \, \mu(dx)\mu(dy)  \nn\\
		&\quad +  2 \int_{B(x_0,R+ 2r) \times B(x_0,R+2 r)}   (f(x)-f(y))^2 J(dx,dy) + \frac{C_2}{r^{\beta s}}\int_{B(x_0,R+2r)} f^2 d\mu.
	\end{align*}	 
\end{lemma}
\pf  Define $U_a:=B(x_0,R+a r)$ for $a\ge 0$. 
By Proposition \ref{p:CS-0}, since $\sF^s=\sF$ and $\wt \sF^s=\wt \sF$, there exist constants $c_1,c_2>0$ independent of $s,x_0,R$ and $r$, and  a cutoff function $\phi \in \sF$ for $U_0 \Subset U_1$ such that for all $f \in \wt \sF$,
\begin{align}\label{e:CS}
	&(1-s)\int_{U_2 \times M}   \frac{f(x)^2(\phi(x)-\phi(y))^2}{V(x,d(x,y))d(x,y)^{\beta s}} \,\mu(dx)\mu(dy) \nn\\
	&\le  c_1(1-s)\int_{U_1 \times U_2} \frac{\phi(x)^2 (f(x) - f(y))^2}{ V(x,d(x,y))d(x,y)^{\beta s}}  \,\mu(dx)\mu(dy) + \frac{c_2}{r^{\beta s}}\int_{U_2} f^2 d\mu. 
\end{align}	 
Since $\phi$ is a cutoff function for $U_0 \Subset U_1$, using  \TJ, we get that for all $f \in \wt \sF$,
\begin{align}\label{e:CS-pre-1}
	&	\int_{U_2 \times U_2^c}   f(x)^2(\phi(x)-\phi(y))^2 J(dx,dy) = 	\int_{U_2 \times U_2^c}   f(x)^2\phi(x)^2 J(dx,dy) \nn\\
	&\le	\int_{U_1}   f(x)^2 J(x,B(x,r)^c) \mu(dx) \le \frac{c_3(1-s)}{r^{\beta s}} \int_{U_1} f^2 d\mu.
\end{align}
Besides, using  Lemma \ref{l:EC-modify} (with $R$ replaced by $R+2r$), since $\phi^2\le 1$,  we see that
\begin{align}\label{e:CS-pre-2}
	&	\int_{U_2 \times U_2}   f(x)^2(\phi(x)-\phi(y))^2 J(dx,dy) \nn\\
	&\le	2\int_{U_2 \times U_2}    \left( (f(x)\phi(x)-f(y)\phi(y))^2   + \phi(y)^2(f(y)-f(x))^2 \right) J(dx,dy) \nn\\
		&\le	c_4(1-s)\int_{U_3 \times U_3}   \frac{(f(x)\phi(x)-f(y)\phi(y))^2}{V(x,d(x,y))d(x,y)^{\beta s}} \mu(dx)\mu(dy)+ \frac{c_5(1-s)}{r^{\beta s}}  \int_{U_2} f^2d\mu\nn\\
		&\quad + 2 \int_{U_2 \times U_2} (f(y)-f(x))^2 J(dx,dy).
\end{align}
Applying \eqref{e:CS} and using $\phi^2\le 1$, we get 
\begin{align}\label{e:CS-pre-3}
	&(1-s)\int_{U_2 \times U_2}   \frac{(f(x)\phi(x)-f(y)\phi(y))^2}{V(x,d(x,y))d(x,y)^{\beta s}} \mu(dx)\mu(dy) \nn\\
	&\le 2(1-s)\int_{U_2 \times U_2}   \frac{(f(x)(\phi(x)-\phi(y))^2 + \phi(y)(f(x)-f(y))^2)}{V(x,d(x,y))d(x,y)^{\beta s}} \mu(dx)\mu(dy) \nn\\
	&\le 2c_1(1-s)\int_{U_2 \times U_2} \frac{ (f(x) - f(y))^2}{ V(x,d(x,y))d(x,y)^{\beta s}} \, \mu(dx)\mu(dy) + \frac{2c_2}{r^{\beta s}} \int_{U_2} f^2 d\mu\nn\\
	&\quad  +   2(1-s)\int_{U_2 \times U_2}   \frac{(f(x)-f(y))^2}{V(x,d(x,y))d(x,y)^{\beta s}} \mu(dx)\mu(dy).
\end{align}
By \ref{e:VD2},  $V(x,d(x,y))\asymp V(y,d(x,y))$ for all $x,y \in M$. Using this and  the fact that  $\phi=0$ in $U_1^c$ in the first inequality below, $\phi^2\le 1$ in the second   and Lemma \ref{l:robust-jumptail} in the third, we obtain
\begin{align}\label{e:CS-pre-4}
	&(1-s)\int_{(U_3 \times U_3) \setminus (U_2 \times U_2)}   \frac{(f(x)\phi(x)-f(y)\phi(y))^2}{V(x,d(x,y))d(x,y)^{\beta s}} \mu(dx)\mu(dy) \nn\\
	&\le c_6(1-s)\int_{U_1 \times U_2^c}   \frac{f(x)^2\phi(x)^2}{V(x,d(x,y))d(x,y)^{\beta s}} \mu(dx)\mu(dy) \nn\\
	&\le  c_6(1-s)\int_{U_1}  f(x)^2 \int_{B(x,r)^c} \frac{\mu(dy)}{V(x,d(x,y))d(x,y)^{\beta s}} \mu(dx)\nn\\
		&\le  \frac{c_7(1-s)}{\beta s_0r^{\beta s}}\int_{U_1} f^2 d\mu.
\end{align}
Combining  \eqref{e:CS-pre-1}, \eqref{e:CS-pre-2}, \eqref{e:CS-pre-3} and \eqref{e:CS-pre-4}, we arrive at the result. \qed

The proof of the next proposition is motivated by \cite[Lemma 2.9]{GHH18} and \cite[Proposition 2.4]{CKW-memo}. However, unlike in \cite{GHH18, CKW-memo}, since we do not assume the existence and pointwise estimates for the density of the jump kernel $J(dx,dy)$, non-trivial modifications are required.

\begin{prop}\label{p:CS-improve}
	Suppose that \TJ \ and \EC \ hold. Then $(\sE,\sF)$ satisfies \CS.
\end{prop}
\pf Let $x_0 \in M$ and $R,r>0$ be such that $R+2r<R_0$. Write $U_a:=B(x_0,R+a r)$ for $a\ge 0$. Fix $\eps>0$.  Let  $\lambda=\lambda(\eps)>0$  be a constant  to be  determined later.  Define $c_0:=1-e^{-\lambda/\beta}$, 
\begin{align*}
		a_n:= c_0e^{-\lambda (n-1)/\beta} \quad \text{and} \quad    b_n:=\sum_{m=1}^n a_m, \qquad n \ge 1.
\end{align*}
Note that $\lim_{n\to \infty} b_n=1$. By Lemma \ref{l:CS-pre}, there exist constants $C_1,C_2>0$ independent of $s,x_0,R,r$ and $\lambda$ such that for each  $n \ge 1$,  there exists a cutoff function $\phi_n\in \sF$ for $U_{b_n-a_{n+1}} \Subset U_{b_{n}}$ so that  for all $f \in \wt \sF_b$,
\begin{align}\label{e:improve-CS-1}
	&\int_{U_{b_{n+1}} \times M}   f(x)^2(\phi_n(x)-\phi_n(y))^2 J(dx,dy) \nn\\
	&\le   C_1(1-s)\int_{U_{b_{n+1}} \times U_{b_{n+1}}} \frac{(f(x) - f(y))^2}{ V(x,d(x,y))d(x,y)^{\beta s}} \, \mu(dx)\mu(dy) \nn\\
	&\quad +  2 \int_{U_{b_{n+1}} \times U_{b_{n+1}} }   (f(x)-f(y))^2 J(dx,dy)   + \frac{C_2e^{s\lambda  n}}{(c_0r)^{\beta s}}\int_{U_{b_{n+1}}} f^2 d\mu. 
\end{align}	
Define
\begin{align*}
	\phi = \sum_{n=1}^\infty \, (e^{-\lambda(n-1)}-e^{-\lambda n})\phi_n.
\end{align*}
Then $\phi \le \sum_{n=1}^\infty (e^{-\lambda(n-1)}-e^{-\lambda n}) =1 $ in $M$, $\phi=0$ in $U_1^c$  and  for every $k \ge 0$, 
\begin{align}\label{e:vp>ek}
	\phi \ge \sum_{n=k+1}^\infty(e^{-\lambda(n-1)}-e^{-\lambda n})\phi_n = \sum_{n=k+1}^\infty(e^{-\lambda(n-1)}-e^{-\lambda n})= e^{-\lambda k} \quad  \text{on}\;\, U_{b_k}.
\end{align}
In particular, $\phi$ is a cutoff function for $U_0\Subset U_1$.  Further, for each  $n \ge 1$,  applying 
\eqref{e:improve-CS-1} with $f=1$, we get
\begin{align*}
	&	\sE(\phi_n,\phi_n)\le 2\int_{U_{b_{n}} \times M}   (\phi_n(x)-\phi_n(y))^2 J(dx,dy)  \le \frac{C_2\mu(U_2) e^{s\lambda  n}}{(c_0r)^{\beta s}}.
\end{align*}	
Using this, since $\phi_n^2\le 1$ and  $s<1$, we get that  for all $N,k\ge 1$, 
\begin{align*}
	&\sE_1\bigg( \sum_{n=N+1}^{N+k} \, (e^{-\lambda(n-1)}-e^{-\lambda n})\phi_n ,\sum_{n=N+1}^{N+k} \, (e^{-\lambda(n-1)}-e^{-\lambda n})\phi_n \bigg)^{1/2} \\
	&\le \sum_{n=N+1}^{N+k} (e^{-\lambda(n-1)}-e^{-\lambda n})  \sE_1 ( \phi_n, \phi_n)^{1/2}  \le \sum_{n=N+1}^{\infty} e^{-\lambda(n-1)} \bigg( \frac{C_2\mu(U_2) e^{s\lambda  n}}{(c_0r)^{\beta s}} + \mu(U_2) \bigg)^{1/2} \\
	& \le  e^\lambda \mu(U_2)^{1/2} \bigg( \frac{C_2}{(c_0r)^{\beta s}} + 1\bigg)^{1/2} \sum_{n=N+1}^{\infty}   e^{-\lambda n/2} = \frac{e^\lambda \mu(U_2)^{1/2}}{e^{\lambda/2}-1} \bigg( \frac{C_2}{(c_0r)^{\beta s}} + 1\bigg)^{1/2}e^{-\lambda N/2}.
\end{align*}
Thus,  $(\sum_{n=1}^{k} \, (e^{-\lambda(n-1)}-e^{-\lambda n})\phi_n )_{k \ge 1}$ is a $\sE_1$-Cauchy sequence, implying $\phi \in \sF$.

Now, we  show that $\phi$ satisfies \eqref{e:rCS}. 
Let $f \in \wt\sF_b$.  Observe that
\begin{align*}
	&\int_{U_2 \times M}   f(x)^2(\phi(x)-\phi(y))^2 J(dx,dy) \\
	&= 2(e^\lambda -1)^2\sum_{n=1}^\infty \sum_{m=n+2}^{\infty} e^{-\lambda (n+m)} \int_{U_{2}\times M}  f(x)^2 (\phi_n(x)-\phi_n(y))(\phi_m(x)-\phi_m(y))J(dx,dy)\nn\\
	&\quad  + 2(e^\lambda -1)^2 \sum_{n=1}^\infty e^{-\lambda (2n+1)}  \int_{U_{2}\times M}  f(x)^2 (\phi_n(x)-\phi_n(y))(\phi_{n+1}(x)-\phi_{n+1}(y)) J(dx,dy)\nn\\
	&\quad + (e^\lambda-1)^2\sum_{n=1}^\infty  e^{-2\lambda n} \int_{U_{2}\times M} f(x)^2 (\phi_n(x)-\phi_n(y))^2J(dx,dy)\nn\\
	&=:I_1+I_2+I_3.
\end{align*}
For  $n \ge 1$ and $m \ge n+2$, we see that $(\phi_n(x)-\phi_n(y))(\phi_m(x)-\phi_m(y)) \neq 0$ only if $x \in U_{b_{n}}$ and $y \in U_{b_{n+1}}^c$, or $x \in U_{b_{n+1}}^c$ and $y \in U_{b_{n}}$. Thus, for any  $n \ge 1$ and $m \ge n+2$, $(\phi_n(x)-\phi_n(y))(\phi_m(x)-\phi_m(y)) \neq 0$ only if $d(x,y) \ge a_{n+1}r$.  Using this property and \TJ,  we obtain
\begin{align}\label{e:I1}
	I_1 	&\le 2(e^\lambda -1)^2 \sum_{n=1}^\infty  \sum_{m=n+2}^\infty e^{-\lambda (n+m)}  \int_{U_2} f(x)^2 J(x,B(x, a_{n+1}r)^c)\mu(dx)\nn\\
	&\le \frac{c_1(1-s)(e^\lambda -1)^2}{(c_0r)^{\beta s}} \sum_{n=1}^\infty  \sum_{m=n+2}^\infty e^{-\lambda (n+m-sn)} \int_{U_2} f^2d\mu\nn\\
	&\le \frac{c_1(e^\lambda -1)^2}{c_0^{\beta }r^{\beta s}} \sum_{n=1}^\infty  \sum_{m=n+2}^\infty e^{-\lambda m} \int_{U_2} f^2d\mu\nn\\
	&=  \frac{c_1(e^\lambda -1)}{e^{\lambda} c_0^{\beta }r^{\beta s}} \sum_{n=1}^\infty  e^{-\lambda n} \int_{U_2} f^2d\mu =  \frac{c_1}{e^\lambda  c_0^{\beta }r^{\beta s}}  \int_{U_2} f^2d\mu.
\end{align}
By the Cauchy inequality,  $I_2 \le 2I_3$. For each $n$, since $\phi_n=0$ in $U_{b_n}^c$,  by   \TJ, 
\begin{align}\label{e:improve-CS-I3-1}
	&\int_{(U_2 \setminus U_{b_{n+1}}) \times M}   f(x)^2(\phi_n(x)-\phi_n(y))^2 J(dx,dy)      =\int_{(U_2 \setminus U_{b_{n+1}}) \times U_{b_{n}}}   f(x)^2\phi_n(y)^2 J(dx,dy)      \nn\\
	&\le  \int_{U_2 \setminus U_{b_{n+1}}} f(x)^2  J(x,B(x, a_{n+1} r)^c) \mu(dx) \le \frac{c_2(1-s)e^{s\lambda n}}{(c_0r)^{\beta s}} \int_{U_{2}} f^2 d\mu.
\end{align}
Further,  applying Lemma \ref{l:EC-modify} with $R$ replaced by $R+b_nr$ and $r$ by $a_{n+2}r$, we see that 
\begin{align}\label{e:improve-CS-I3-2}
	&(1-s)\int_{U_{b_{n+1}} \times U_{b_{n+1}}} \frac{(f(x) - f(y))^2}{ V(x,d(x,y))d(x,y)^{\beta s}} \, \mu(dx)\mu(dy) \nn\\
	&\le c_3\int_{U_{b_{n+2}} \times U_{b_{n+2}}}  (f(x) - f(y))^2 J(dx,dy)  + \frac{(1-s)e^{s\lambda(n+1)}}{(c_0r)^{\beta s}} \int_{U_2} f^2 d\mu 
\end{align}
Combining  \eqref{e:improve-CS-1}, \eqref{e:improve-CS-I3-1} and \eqref{e:improve-CS-I3-2}, we deduce that for all $n \ge 1$,
\begin{align*}
	&\int_{U_{2}\times M} f(x)^2 (\phi_n(x)-\phi_n(y))^2J(dx,dy) \nn\\
		& \le c_4 \bigg( \int_{U_{b_{n+2}} \times U_{b_{n+2}} }   (f(x)-f(y))^2 J(dx,dy)   + \frac{e^{s\lambda  (n+1)}}{(c_0r)^{\beta s}}\int_{U_{2}} f^2 d\mu\bigg). 
\end{align*}
Therefore, we obtain
\begin{align*}
	I_3&\le   c_4(e^\lambda-1)^2\sum_{n=1}^\infty  e^{-2\lambda n}  \int_{U_{b_{n+2}} \times U_{b_{n+2}} }   (f(x)-f(y))^2 J(dx,dy)  \\
	&\quad +  c_4e^{s\lambda }(e^\lambda-1)^2\sum_{n=1}^\infty   \frac{e^{-(2-s)\lambda n}}{(c_0r)^{\beta s}}   \int_{U_{2}} f^2 d\mu   \\
	&=:I_{3,1}+I_{3,2}.
\end{align*}
Note that 
\begin{align}\label{e:I3-2}
	I_{3,2}&\le \frac{c_4 e^\lambda (e^\lambda -1)^2  }{ c_0^{\beta}  r^{\beta s}} \sum_{n=1}^\infty e^{-\lambda n}\int_{U_2} f^2 d\mu = \frac{c_4e^\lambda (e^\lambda -1) }{ c_0^{\beta } r^{\beta s}}\int_{U_2} f^2d\mu.
\end{align}
For $I_{3,1}$,   we have
\begin{align}\label{e:I3-1}
	I_{3,1}	&\le c_4(e^\lambda-1)^2\sum_{n=1}^\infty  e^{-2\lambda n}  \int_{U_{0} \times U_{b_{n+2}} }    (f(x)-f(y))^2 J(dx,dy)\nn\\
	&\quad + c_4(e^\lambda-1)^2\sum_{n=1}^\infty  \sum_{k=1}^{n+2} e^{-2\lambda n}  \int_{(U_{b_k} \setminus U_{b_{k-1}}) \times U_{b_{n+2}} }    (f(x)-f(y))^2 J(dx,dy)\nn\\
	&=:c_4(I_{3,1}' + I_{3,2}' ).
\end{align}
Since $\phi=1$ in $U_0$, it holds that 
\begin{align}\label{e:I3-1'}
	I_{3,1}'&\le (e^\lambda-1)^2\sum_{n=1}^\infty  e^{-2\lambda n}\int_{U_{0} \times U_{2} }    \phi(x)^2(f(x)-f(y))^2 J(dx,dy)\nn\\
	&= \frac{(e^\lambda-1)}{e^\lambda + 1}\int_{U_{0} \times U_{2} }    \phi(x)^2(f(x)-f(y))^2 J(dx,dy).
\end{align}
Moreover, by \eqref{e:vp>ek},  we have
\begin{align}\label{e:I3-2'}
	I_{3,2}' &\le   (e^\lambda-1)^2\sum_{n=1}^\infty  \sum_{k=1}^{n+2} e^{-2\lambda (n-k)}  \int_{(U_{b_k} \setminus U_{b_{k-1}}) \times U_{b_{n+2}} }  \phi(x)^2  (f(x)-f(y))^2 J(dx,dy)\nn\\
	& \le    (e^\lambda -1)^2   \sum_{k=1}^\infty e^{2\lambda k}  \int_{(U_{b_{k}} \setminus U_{b_{k-1}}) \times U_{2}} \phi(x)^2  (f(x)-f(y))^2 J(dx,dy) \sum_{n=k-2}^\infty e^{-2\lambda n}\nn\\
	& =   \frac{ (e^\lambda -1)e^{6\lambda}}{e^{\lambda}+1} \sum_{k=1}^\infty  \int_{(U_{b_{k}} \setminus U_{b_{k-1}}) \times U_{2}} \phi(x)^2  (f(x)-f(y))^2 J(dx,dy)  \nn\\
	& =  \frac{ (e^\lambda -1)e^{6\lambda}}{e^{\lambda}+1}  \int_{(U_{1} \setminus U_{0}) \times U_{2}} \phi(x)^2  (f(x)-f(y))^2 J(dx,dy)  .
\end{align}
Since $I_2\le 2I_3$, by  \eqref{e:I3-2},   \eqref{e:I3-1}, \eqref{e:I3-1'} and \eqref{e:I3-2'}, we deduce that
\begin{align*}
	I_2 + I_3\le \frac{ 3c_4(e^\lambda -1)e^{6\lambda}}{e^{\lambda}+1}  \int_{U_{1} \times U_{2}} \phi(x)^2  (f(x)-f(y))^2 J(dx,dy)  +   \frac{3c_4e^\lambda(e^\lambda -1) }{ c_0^{\beta } r^{\beta s}}\int_{U_2} f^2d\mu.
\end{align*}
Combining this with \eqref{e:I1}, we arrive at
\begin{align*}
	&\int_{U_2 \times M}   f(x)^2(\phi(x)-\phi(y))^2 J(dx,dy) \\
	&\le \frac{ 3c_4(e^\lambda -1)e^{6\lambda}}{e^{\lambda}+1}  \int_{U_{1} \times U_{2}} \phi(x)^2  (f(x)-f(y))^2 J(dx,dy)  +   \frac{c_1+3c_4e^\lambda(e^{\lambda} -1) }{ c_0^{\beta } r^{\beta s}}\int_{U_2} f^2d\mu .
\end{align*}
Note that $\lim_{a \to 0} 3c_4(e^a -1)e^{6a}/(e^a+1)=0$.
By choosing $\lambda$  such that   $ 3c_4(e^\lambda -1)e^{6\lambda}/(e^\lambda+1)=\eps$,  we conclude that  \eqref{e:rCS} holds. The proof is complete.\qed

\smallskip 

\noindent \textbf{Proof of Theorem \ref{t:main-1}.} The result follows from Corollary \ref{c:domain=domain}, and Propositions \ref{p:FK-PI} and \ref{p:CS-improve}. \qed

\section{$L^2$-mean value inequality for $(\sE,\sF)$}\label{section:mean-value-inequality}

In Sections \ref{section:mean-value-inequality} and \ref{section:WEHI-EHR}, we assume that $(\sE,\sF)$ is a pure-jump type regular Dirichlet form on $L^2(M)$ with the representation \eqref{e:def-sE}. In this section, we establish the $L^2$-mean value inequality for subsolutions associated with $(\sE,\sF)$  (Proposition \ref{p:L2mean}). We mainly follow the strategy of \cite[Section 4]{CKW-memo}.

The next lemma follows from \cite[Lemma 3.2(i)]{GHH18}.
\begin{lemma}\label{l:supersolution}
	Let $D \subset M$ be a non-empty bounded open set and  $F:\R\to \R$ be a twice differentiable function such that  $F'' \ge 0$ and  $	\sup_\R |F'| + \sup_\R F''<\infty$.
	Suppose  that 	$u \in \sF^{\rm loc}_{D}$ is locally bounded in $D$ and satisfies \eqref{e:harmonicity-1}.  Then for any $0\le \phi \in \sF\cap C_c(D)$,  $	\sE(F\circ u, \phi)$ and $\sE(u, (F'\circ u)\phi )$
	are absolutely convergent and
	\begin{align}\label{e:supersolution}
		\sE(F\circ u, \phi) \le \sE(u, (F'\circ u)\phi ).
	\end{align}   
\end{lemma}

We establish a robust version of   Caccioppoli-type inequality. Note that a non-robust version was previously  established in \cite[Lemma 4.6]{CKW-memo}.

 Recall that the nonlocal tail $\NT(u, D_1, D_2)$ of $u$ is defined as \eqref{e:def-tail}.
\begin{lemma}\label{l:Caccio}
	Suppose that \TJ \ 
	and \CS \ hold. There exists $C>0$ independent of $s$ such that the following holds: For any  $x_0 \in M$ and $R,r>0$ satisfying $R+2r<R_0$, there exists a cutoff function $\phi\in \sF$ for  $B(x_0,R) \Subset B(x_0,R+r)$ such that  if $u$ is   bounded in $B(x_0,R+2r)$ and  $-\sL u \le f$ in $B(x_0,R+2r)$ for  $f \in L^\infty(B(x_0,R+2r))$, then for all $\theta_0\ge 0$ and $\theta_1> \theta_0 +  r^{\beta s} \lVert f_+ \rVert_{L^\infty(B(x_0,R+2r))}$,
	\begin{align*}
		&\int_{B(x_0,R+2r)\times M} \big(\phi(x)  (  u(x) - \theta_1)_+ - \phi(y) ( u(y) - \theta_1)_+  \big)^2 J(dx,dy)\nn\\
		& \le   \bigg(   \frac{C}{r^{\beta s}} +  \frac{18	\NT(u_+,B(x_0,R+r), B(x_0,R+2r)^c) }{\theta_1-  \theta_0} 
		\bigg)\int_{B(x_0,R+2r)}  ( u(x) - \theta_0)^2_+ d\mu.
	\end{align*}
\end{lemma}
\pf  Define  $U_a:=B(x_0,R+ar)$ for $a\ge 0$. Let $\phi \in \sF$ be a cutoff function for $U_0\Subset U_1$ satisfying \eqref{e:rCS} with $\eps=1/8$.   Denote $T:=	\NT(u_+,U_1, U_2^c)$, $v:=(u -\theta_1)_+$ and $w:=(u-\theta_0)_+$. Since $u \in \sF^{\rm loc}_{U_2} \cap L^\infty(U_2)$, we have $v \in \sF^{\rm loc}_{U_2} \cap L^\infty(U_2)$. Hence, there exists $\wt v \in \sF_b$ such that $v=\wt v$ in $U_1$. By \cite[Theorem 1.4.2]{FOT},   $v \phi^2  = \wt v\phi^2 \in \sF_{U_1}\cap L^\infty(M)$. Since $-\sL u \le f$ in $U_{2}$,  we obtain
\begin{align}\label{e:CKWp.36-0}
	\int_{U_1}  f v\phi^2 d\mu \ge \sE( u,  v\phi^2)	 & = \int_{U_{2} \times U_{2}} (u(x)-u(y)) (v(x)\phi(x)^2 - v(y)\phi(y)^2)J(dx,dy)
	\nn\\
	&\quad + 2 \int_{U_{1} \times U_{2}^c} ( u(x)- u(y))v(x) \phi(x)^2 J(dx,dy)\nn\\
	&=:I_1+I_2.
\end{align}
For all $x \in U_1$, if $u(x) \le \theta_1$, then $f(x)v(x)\phi(x)^2=0$ and if $u(x) >\theta_1$, then 
$$f(x)v(x) \phi(x)^2 \le \lVert f_+\rVert_{L^\infty(U_2)}v(x) \le  r^{-\beta s} (\theta_1-\theta_0)(u(x)-\theta_1)  \le r^{-\beta s}w(x)^2.$$
Thus, $f(x)v(x)\phi(x)^2 \le r^{-\beta s}w(x)^2$ for all $x \in U_1$, implying that
\begin{align}\label{e:CKWp.36-00}
	\int_{U_1}  f v\phi^2 d\mu  \le \frac{1}{r^{\beta s}}\int_{U_1} w^2 d\mu.
\end{align}
For  $I_1$, following the argument in  \cite[p.36]{CKW-memo}, we get
\begin{align}\label{e:estimate-I1}
	I_1 \ge \frac12 \int_{U_2\times U_2} \phi(x)^2 (v(x)-v(y))^2 J(dx,dy) - 2 \int_{U_2 \times U_2} v(x)^2(\phi(x)-\phi(y))^2 J(dx,dy).
\end{align}
For $I_2$, we note that for all $x \in U_1$ and $y \in U_2^c$, 
if  $u(x) \ge u(y)$ or $u(x)<\theta_1$, then
$
(u(x)-u(y))v(x)\phi(x)^2  \ge0,
$
and if $u(y)>u(x) \ge \theta_1$, then 
\begin{align*}
	(u(x)-u(y))v(x)\phi(x)^2  \ge  - u(y)(u(x)-\theta_1) \ge - \frac{u_+(y)w(x)^2}{\theta_1 -  \theta_0  }.
\end{align*}
Consequently, it holds that 
\begin{align}\label{e:estimate-I2}
	I_2&\ge   -2\int_{U_1} \frac{w(x)^2}{\theta_1-  \theta_0}  \int_{ U_{2}^c} u_+(y) J(x,y) \mu(dy)\,\mu(dx)\ge  -\frac{2T}{\theta_1-\theta_0}\int_{U_1}w^2d\mu.
\end{align}
Combining  \eqref{e:CKWp.36-0}, \eqref{e:CKWp.36-00}, \eqref{e:estimate-I1} and \eqref{e:estimate-I2},   we obtain
\begin{align}\label{e:Caccio-1}
	& \frac12 \int_{U_2\times U_2} \phi(x)^2 (v(x)-v(y))^2 J(dx,dy)\nn\\
	&\le   2 \int_{U_2 \times U_2} v(x)^2(\phi(x)-\phi(y))^2 J(dx,dy)+ \bigg( \frac{1}{r^{\beta s}}  +  \frac{2T}{\theta_1-\theta_0 } \bigg) \int_{U_1} w^2 d\mu.
\end{align}
By the symmetry of $J$, since $\phi=0$ in $U_2^c$, we have
\begin{align*}
	&\frac19 \int_{U_{2}\times M} (\phi(x)v(x)-\phi(y)v(y))^2 J(dx,dy)\nn\\
	&= \frac29\int_{U_{2}\times M} v(x)^2(\phi(x)-\phi(y))^2 J(dx,dy)  + \frac29\int_{U_{2} \times U_2} \phi(x)^2 (v(x)-v(y))^2J(dx,dy).
\end{align*}
Using this in the first inequality below,  \eqref{e:rCS}  (with  $f=v$ and  $\eps=1/8$)  in the second, and \eqref{e:Caccio-1}  and  $v^2\le w^2$ in the last, we arrive at
\begin{align*}
	&\frac19 \int_{U_{2}\times M} (\phi(x)v(x)-\phi(y)v(y))^2 J(dx,dy)\nn\\
	&\le \frac{2}{9}\int_{U_{2}\times M} v(x)^2(\phi(x)-\phi(y))^2 J(dx,dy)+ \frac29\int_{U_{2} \times U_2} \phi(x)^2 (v(x)-v(y))^2J(dx,dy) \nn\\
	&\le  \frac12\int_{U_{2} \times U_2} \phi(x)^2 (v(x)-v(y))^2J(dx,dy)  +  \frac{c_1}{r^{\beta s}}\int_{U_2} v^2 d\mu -2\int_{U_{2}\times M} v(x)^2(\phi(x)-\phi(y))^2 J(dx,dy) \nn\\
	&\le \bigg( \frac{c_1+1}{r^{\beta s}} + \frac{2T}{\theta_1-\theta_0} \bigg) \int_{U_2}w^2d\mu.
\end{align*}
The proof is complete.\qed

The proof of the next lemma is originally due to  \cite[Lemma 3.2]{Gr91}.

\begin{lemma}\label{l:Gr0-Lemma3.2}
	Suppose that \TJ, \FK \ and \CS \ hold. 	There exists $C>0$ independent of $s$ such that the following holds: Let  $x_0 \in M$ and $0<r\le R$ be such that  $R+2r<R_0/2$. Suppose that  	$u$ is  bounded  in $B(x_0,R+2r)$ and  $-\sL u \le f$ in $B(x_0,R+2r)$ for $f \in L^\infty(B(x_0,R+2r))$. For  given  $\theta_0\ge 0$ and  $\theta_1 >\theta_0+r^{\beta s}\lVert f_+ \rVert_{L^\infty(B(x_0,R+2r))}$, define
	\begin{align*}
		\sI_0:= \int_{B(x_0,R+2r)} (u-\theta_0)_+^2 d\mu \quad \text{and} \quad \sI_1:= \int_{B(x_0,R)} (u-\theta_1)_+^2 d\mu.
	\end{align*}
	Then we have
	\begin{align*}
		\sI_1\le \frac{CR^{\beta s} }{(\theta_1 -  \theta_0)^{2\beta s/ \alpha} V(x_0,R+2r)^{\beta s/ \alpha}} \bigg[    \frac{1}{r^{\beta s}} +  \frac{\NT(u_+,B(x_0,R+r),B(x_0,R+2r)^c)}{\theta_1 -  \theta_0}  
		\bigg]
		\, \sI_0^{1+\beta s / \alpha}.
	\end{align*}
\end{lemma}
\pf  Write   $U_a:=B(x_0,R+ar)$ for $a\ge 0$. Denote  $w:=(u-\theta_0)_+$, $v:=(u-\theta_1)_+$ and $T:=	\NT(u_+,U_1,U_2^c) $. 
By Lemma \ref{l:Caccio}, there exists a cutoff function $\phi$ for $U_0 \Subset U_1$ such that 
\begin{equation}\label{e:Gr0-Lemma3.2-1}
	\int_{U_2\times M} \big(\phi(x)v(x)  - \phi(y)v(y)  \big)^2 J(dx,dy)\le   \bigg(  \frac{c_1}{r^{\beta s}} +  \frac{18T}{\theta_1 - \theta_0} 
	\bigg)\, \sI_0.
\end{equation}
Let
\begin{align*}
	E:=\left\{x \in U_{1}:v(x)>0\right\}=\left\{x \in U_{1}:u(x)>\theta_1\right\}.
\end{align*} 
If $\mu(E)=0$, then $\sI_1=0$ and the desired inequality is evident. Suppose that $\mu(E)>0$. By the outer regularity of    $\mu$, there is an open set $D$ such that  $E\subset D \subset U_2$ and  $\mu(D) \le 2\mu(E)$. 
By the Markov's inequality, it follows that
\begin{align}\label{e:Gr0-Lemma3.2-2}
	\mu(D) \le  2\mu(E) \le 2\int_{E} \frac{(u(x)-\theta_0)_+^2}{(\theta_1-\theta_0)^2} \mu(dx) \le \frac{2\sI_0}{(\theta_1-\theta_0)^2}.
\end{align}
Since  $\phi v \in \sF_D$ and $(\sE,\sF)$ satisfies \FK \ (with $K_2\ge 1$), we have
\begin{align*}
	&\lVert \phi v \rVert_2^2\left( \frac{V(x_0,R+2r)}{\mu(D)}\right)^{\beta s/ \alpha}	\nn\\
	& \le c_2\bigg((3R)^{\beta s}\int_{B(x_0,K_2R+2K_2r)\times B(x_0,K_2R+2K_2r)} \big(\phi(x)v(x)  - \phi(y)v(y)  \big)^2 J(dx,dy)  + \lVert \phi v \rVert_2^2\bigg).
\end{align*}
Thus, since $\phi=0$ in $U_1^c$, using the symmetry of $J$, \eqref{e:Gr0-Lemma3.2-1} and \TJ, we obtain
\begin{align}\label{e:Gr0-Lemma3.2-3}
	&\lVert \phi v \rVert_2^2\left( \frac{V(x_0,R+2r)}{\mu(D)}\right)^{\beta s/ \alpha}	\nn\\
	& \le c_2\bigg[ (3R)^{\beta s}\int_{U_2\times U_2} \big(\phi(x)v(x)  - \phi(y)v(y)  \big)^2 J(dx,dy)\nn\\
	&\qquad\quad   + 2(3R)^{\beta s}\int_{U_1} \phi(x)^2 v(x)^2 \int_{U_2^c} J(x,dy) \mu(dx)  + \lVert \phi v \rVert_2^2\bigg] \nn\\
	& \le 3^{\beta}c_2\bigg[ R^{\beta s}\bigg(  \frac{c_1}{r^{\beta s}} +  \frac{18T}{\theta_1 - \theta_0} 
	\bigg)\, \sI_0 + 2R^{\beta s}\int_{U_1} \phi(x)^2 v(x)^2 \int_{B(x,r)^c} J(x,dy) \mu(dx)  +  \lVert \phi v \rVert_2^2 \bigg] \nn\\
	&\le 3^\beta c_2 R^{\beta s}\bigg[ \bigg(  \frac{c_1}{r^{\beta s}} +  \frac{18T}{\theta_1 - \theta_0} 
	\bigg)\, \sI_0   + \bigg( \frac{c_3(1-s)+1}{r^{\beta s}}   \bigg)  \lVert \phi v \rVert_2^2\bigg].
\end{align}
Since $\phi=1$ in $U_0$ and $\phi^2v^2\le w^2$, we have
\begin{align}\label{e:Gr0-Lemma3.2-4}
	\sI_1= \int_{U_0}  \phi^2 v^2 d\mu \le  \lVert \phi v \rVert_2^2 \le \sI_0.
\end{align}
Combining  \eqref{e:Gr0-Lemma3.2-2}, \eqref{e:Gr0-Lemma3.2-3} and  \eqref{e:Gr0-Lemma3.2-4}, we arrive at
\begin{align*}
	\sI_1 	&\le  \lVert \phi v \rVert_2^2\le 3^\beta c_2 R^{\beta s} \left( \frac{\mu(D)}{V(x_0,R+2r)}\right)^{\beta s/ \alpha} \bigg[ \bigg(  \frac{c_1}{r^{\beta s}} +  \frac{18T}{\theta_1 - \theta_0} 
	\bigg)\, \sI_0   + \bigg( \frac{c_3(1-s)+1}{r^{\beta s}}   \bigg)  \lVert \phi v \rVert_2^2\bigg] \\
	&\le 3^\beta c_2R^{\beta s} \bigg(\frac{2\sI_0}{(\theta_1-\theta_0)^2V(x_0,R+2r)} \bigg)^{\beta s/\alpha} \bigg[  \frac{c_1+c_3+1}{r^{\beta s}} +  \frac{18T}{\theta_1 - \theta_0}   \bigg] \,\sI_0 \\
	&\le \frac{2^{\beta/\alpha}3^\beta c_2R^{\beta s} }{(\theta_1 -  \theta_0)^{2\beta s/ \alpha} V(x_0,R+2r)^{\beta s/ \alpha}} \bigg[    \frac{c_1+c_3+1}{r^{\beta s}} +  \frac{18T}{\theta_1 - \theta_0}  
	\bigg]
	\, \sI_0^{1+\beta s / \alpha} .
\end{align*} \qed

We recall the following elementary iteration lemma from \cite[Lemma 4.9]{CKW-memo}.
\begin{lemma}\label{l:iteration}
	Let $(a_n)_{n \ge 0}$ be a sequence of positive real numbers such that for all $n \ge 0$,
	$$
	a_{n+1} \le c_0 b^n a_n^{1+\eps},
	$$
	for some constants $\eps>0$, $b>1$ and $c_0>0$. If $a_0 \le c_0^{-1/\eps}\,b^{-1/\eps^2}$, 
	then $\lim_{n \to \infty}a_n=0$.
\end{lemma}

We now establish the $L^2$-mean value inequality for subsolutions.

\begin{prop}\label{p:L2mean} 
	Suppose that \TJ, \FK \ and \CS \ hold.  There exists $C>0$ independent of $s$ such that the following holds: Let  $x_0 \in M$ and $0<R<R_0/2$. Suppose that $u$ is bounded in $B(x_0,R)$ and $-\sL u \le f$ in $B(x_0,R)$ for $f \in L^\infty(B(x_0,R))$. Then for all $\delta>0$ and $q \in [0,\beta]$, it holds that
	\begin{align*}
		\esssup_{ B(x_0, R/2)} u& \le C \left[  \left(1+\frac1\delta\right)^{ \alpha/(2\beta s)} \sI^{1/2}+   \delta  \sT_q + R^{\beta s}\lVert f_+ \rVert_{L^\infty(B(x_0,R))}\right],
	\end{align*}
	where 
	\begin{align}
		\sI=\sI(u,x_0,R)&:= \frac{1}{V(x_0,R)}\int_{B(x_0,R)} u^2  d\mu, \nn\\
		\sT_q=\sT_q(u,x_0,R)&:=\sup \left\{   \frac{a^q\sT(u_+, B(x_0,r), B(x_0,r+a)^c)}{R^{q-\beta s}}\,: r \in [R/2,R],\, a \in (0,R/8]\right\}.\nn\\\label{e:L2-Caccio}
	\end{align}
\end{prop}
\pf    Fix $\delta>0$ and  $q \in [0,\beta]$.  Set $\nu:=\beta s/ \alpha$ and  let $\theta>R^{\beta s} \lVert f_+ \rVert_{L^\infty(B(x_0,R))}$ be a constant  to be determined later.  For  $n \ge 0$, define 
\begin{align*}
	r_n=\frac{1}{2}(1+2^{-n})R, \quad\;\; \theta_n=\sum_{j=0}^n   2^{-j\beta s}\theta   \quad \text{and} \quad a_n= \int_{B(x_0,r_n)} (u-\theta_n)_+^2 d\mu.
\end{align*}
For all $n \ge 0$,  we have
$$\theta_{n+1}-\theta_n = 2^{-(n+1)\beta s} \theta > 2^{-(n+1)\beta s}R^{\beta s} \lVert f_+ \rVert_{L^\infty(B(x_0,R))} > (r_n-r_{n+1})^{\beta s}\lVert f_+ \rVert_{L^\infty(B(x_0,R))}.$$
Applying Lemma \ref{l:Gr0-Lemma3.2} (with $R=r_{n+1}$, $r=(r_{n}-r_{n+1})/2$, $\theta_1 = \theta_{n+1}$ and $ \theta_0 = \theta_{n}$) and using \ref{e:VD},   we get that for all $n \ge 0$,
\begin{align}\label{e:iteration-1}
	a_{n+1}&\le \frac{c_1 r_{n+1}^{\beta s} }{(\theta_{n+1} -  \theta_n)^{2\nu} V(x_0,r_n )^{\nu}}  \nn\\
	&\qquad \times \bigg[    \frac{2^{\beta s}}{(r_n-r_{n+1})^{\beta s}} +  \frac{\NT(u_+,B(x_0,(r_n+r_{n+1})/2),B(x_0,r_n)^c)}{\theta_{n+1}-\theta_n}  
	\bigg]
	a_n^{1+\nu}\nn\\
	&\le \frac{c_1 r_{n+1}^{\beta s} }{(\theta_{n+1} -  \theta_n)^{2\nu} V(x_0,R/2)^{\nu}} \bigg[    \frac{2^{\beta s}}{(r_n-r_{n+1})^{\beta s}} +  \frac{R^{q-\beta s}\NT_q}{(\theta_{n+1}-\theta_n) ((r_n-r_{n+1})/2)^q}  
	\bigg]
	a_n^{1+\nu}\nn\\
	&\le \frac{2^{2\nu \alpha+2\nu \beta s(n+1)}c_2 R^{\beta s} }{ \theta^{2\nu} V(x_0,R)^{\nu}} \bigg[    \frac{2^{\beta s (n+3)}}{R^{\beta s}} +  \frac{2^{\beta s (n+1)+q(n+3)}  \sT_q}{  \theta R^{\beta s}}  
	\bigg]
	a_n^{1+\nu}\nn\\
	&\le  \frac{2^{2\nu\alpha + 2\nu \beta s+3\beta s  + 3q}c_2 }{\theta ^{2\nu} V(x_0,R)^{\nu}} \bigg[    1+  \frac{  \sT_q}{\theta }  
	\bigg] 2^{(2\nu \beta s+\beta  s+ q)n}
	a_n^{1+\nu}.
\end{align}
Without loss of generality, we assume $c_2\ge 1$. Set $b:= 2^{2\nu \beta s+\beta  s+ q}$, $c_3:=2^{2\nu\alpha+2\nu \beta s+3\beta s  + 3q}c_2  $,
\begin{align*}
	\theta:=  \bigg(\frac{c_3^{1/\nu}b^{1/\nu^2} (1+\delta^{-1})^{1/\nu} a_0}{V(x_0,R)}\bigg)^{1/2} +\delta  \sT_q  + R^{\beta s}\lVert f_+ \rVert_{L^\infty(B(x_0,R))}
\end{align*}
and
$$
c_0:=\frac{c_3}{\theta^{2\nu} V(x_0,R)^{\nu}} \bigg[    1+  \frac{\sT_q}{\theta }  
\bigg].
$$
By \eqref{e:iteration-1},  $a_{n+1} \le c_0 b^n a_n^{1+\nu}$ for all $n \ge 0$. Moreover, it holds that 
\begin{align*}
	&c_0^{-1/\nu} b^{-1/\nu^2} = c_3^{-1/\nu} b^{-1/\nu^2} \frac{ \theta^{2} V(x_0,R)}{(1+ \sT_q/\theta)^{1/\nu}}
	\ge \frac{(1+\delta^{-1})^{1/\nu}a_0}{ (1 + \sT_q/\theta)^{1/\nu}} 	\ge \frac{(1+\delta^{-1})^{1/\nu}a_0}{ (1 +\delta^{-1})^{1/\nu}}=a_0.
\end{align*}
Thus, by Lemma \ref{l:iteration}, we obtain $\lim_{n \to \infty} a_n=0$ which implies 
$	\esssup_{B(x_0,R/2)} u \le \sum_{j=0}^\infty 2^{-j\beta s}\theta =2^{\beta s}\theta/(2^{\beta s}-1).$
Since $a_0 \le \sI^2V(x_0,R)$,  by the definition of $\theta$, we arrive at
\begin{align*}
	&	\esssup_{B(x_0,R/2)} u \\
	&\le \frac{2^{\beta }}{2^{\beta s_0}-1} \bigg[\Big( 2^{7\alpha + 2\beta  + 3\alpha/s_0  + \alpha^2/(\beta s_0) + \alpha^2/(\beta s_0^2)}(1+\delta^{-1})^{ 1/\nu}\Big)^{\frac12} \sI  +\delta \sT_q + R^{\beta s}\lVert f_+ \rVert_{L^\infty(B(x_0,R))} \bigg].
\end{align*}   \qed

\begin{cor}\label{c:L2mean-diff-form}
	Suppose that \TJ, \FK \ and \CS \ hold.  There exists $C>0$ independent of $s$ such that the following holds:   Let  $x_0 \in M$ and   $0<R<R_0/2$. Suppose that $u$ is bounded in $B(x_0,R)$ and  $-\sL u \le f$ in $B(x_0,R)$ for $f \in L^\infty(B(x_0,R))$. Then for all $q \in [0,\beta]$, 
	\begin{equation}\label{e:L2mean-diff-form-1}
		\esssup_{B(x_0, R/2)} u  \le C\left[ \sI^{2\beta s/(\alpha+2\beta s)} \left( \sI \vee \sT_q \right)^{\alpha/(\alpha+2\beta s)}  + R^{\beta s} \lVert f_+ \rVert_{L^\infty(B(x_0,R))} \right], 
	\end{equation}
	where $\sI$ and  $\sT_q$  are defined as  \eqref{e:L2-Caccio}. 
	Additionally, if $u_+$ is  bounded in $M$, then 
	\begin{equation}\label{e:L2mean-diff-form-2}
		\esssup_{B(x_0, R/2)} u  \le C  \bigg[ \sI^{2\beta s/(\alpha+2\beta s)}\bigg( \sI \vee \Lambda(1-s) \esssup_{B(x_0,R/2)^c} u_+ \bigg)^{\alpha/(\alpha+2\beta s)}  + R^{\beta s} \lVert f_+ \rVert_{L^\infty(B(x_0,R))} \bigg],
	\end{equation}
	where $\Lambda> 0$ is the constant in \TJ.
\end{cor}
\pf Set $\nu':=2\beta s/(\alpha+2\beta s)$. For  $q \in [0,\beta]$, applying  Proposition \ref{p:L2mean} with $\delta=(\sI/\sT_q)^{\nu'}$, we get that 
\begin{align*}
	\esssup_{ B(x_0, R/2)} u \, &\le c_1   \big(  (1+\delta^{-1})^{ \alpha/(2\beta s)} \sI  +  \delta   \sT_q + R^{\beta s}\lVert f_+ \rVert_{L^\infty(B(x_0,R))} \big) \\[-2mm]
	&\le c_1R^{\beta s}\lVert f_+ \rVert_{L^\infty(B(x_0,R))} +  c_1 \begin{cases}
		2^{ \alpha/(2\beta s)} \sI  + \sI^{\nu'} \sT_q^{1-\nu'}	&\mbox{ if }  \delta \ge 1,\\[2mm]
		2^{ \alpha/(2\beta s)}  (\sT_q/\sI)^{\nu'  \alpha/(2\beta s)}\sI + \sI^{\nu'} \sT_q^{1-\nu'}		&\mbox{ if }   \delta < 1
	\end{cases}\\[1mm]
	&\le c_1R^{\beta s}\lVert f_+ \rVert_{L^\infty(B(x_0,R))} +  (2^{ \alpha/(2\beta s_0)}+1)c_1 \sI^{\nu'} \left(\sI \vee \sT_q \right)^{1-\nu'}.
\end{align*}
This proves \eqref{e:L2mean-diff-form-1}. For \eqref{e:L2mean-diff-form-2}, we assume that $K:= \esssup_{B(x_0,R/2)^c} u_+<\infty$ and take $q=\beta s$. By \TJ,   for all $r \in [R/2,R]$ and $a \in (0,R/8]$, we have
\begin{align*}
	a^{\beta s} \sT(u_+, B(x_0,r), B(x_0,r+a)^c) \le  Ka^{\beta s}\sup_{x \in B(x_0,r)} \int_{B(x,a)^c} J(x,dy)  \le \Lambda(1-s) K.
\end{align*}
Thus, $\sT_{\beta s}\le \Lambda (1-s)K$. From \eqref{e:L2mean-diff-form-1}, we conclude that \eqref{e:L2mean-diff-form-2} holds.	\qed

\section{Regularity estimates  for $(\sE,\sF)$}\label{section:WEHI-EHR}

In this section, we continue to assume that $(\sE,\sF)$ is a regular Dirichlet form on $L^2(M)$ with the representation \eqref{e:def-sE}. The goal of this section is to   establish a weak elliptic Harnack inequality (Proposition \ref{p:WEHI}) and elliptic H\"older regularity for $(\sE,\sF)$ (Corollary \ref{c:EHR}). For this, we mainly follow the framework of \cite{CKW-elp}.

\subsection{\TJ \ $+$ \FK \ $+$ \CS \ $+$ \PI \ $\Rightarrow$ \WEHI}

\smallskip

\begin{lemma}\label{l:WEHI-1+}
	Suppose that \PI \ holds with $K_1$.	There exists $C>0$ independent of $s$ such that the following holds:	Let $x_0 \in M$, $R\in (0,R_0)$ and   $r \in (0,R/(2K_1))$. Suppose that $u$  is  bounded, non-negative and  $-\sL u \ge f$ in $B(x_0,R)$ for  $f \in L^\infty(B(x_0,R))$. For given  $a,h>0$ and $b>1$, let
	\begin{align*}
		v=v(u,a,h,b):=\left[ \log \frac{a+h}{u+h}\right]_+ \wedge \log b.
	\end{align*}
	Then we have
	\begin{align*}
		&	 \int_{B(x_0,r)}  (v-\overline{v}_{B(x_0,r)})^2  \,d\mu \\
		&\le CV(x_0,r) \bigg( 1 + \frac{r^{\beta s}\,\NT(u_-,B(x_0,2K_1r), B(x_0,R)^c)+ r^{\beta s}\lVert f_- \rVert_{L^\infty(B(x_0,R))}}{h} \bigg).
	\end{align*}
\end{lemma}
\pf Following the proof of  \cite[Proposition 4.13]{CKW-jems} and noting that the constants appearing in the proof are independent of $s$, we obtain the result. \qed

Using Lemma \ref{l:WEHI-1+}, one can follow the proof of \cite[Lemma 3.3]{CKW-elp} with careful considerations of constants and obtain
\begin{lemma}\label{l:WEHI-2}
	Suppose that \PI \ holds with $K_1$. There exists $C>0$ independent of $s$ such that the following holds:	Let $x_0 \in M$,  $R\in (0,R_0)$ and $r \in (0,R/(4K_1))$. Suppose that $u$ is  bounded, non-negative and  $-\sL u \ge f$ in $B(x_0,R)$ for   $f \in L^\infty(B(x_0,R))$. 
	If there exist $a > 0$ and $\delta \in (0,1)$ satisfying
	\begin{align}\label{e:WEHI-2-1}
		\mu\left( B(x_0,r) \cap \{u < a\}\right) \le \delta V(x_0,r),
	\end{align}
	then for all $\eps \in (0,1)$,
	\begin{align*}
		\mu\left( B(x_0,2r) \cap \left\{ u \le  \eps  a-r^{\beta s} \big(  \NT(u_-, B(x_0,4K_1r), B(x_0,R)^c)  + \lVert f_-\rVert_{L^\infty(B(x_0,R))}\big) \right\} \right)
		\le \frac{CV(x_0,2r)}{(1-\delta) |\log \eps|}.
	\end{align*}
\end{lemma}

\begin{lemma}\label{l:WEHI-3}
	Suppose that \TJ, \FK, \CS \ and \PI \ hold with $K_1$. Let $x_0 \in M$,  $R\in (0,R_0)$ and  $r \in (0,R/(4K_1))$. Suppose that $u$ is  bounded, non-negative and  $-\sL u \ge f$ in $B(x_0,R)$ for  $f \in L^\infty(B(x_0,R))$. If \eqref{e:WEHI-2-1} holds for    $a>0$ and $\delta \in (0,1)$,  then there exists  $\eps_0=\eps_0(\delta) \in (0,1)$ depending on $\delta$ but  independent of $s,x_0,R,r,f,u$ and $a$ such that 
	\begin{align}\label{e:PLG-claim}
		\essinf_{B(x_0,r)} u \ge \eps_0 a -r^{\beta s} \left( \NT(u_-, B(x_0, 4K_1r), B(x_0,R)^c)+ \lVert f_-\rVert_{L^\infty(B(x_0,R))}  \right).
	\end{align} 
\end{lemma}
\pf    Define $B_l:=B(x_0,l)$ for $l>0$. Set
$$h:=r^{\beta s} \big(  \NT(u_-, B(x_0,4K_1r), B(x_0,R)^c)  + \lVert f_-\rVert_{L^\infty(B_R)}\big).$$ 
Let $k>0$ be a constant to be determined later and  $F:\R\to \R$ be a twice differentiable function with the following properties:
\begin{align}\label{e:condition-F-2}
	F '\le 0, \quad F''\ge 0, \quad \sup_{t\in \R} |F '(t)| +\sup_{t\in \R}  F''(t)<\infty,
\end{align}
\begin{equation}\label{e:condition-F-1}
	F(t)=\frac{1}{t+k}\quad \text{for all $t>-\frac{k}{2}$} \quad \text{ and } \quad 	F (t)\le -\frac{10t}{k^2}  \quad \text{for all $t\le -\frac{k}{2}$}.
\end{equation}
Since $-\sL u \ge f$ in $B_R$ and $(F' \circ u ) \le 0$, by Lemma  \ref{l:supersolution}, for any $0\le \phi \in \sF \cap C_c({B_R})$,
\begin{align*}
	\sE(F\circ u, \phi) \le \sE(u, (F' \circ u)\phi) \le \la (F' \circ u) f  \phi \ra,
\end{align*}
that is, $-\sL(F\circ u) \le(F' \circ u)f $ in $B_R$. 
Applying \eqref{e:L2mean-diff-form-1} to $F\circ u$, since $u \ge 0$ in $B_{r}$, we obtain
\begin{equation}\label{e:WEHI-eq-0}
	\Big(\essinf_{B_{r}} u + k \Big)^{-1} = \esssup_{B_{r}}F\circ u \le c_1\left[ \sI^{\nu_0} \left( \sI \vee  \sT_{\beta s}\right)^{1-\nu_0} + (2r)^{\beta s} \lVert [(F'\circ u)f]_+\rVert_{L^\infty(B_R)} \right], 
\end{equation}
where 
\begin{align*}
	\sI&:=\bigg(  \frac{1}{V(x_0,2r)}\int_{B_{2r}} (F\circ u)^2 d\mu  \bigg)^{1/2}, \\
	\sT_{\beta s}&:=\sup \left\{  b^{\beta s} \sT((F\circ u)_+, B(x_0,l), B(x_0,l+b)^c)\,: l \in [r,2r],\, b \in (0,r/4]\right\}
\end{align*} 
and $\nu_0:=2\beta s_0/( \alpha+2\beta s_0)$.  Since $u \ge 0$ in $B_R$, we have $-F'\circ u = (u+k)^{-2}\le  k^{-2}$ and $F\circ u = (u+k)^{-1}\le k^{-1}$  in $B_{R}$. Hence, $\sI \le k^{-1}$ and
\begin{equation}\label{e:WEHI-eq-f}
	(2r)^{\beta s}	\lVert [(F'\circ u)f]_+\rVert_{L^\infty(B_R)}    \le  2^{\beta} r^{\beta s} k^{-2} \lVert f_-\rVert_{L^\infty(B_R)}    \le 2^{\beta} hk^{-2}.
\end{equation}
Let $l\in [r,2r]$ and $b\in (0,r/4]$. For all  $x \in B_{l}$, using the fact that  $u \ge 0$ in $B_R$ and \TJ, we obtain
\begin{align*}
	b^{\beta s} \int_{B_{l+b}^c} F(u(y)) J(x,dy) &\le b^{\beta s}\bigg( \int_{ B_{l+b}^c \cap \{u> -k/2\}} \frac{J(x,dy)}{u(y)+k} + \frac{10}{k^2}\int_{B_{l+b}^c \cap \{u\le -k/2\}} u_-(y) J(x,dy)\bigg)\\
	&\le b^{\beta s}\bigg( \frac2k\int_{ B_{l+b}^c \cap \{u> -k/2\}}  J(x,dy) + \frac{10}{k^2}\int_{B_{R}^c \cap \{u\le -k/2\}} u_-(y) J(x,dy)\bigg)\\
	&\le b^{\beta s}\bigg( \frac2k\int_{B(x, b)^c}  J(x,dy) + \frac{10}{k^2} \sT(u_-, B_{2r}, B_R^c)\bigg)\le \frac{2\Lambda(1-s)}{k}  + \frac{10h}{k^2} .
\end{align*}
Thus, since $\sI\le k^{-1}$, we get
\begin{align}\label{e:WEHI-eq-2}
	\sI \vee \sT_{\beta s} \le  \frac{1+2\Lambda}{k}  + \frac{10h}{k^2}.
\end{align}
Let $\eps \in (0,1)$ be a constant  whose exact value  to be also  determined later. Since $u\ge 0$ in $B_{2r}$, by Lemma \ref{l:WEHI-2}, we have
\begin{align}\label{e:WEHI-eq-3}
	\sI^2 &= \frac{1}{V(x_0,2r)}\int_{B_{2r}  \cap \{ u > \eps a-h\}} \frac{\mu(dx)}{(u(x)+k)^2} + \frac{1}{V(x_0,2r)}\int_{B_{2r} \cap \{ u \le \eps a - h\}} \frac{\mu(dx)}{(u(x)+k)^2} \nn\\
	&\le \frac{1}{((\eps a-h)_+ +k)^2}  + \frac{\mu(B_{2r} \cap \{ u \le \eps a - h\})}{k^2 V(x_0,2r)} \le \frac{1}{((\eps a-h)_+ +k)^2} + \frac{c_2}{k^2(1-\delta)|\log \eps|}.
\end{align}
Combining \eqref{e:WEHI-eq-0}, \eqref{e:WEHI-eq-f}, \eqref{e:WEHI-eq-2} and \eqref{e:WEHI-eq-3}, and using the inequality $(a+b)^{-1} \ge (a^{-1} \wedge b^{-1})/2$ for $a,b>0$, we obtain
\begin{align}\label{e:WEHI-eq-4}
	&\essinf_{B_{r}} u  \ge \frac{1}{c_1}\left( \sI^{\nu_0} (\sI \vee \sT_{\beta s})^{1-\nu_0}  + 2^{\beta} hk^{-2} \right)^{-1} - k\nn\\
	& \ge \frac{1}{2c_1} \left[  \frac{ k^2 }{2^\beta h}\wedge  \bigg( \frac{1}{((\eps a-h)_+ +k)^2} + \frac{c_2}{k^2(1-\delta)|\log\eps|} \bigg)^{-\nu_0/2} \bigg( \frac{1+2\Lambda}{k} + \frac{10h}{k^2} \bigg)^{\nu_0-1}  \right] -k\nn\\
	& \ge  c_3 \left[ \frac{k^2}{h} \wedge \bigg( \frac{1}{((\eps a-h)_+ +k)^2} + \frac{1}{k^2(1-\delta)|\log\eps|} \bigg)^{-\nu_0/2} \bigg(\frac1k + \frac{h}{k^2} \bigg)^{\nu_0-1}  \right] - k,
\end{align}
where $c_3\in (0,1)$ is a constant independent of $s,x_0,R,r,f,u$ and $a$.

Now, we let
$
\gamma:= ( 2^{\nu_0/2-2}c_3)^{1/\nu_0} \in (0,1),
$
and take
\begin{align*}
	\eps = \exp\left( - \frac{1}{\gamma^2(1-\delta)}\right) \quad \text{ and } \quad k=\gamma \eps a.
\end{align*}
Set $\eps_0:= 2^{-1}c_3\gamma \eps$. 
If $h \ge \eps_0a$, then since $\essinf_{B_r}u \ge 0$,
 \eqref{e:PLG-claim} holds. Assume that $h<\eps_0a$. Then $c_3k^2/h = 2\eps_0 a k/h > 2k$. Moreover, since  $h<k<\eps a$, we have
\begin{align*}
	&c_3\bigg( \frac{1}{((\eps a-h)_+ +k)^2} + \frac{1}{k^2(1-\delta)|\log\eps|} \bigg)^{-\nu_0/2} \bigg(\frac1k + \frac{h}{k^2} \bigg)^{\nu_0-1} \\
	&= c_3k\bigg( \frac{k^2}{(\eps a-h+k)^2} + \frac{1}{(1-\delta)|\log\eps|} \bigg)^{-\nu_0/2} \bigg(1 + \frac{h}{k} \bigg)^{\nu_0-1}\\
	&\ge 2^{\nu_0-1}c_3k\bigg( \frac{(\gamma \eps a)^2}{(\eps a)^2} + \gamma^2 \bigg)^{-\nu_0/2} = 2^{\nu_0/2-1} c_3 k \gamma^{-\nu_0} = 2k.
\end{align*}
Consequently, we deduce from \eqref{e:WEHI-eq-4} that $\essinf_{B_{r}} u  \ge    2k-k =k \ge  \eps_0 a$, 
proving that \eqref{e:PLG-claim} holds.  The proof is complete. \qed

We  recall a Krylov-Safonov type covering lemma from \cite[Lemma 3.8]{CKW-elp}, which is originally due to \cite[Lemma 7.2]{KS01}.

\begin{lemma}\label{l:covering}
	Let $x_0 \in M$ and $r>0$. For a measurable set $E \subset B(x_0,r)$ and  $\eps \in (0,1)$, define
	\begin{equation}\label{e:[E]}
		[E]_\eps= \bigcup_{l \in (0,r)} \left\{ B(x,5l) \cap B(x_0,r): x \in B(x_0,r) \;\text{ and } \;\frac{\mu(E \cap B(x,5l))}{V(x,l)}>\eps \right\}. 
	\end{equation}
	Then we have either {\rm (1)} $ [E]_\eps = B(x_0,r) $ or {\rm (2)} $\mu([E]_\eps) \ge \eps^{-1}\mu(E).$
\end{lemma}

\begin{prop}\label{p:WEHI}
	Suppose that \TJ, \FK, \CS \ and \PI \ hold with $K_1$. There exist constants $\delta,C >0$  independent of $s$ such that for any $x_0 \in M$, $R \in (0,R_0)$, $r \in (0,R/(20K_1+3))$, and any Borel function $u$ that is bounded, non-negative and $-\sL u \ge f$ in $B(x_0,R)$ for  $f \in L^\infty(B(x_0,R))$, 
	\begin{align*}
		&\bigg(\frac{1}{V(x_0,r)}\int_{B(x_0,r)} u^\delta d\mu \bigg)^{1/\delta} \\
		&\le C \left[\, \essinf_{B(x_0,r)} u +  r^{\beta s}\Big(   \NT\left(u_-, B(x_0,(10K_1+1)r), B(x_0,R-2r)^c \right) + \lVert f_-\rVert_{L^\infty(B(x_0,R))}\Big) \,\right].
	\end{align*}
	Thus, \WEHI \ holds.
\end{prop}
\pf  Write $B_a:=B(x_0,a)$ for $a>0$ and let
$$h:= (5/2)^{\beta s}r^{\beta s}\left( \NT(u_-, B_{(10K_1+1)r}, B_{R-2r}^c) + \lVert f_-\rVert_{L^\infty(B_R)}\right).$$ By \ref{e:VD}, there exists $c_1 \in (0,1/2)$ such that
\begin{equation}\label{e:WEHI-VD}
	V(x,a) \ge 2c_1 V(x,5a) \quad \text{for all} \;\, x \in M \text{ and } a>0. 
\end{equation}
Let $\eps_0\in (0,1)$ be the constant from Lemma \ref{l:WEHI-3}, associated with the parameter $\delta = 1 - c_1$. 

Let $b>0$ be an arbitrary positive constant. Define for $n \ge 0$,
\begin{align*}
	E_n(b)&:= \left\{ x \in B_r: u(x)\ge b\eps_0^n - (1-\eps_0)^{-1}h\right\},\\
	H_n(b)&:=\left\{ (x,l) \in B_r \times (0,r) : \mu(E_n(b) \cap B(x,5l)) > 2^{-1} V(x,l)\right\}.
\end{align*}
We also define $[E_n(b)]_{2^{-1}}$ as \eqref{e:[E]}. By the definition, we have    
\begin{align}\label{e:WEHI-sets}
	[E_n(b)]_{2^{-1}}=\cup_{(x,l) \in H_n(b)} (B(x,5l) \cap B_r) \subset \cup_{(x,l)\in H_n(b)} B(x,5l).
\end{align}
We will prove at the end of this proof that
\begin{align}\label{e:WEHI-claim}
	\mu(E_{n+1}(b)) \ge \mu(	[E_n(b)]_{2^{-1}}) \quad \text{for all} \;\, n \ge 0.
\end{align}
Assume for the moment that \eqref{e:WEHI-claim} holds. Let $n_0=n_0(b) \in \N$ be such that
\begin{align}\label{e:WEHI-n0}
	2^{-n_0}V(x_0,r) < \mu(E_0(b)) \le 2^{-n_0+1}V(x_0,r).
\end{align}
Suppose that  $\mu(E_{n_0}(b)) <V(x_0,r)$.
Then by  \eqref{e:WEHI-claim} and Lemma \ref{l:covering},
\begin{align*}
	\mu(E_{k}(b)) \ge \mu ([E_{k-1}(b)]_{2^{-1}}) \ge 2 \mu (E_{k-1}(b)) \quad \text{for all} \;\, 1\le k \le n_0,
\end{align*}
and therefore, $	V(x_0,r)>	\mu(E_{n_0}(b))   \ge 2^{n_0}\mu(E_0(b)) >V(x_0,r).$
This gives a contradiction. Hence, $\mu(E_{n_0}(b)) = V(x_0,r)$. Using this in the first inequality below and  \eqref{e:WEHI-n0} in the second, we obtain
\begin{align}\label{e:WEHI-claim2}
	\essinf_{B_r} u&\ge b\eps_0^{n_0} - (1-\eps_0)^{-1}h \ge b \bigg( \frac{\mu(E_0(b))}{2V(x_0,r)}\bigg)^{|\log \eps_0|/ \log 2} -(1-\eps_0)^{-1}h.
\end{align}
Set $\delta_1:= (\log 2)/|\log \eps_0|$ and  
$$K:=\essinf_{B_{r}} u  + (1-\eps_0)^{-1}h.$$ By \eqref{e:WEHI-claim2}, we get that  for all $b>0$,
\begin{align}\label{e:WEHI-last}
	&\mu\left( \{x \in B_r: u(x) \ge b \}\right) \le  	\mu(E_0(b)) \le 2 (K/b)^{\delta_1}V(x_0,r) .
\end{align}
Using \eqref{e:WEHI-last}, we conclude that for any $\delta \in (0,\delta_1)$,
\begin{align*}
	&	\frac{1}{V(x_0,r)}\int_{B_r} u^\delta d\mu  \\
	&= \delta \int_0^\infty b^{\delta-1} \frac{\mu( \{x\in B_r:u(x)\ge b\})}{V(x_0,r)} db\le \delta \int_0^K b^{\delta-1} db + 2\delta K^{\delta_1} \int_K^\infty b^{\delta-1-\delta_1} db =    \frac{\delta_1+\delta}{\delta_1-\delta} K^\delta,
\end{align*}
which implies the desired result.

Now, we prove \eqref{e:WEHI-claim}. Let $n \ge 0$ and $\eps>0$. By the inner regularity of $\mu$ and  \eqref{e:WEHI-sets}, there exist a finite collection $\{B(x_i,5l_i)\}_{1\le i\le N}$  with $(x_i,l_i) \in H_n(b)$ such that
\begin{align}\label{e:cover-En(b)}
	\mu([E_n(b)]_{2^{-1}}) \le \mu( \cup_{i=1}^N B(x_i,5l_i)) + \eps.
\end{align}  
For all $1\le i \le N$, using  \eqref{e:WEHI-VD},   we see that
\begin{align*}
	&\mu\left( B(x_i,5l_i) \cap \{u < b\eps_0^n-(1-\eps_0)^{-1}h\}\right) \\
	&\le \mu\left( B(x_i,5l_i) \cap E_n^c\right)\le V(x_i,5l_i)- 	2^{-1} V(x_i,l_i) \le  (1-c_1) V(x_i,5l_i).
\end{align*}
Hence, applying Lemma \ref{l:WEHI-3},   we obtain
\begin{align*}
	\essinf_{B(x_i,5l_i)} u &\ge \eps_0 ( b\eps_0^n - (1-\eps_0)^{-1})h  - (5l_i)^{\beta s} ( \NT(u_-, B(x_i, 20K_1 l_i), B(x_i, R-r)^c)  + \lVert f_-\rVert_{L^\infty(B_R)})  \\
	&\ge b\eps_0^{n+1} - \eps_0(1-\eps_0)^{-1}h - (5/2)^{\beta s} r^{\beta s} (  \NT(u_-, B_{(10K_1+1)r}, B_{R-2r}^c) + \lVert f_-\rVert_{L^\infty(B_R)}) \\
	&\ge b\eps_0^{n+1} - \eps_0(1-\eps_0)^{-1}h - h= b\eps_0^{n+1} - (1-\eps_0)^{-1}h.
\end{align*}It follows that  $\mu(\cup_{i=1}^NB(x_i,5l_i) \setminus E_{n+1}(b))=0$. Combining this with \eqref{e:cover-En(b)}, we arrive at
\begin{align*}
	&\mu(	[E_n(b)]_{2^{-1}}) 	\le \mu(\cup_{i=1}^N B(x_i, 5l_i) ) +\eps\le \mu(E_{n+1}(b)) + \eps.
\end{align*}
Since $\eps>0$ is an arbitrary positive constant,   \eqref{e:WEHI-claim} holds true.
The proof is complete. \qed 

\subsection{\TJ  \ $+$ \WEHI \ $\Rightarrow$ \EHR}

\smallskip

For an open set $D \subset M$ and $u \in L^\infty(D)$, define
$$\essosc_D u:=\esssup_D u - \essinf_D u.$$

By a standard argument,  we establish the next oscillation inequality for weak solutions of the Poisson equation \eqref{e:Poisson} under \WEHI \ and \TJ.

\begin{prop}\label{p:EHR-osc}
	Suppose that  \WEHI \ and \TJ \ hold. There exist constants $\gamma \in (0,1)$ and $C>0$ independent of $s$ such that for any $x_0 \in M$,  $0<r\le R<R_0$, and any Borel function $u$ that is bounded in $M$ and $-\sL u = f$ in $B(x_0,R)$ for  $f \in L^\infty(B(x_0,R))$,
	\begin{align*}
		\essosc_{B(x_0,r)} u \le C  \left((r/R)^\gamma  \lVert u \rVert_{L^\infty(M)}  + R^{\beta s }\lVert f \rVert_{L^\infty(B(x_0,R))}\right).
	\end{align*} 
\end{prop} 
\pf Let $x_0 \in M$, $R \in (0,R_0)$ and $f\in L^\infty(B(x_0,R))$. Write $B_a:=B(x_0,a)$ for $a>0$. Let $K\ge 1$ be  the constant in  \WEHI \  and set $\eps:=1/(3K+2)$.  Define 
$$W_{-1}:=M \quad \text{and}  \quad W_n:=B_{\eps^n R} \text{ \;for $n \ge 0$}.$$
Let  $\gamma \in (0, \beta s_0/2)$  be a constant  to be determined later. In the following, we construct a non-increasing sequence $(b_n)_{n \ge -1}$ and a non-decreasing sequence $(a_n)_{n \ge -1}$ with the following properties: For all $n \ge -1$,
\begin{equation}\label{e:EHR-construct}
	\esssup_{W_n} u \le b_n +R^{\beta s} \lVert f \rVert_{L^\infty(B_R)}\sum_{0\le i \le n}\eps^{i\beta s}  , \qquad \essinf_{W_n} u \ge a_n -R^{\beta s} \lVert f \rVert_{L^\infty(B_R)} \sum_{0\le i \le n}\eps^{i\beta s}  
\end{equation}
and
\begin{equation}\label{e:EHR-construct-2}
	b_n-a_n= 2\eps^{n\gamma}\lVert u \rVert_{L^\infty(M)}.
\end{equation} 
Define $b_{-1}=b_0=\lVert u \rVert_{L^\infty(M)}$, $a_{-1}= -(2\eps^{-\gamma}-1)\lVert u \rVert_{L^\infty(M)}$ and  $a_0=-\lVert u \rVert_{L^\infty(M)}$.  Then  \eqref{e:EHR-construct} and \eqref{e:EHR-construct-2} hold for $n=-1,0$.

Suppose that monotone sequences  $(b_n)_{-1 \le n \le k}$ and $(a_n)_{-1 \le n \le k}$ are constructed to satisfy \eqref{e:EHR-construct} and \eqref{e:EHR-construct-2} for all $0\le n \le k$, for some $k \ge 0$. Consider the function
\begin{align*}
	v_k(x):=\eps^{-k\gamma}\bigg(u(x) - \frac{a_k+b_k}{2} \bigg).
\end{align*}
For all $0 \le j \le k+1$, using the induction hypothesis, we see that
\begin{align}\label{e:EHR-1}
	\esssup_{W_{k-j}} v_k &\le \eps^{-k\gamma }\bigg(b_{k-j} -  \frac{a_k+b_k}{2} \bigg)  + \eps^{-k\gamma} R^{\beta s} \lVert f \rVert_{L^\infty(B_R)}\sum_{0\le i \le k-j}\eps^{i\beta s}   \nn\\
	&\le \eps^{-k\gamma }\bigg(b_{k-j} - a_{k-j} -\frac{b_{k}-a_k}{2} \bigg) + \eps^{-k\gamma} R^{\beta s} \lVert f \rVert_{L^\infty(B_R)}\sum_{0\le i \le k-j}\eps^{i\beta s}  \nn\\
	& = (2\eps^{-j\gamma} - 1) \lVert u \rVert_{L^\infty(M)}+ \eps^{-k\gamma} R^{\beta s} \lVert f \rVert_{L^\infty(B_R)}\sum_{0\le i \le k-j}\eps^{i\beta s}  
\end{align}
and
\begin{align}\label{e:EHR-2}
	\essinf_{W_{k-j}} v_k& \ge - \eps^{-k\gamma }\bigg(-a_{k-j}+\frac{a_{k}+b_k}{2} \bigg)  - \eps^{-k\gamma} R^{\beta s} \lVert f \rVert_{L^\infty(B_R)}\sum_{0\le i \le k-j}\eps^{i\beta s}  \nn\\
	& \ge -(2\eps^{-j\gamma} - 1) \lVert u \rVert_{L^\infty(M)}- \eps^{-k\gamma} R^{\beta s} \lVert f \rVert_{L^\infty(B_R)}\sum_{0\le i \le k-j}\eps^{i\beta s}  .
\end{align}
There are two cases.

\smallskip

Case 1: Suppose that $\mu(W_{k+1} \cap \{v_k \ge 0\}) \ge 2^{-1} \mu(W_{k+1})$. Consider 
$$w_k:=v_k +  \lVert u \rVert_{L^\infty(M)} +\eps^{-k\gamma} R^{\beta s} \lVert f \rVert_{L^\infty(B_R)} \sum_{0\le i \le k}\eps^{i\beta s} .$$
By \eqref{e:EHR-2}, $\essinf_{W_{k}} w_k  \ge 0$.  Further, we have $-\sL w_k \ge -\eps^{-k\gamma} |f|$ in $W_k$.  Applying \WEHI \  (with $R$ replaced by $\eps^k R$ and $r$ replaced by $\eps^{k+1}R$) to $w_k$, we obtain
\begin{align}\label{e:EHR-onestep}
	&	\bigg(\frac{1}{\mu(W_{k+1})}\int_{W_{k+1}} w_k^\delta \, d\mu \bigg)^{1/\delta}\nn\\
	&\le c_1\left( \essinf_{W_{k+1}} w_k + (\eps^{k+1}R)  ^{\beta s}\Big[ \NT\big((w_k)_-, B_{K\eps^{k+1}R}, B_{(1- 2\eps)\eps^kR}^c \big) +  \eps^{-k\gamma}\lVert f\rVert_{L^\infty(B_R)}\Big] \right), 
\end{align}
where $c_1$ and $\delta$ are positive constants independent of $s,x_0,r,R,f,u$ and $k$.   
By the assumption $\mu(W_{k+1} \cap \{v_k \ge 0\}) \ge 2^{-1} \mu(W_{k+1})$, we have
\begin{equation}\label{e:EHR-main1}
	\bigg(\frac{1}{\mu(W_{k+1})}\int_{W_{k+1}} w_k^\delta \, d\mu \bigg)^{1/\delta} \ge \bigg(\frac{1}{\mu(W_{k+1})}\int_{W_{k+1} \cap \{v_k \ge 0\}} w_k^\delta \, d\mu \bigg)^{1/\delta}\ge  2^{-1/\delta}\lVert u \rVert_{L^\infty(M)}.
\end{equation}
Moreover, since $K\eps <1/3$ and  $\essinf_{W_k}w_k\ge 0$, we  see that
\begin{align}\label{e:EHR-main2}
&	\NT\big((w_k)_-, B_{K\eps^{k+1}R}, B_{(1- 2\eps)\eps^kR}^c \big) \nn\\
	&\le  \NT\big((w_k)_-, B_{\eps^{k}R/3}, W_k^c \big)= \sup_{x \in B_{\eps^k R/3}} \sum_{j=1}^{k+1} \int_{W_{k-j} \setminus W_{k-j+1}} (w_k)_-(y) J(x,dy).
\end{align}
Let  $x \in B_{\eps^k R/3}$ and $1\le j \le k+1$. Using \eqref{e:EHR-2} and \TJ, we see that
\begin{align*}
	\int_{W_{k-j} \setminus W_{k-j+1}}  (w_k)_-(y)J(x,dy) &\le 2 (\eps^{-j\gamma} -1)\lVert u \rVert_{L^\infty(M)}  J(x,W_{k-j+1}^c) \\
	&\le 2 (\eps^{-j\gamma} -1)\lVert u \rVert_{L^\infty(M)}  J(x,B(x,\,2\eps^{k-j+1}R/3)^c)\\
	&\le \frac{2\Lambda(1-s) (\eps^{-j\gamma} -1)\lVert u \rVert_{L^\infty(M)}}{(2\eps^{k-j+1}R/3)^{\beta s}}.
\end{align*}
Hence, from \eqref{e:EHR-main2}, using  the inequality $a^b-1\le ba^b\log a$ for $a> 1$ and $b>0$ in the second inequality below,  and $\gamma<\beta s_0/2$ and  $\sup_{0<a<1}a^{\beta s_0/2} |\log a|<\infty$ in the fourth, we obtain
\begin{align*}
&	(\eps^{k+1}R)^{\beta s} \,\NT\big((w_k)_-, B_{K\eps^{k+1}R}, B_{(1- 2\eps)\eps^kR}^c \big)\le \frac{2\Lambda (1-s)\lVert u \rVert_{L^\infty(M)}}{(2/3)^\beta} \sum_{j=1}^{k+1}  \eps^{j\beta s}(\eps^{-j\gamma} -1) \\
	&\le \frac{2\gamma\Lambda (1-s) |\log \eps|\lVert u \rVert_{L^\infty(M)}}{(2/3)^\beta} \sum_{j=1}^{\infty} j\eps^{j(\beta s-\gamma)}   \le  \frac{2\gamma\Lambda   \eps^{\beta s_0/2}|\log \eps|\lVert u \rVert_{L^\infty(M)} }{(2/3)^\beta (1- \eps^{\beta s_0/2})^2}.
\end{align*} Combining this with \eqref{e:EHR-onestep} and  \eqref{e:EHR-main1}, we arrive at
\begin{equation}\label{e:EHR-main6}
	\essinf_{W_{k+1}} w_k\ge 	\bigg(2^{-1/\delta}c_1^{-1} -  \frac{2\gamma\Lambda   \eps^{\beta s_0/2}|\log \eps|}{(2/3)^\beta (1- \eps^{\beta s_0/2})^2}\bigg) \lVert u \rVert_{L^\infty(M)}  - \eps^{(k+1)\beta s -k\gamma} R^{\beta s} \lVert f \rVert_{L^\infty(B_R)}.
\end{equation}
Now we choose $\gamma \in (0,\beta s_0/2)$ small enough to satisfy
\begin{align*}
	2^{-1/\delta}c_1^{-1} -  \frac{2\gamma\Lambda   \eps^{\beta s_0/2}|\log \eps|}{(2/3)^\beta (1- \eps^{\beta s_0/2})^2} \ge  2^{-2/\delta} c_1^{-1} \quad \text{and} \quad 2\eps^\gamma \ge 2-2^{-2/\delta} c_1^{-1}.
\end{align*}
Then, by  \eqref{e:EHR-main6}, since $a_k=b_k - 2\eps^{k\gamma} \lVert u \rVert_{L^\infty(M)}$ by the induction hypothesis, we deduce that
\begin{align*}
	&\essinf_{W_{k+1}} u = \frac{a_k+b_k}{2} + \eps^{k\gamma }\essinf_{W_{k+1}} \bigg[w_k- \lVert u \rVert_{L^\infty(M)} -\eps^{-k\gamma} R^{\beta s} \lVert f \rVert_{L^\infty(B_R)} \sum_{0\le i \le k}\eps^{i\beta s}\bigg]\\
	&\ge \frac{a_k+b_k}{2} +  \eps^{k\gamma} \bigg[	( 2^{-2/\delta}c_1^{-1} -1) \lVert u \rVert_{L^\infty(M)} - \eps^{-k\gamma} R^{\beta s} \bigg(   \eps^{(k+1)\beta s}  +   \sum_{0\le i \le k}\eps^{i\beta s} \bigg) \lVert f \rVert_{L^\infty(B_R)}    \bigg] \\
	&\ge b_k   - \eps^{k\gamma} (2- 2^{-2/\delta}c_1^{-1})\lVert u \rVert_{L^\infty(M)} - R^{\beta s} \lVert f \rVert_{L^\infty(B_R)}   \sum_{0\le i \le k+1}\eps^{i\beta s}\\
	&\ge b_k - 2\eps^{(k+1)\gamma} \lVert u \rVert_{L^\infty(M)}   - R^{\beta s} \lVert f \rVert_{L^\infty(B_R)}   \sum_{0\le i \le k+1}\eps^{i\beta s}.
\end{align*}
By letting $b_{k+1}=b_{k}$ and  $a_{k+1}=b_{k}-2\eps^{(k+1)\gamma} \lVert u \rVert_{L^\infty(M)}  $, we conclude that  \eqref{e:EHR-construct} and \eqref{e:EHR-construct-2} hold for $k+1$.

Case 2: Suppose that  $\mu(W_{k+1} \cap \{v_k\ge 0\}) < 2^{-1} \mu(W_{k+1})$. Define
$$
\wt w_k = -v_k+\lVert u \rVert_{L^\infty(M)} +\eps^{-k\gamma} R^{\beta s} \lVert f \rVert_{L^\infty(B_R)} \sum_{0\le i \le k}\eps^{i\beta s}.
$$
Following the argument for Case 1, using the function $\wt w_k$ instead of $w_k$ and the inequality  \eqref{e:EHR-1} instead of \eqref{e:EHR-2}, one can deduce that  \eqref{e:EHR-construct} and \eqref{e:EHR-construct-2} hold for $k+1$ with $a_{k+1}=a_k$ and  $b_{k+1}=a_{k} + 2\eps^{(k+1)\gamma} \lVert u \rVert_{L^\infty(M)} $.

\smallskip

By induction, we get monotone sequences $(b_n)_{n \ge -1}$ and $(a_n)_{n \ge -1}$ satisfying \eqref{e:EHR-construct} and \eqref{e:EHR-construct-2} for all $n \ge -1$. Consider any $r \in (0,R]$, and  let $n_0 \ge 0$ be such that $\eps^{n_0+1}<r/R\le \eps^{n_0}$. By  \eqref{e:EHR-construct} and \eqref{e:EHR-construct-2}, we arrive at
\begin{align*}
&	\essosc_{B_r} u \le 	\essosc_{W_{n_0}} u \le  b_{n_0} - a_{n_0}  + 2R^{\beta s} \lVert f \rVert_{L^\infty(B_R)} \sum_{i \ge 0} \eps^{i\beta s}\\
	&\le 2\eps^{n_0\gamma} \lVert u \rVert_{L^\infty(M)}  + \frac{2R^{\beta s}}{1-\eps^{\beta s}} \lVert f \rVert_{L^\infty(B_R)}  \le 2\eps^{-\gamma}\bigg(\frac{r}{R}\bigg)^\gamma \lVert u \rVert_{L^\infty(M)} + \frac{2R^{\beta s}}{1-\eps^{\beta s_0}} \lVert f \rVert_{L^\infty(B_R)}.
\end{align*}
The proof is complete. \qed

\begin{cor}\label{c:EHR}
	Suppose that \WEHI \ and \TJ \  hold. Then \EHR \ holds.
\end{cor} 
\pf   For $\mu$-a.e. $x,y \in B(x_0,R/4)$,  we get from Proposition \ref{p:EHR-osc} that 
\begin{align*}
	|u(x)-u(y)| \le 	\essosc_{B(x,3d(x,y)/2)} u \le c_1 \bigg(\frac{3d(x,y)/2}{3R/4}\bigg)^\gamma\lVert u \rVert_{L^\infty(M)},
\end{align*} 
where $\gamma>0$ is the constant in Proposition \ref{p:EHR-osc}.  \qed

\smallskip 

\noindent \textbf{Proof of Theorem \ref{t:main-2}.} The result follows from Theorem \ref{t:main-1}, Proposition \ref{p:WEHI} and Corollary \ref{c:EHR}. \qed

\section{Proof  of Theorem \ref{t:coercivity}}\label{section:Thm1}

In this section, we give the proof of Theorem \ref{t:coercivity}.
 Throughout this section, we assume that $J(dx,dy)=J(x,y)\mu(dx)\mu(dy)$ in $M\times M$ and  \eqref{e:ass-coercivity} holds with $R_0,\delta_0,\sigma$ and $\theta$. By \ref{e:VD}, there exists $\eps_0 \in (0,1)$ such that
 \begin{align}\label{e:VD-5}
 	V(x,R)\le \eps_0 V(x, \frac{\delta_0}{5}R) \quad \text{for all $x \in M$ and $R>0$.}
 \end{align}
 
  We begin with two lemmas about  geometric properties of  $(M,d,\mu)$.

 \begin{lemma}\label{l:coercivity-0}
 	Let $x \in M$, $R>0$ and $E$ be a measurable set. If
 	\begin{align*}
 		\mu\left(E \cap B(x,R)\right) \ge (1-\eps_0^2 a) V(x,R)
 	\end{align*}
 	for some $a\in (0,1)$, then for any $r \in (0,(1-\delta_0/5)R]$, there exists $ z \in  B(x,R-r)$ such that
 	\begin{align*}
 		\mu(E \cap B(z,r)) \ge (1- a)V(z,r).
 	\end{align*}
 \end{lemma}
 \pf Suppose that $\mu(E \cap B(z,r))<(1-a)V(z,r)$ for all $B(z,r ) \subset B(x_0,R)$. By the Vitali covering lemma, there exists a collection $\{B(z_i,r)\}_{i=1}^\infty$ of pairwise disjoint open balls with $z_i \in B(x_0,R-r)$ such that $B(x_0,R-r) \subset \cup_{i=1}^\infty B(z_i,5r)$. By assumption, $\mu\left(E^c \cap B(z_i,r)\right) > aV(z_i,r)$ for all $i \ge 1$. Using this and \eqref{e:VD-5}, we get
 \begin{align*}
 	&	\mu\left(E^c \cap B(x_0,R)\right) \ge \sum_{i=1}^\infty	\mu\left(E^c \cap B(z_i,r)\right)\\
 	& > a\sum_{i=1}^\infty V(z_i,r)  \ge \eps_0 a\sum_{i=1}^\infty V(z_i,5r) \ge \eps_0 a V(x_0,R-r) \ge \eps_0^2 a V(x_0,R),
 \end{align*}
which contradicts the assumption that $	\mu\left(E \cap B(x,R)\right) \ge (1-\eps_0^2 a) V(x,R)$. \qed

 \begin{lemma}\label{l:coercivity-0+}
 	Let $E$ be a measurable set. The function
 	\begin{align*}
 	(x,r) \mapsto \frac{\mu(E\cap B(x,r))}{V(x,r)}
 	\end{align*}
 	is jointly continuous in $M \times (0,R_0)$.
 \end{lemma}
 \pf Let $x \in M$ and $r\in (0,R_0)$. For any $(y,r') \in B(x,r/4) \times (r/2,R_0)$, we have
 \begin{align*}
 	&\left| \frac{\mu(E\cap B(x,r))}{V(x,r)} -  \frac{\mu(E\cap B(y,r'))}{V(y,r')} \right|  \nn\\
 	&\le \frac{ V(y,r')\left| \mu(E\cap B(x,r)) - \mu(E\cap B(y,r')) \right|   + \mu(E\cap B(y,r')) \left| V(y,r')  - V(x,r) \right|  }{V(x,r)V(y,r')}\nn\\
 		&\le \frac{ \mu(E \cap  B(x, r\vee r' + d(x,y)))  - \mu\left(E \cap  B(x, r\wedge r'-d(x,y) ) \right)   }{V(x,r)}\nn\\
 		&\quad + \frac{ V(x,r\vee r' + d(x,y)) - V(x, r\wedge r'-d(x,y) ) }{V(x,r)} \nn\\
 		&\le \frac{ 2(V(x,r\vee r' + d(x,y)) - V(x, r\wedge r'-d(x,y) )) }{V(x,r)} .
 \end{align*}
 Since $\lim_{(y,r') \to (x,r)} (V(x,r\vee r' + d(x,y)) - V(x, r\wedge r'-d(x,y) )) =  \mu(\{z \in M: d(x,z)=r\})=0$ by the outer regularity of $\mu$ and \eqref{e:surface}, the result follows.
  \qed

 We  define $J_{k}(x,y)$ and $N_{k}(x)$ inductively as follows.  Let $J_0(x,y):=J(x,y)$ and
\begin{align*}
	N_0(x):= \left\{ z \in M  : J_0(x,z) \ge \frac{\theta(1-s)}{ V(x,d(x,z)) d(x,z)^{\beta s}}  \right\}.
\end{align*}   Let   $\lambda\in (0,1)$ be a constant   to be  determined later. Suppose that $J_k(x,y)$ and  $N_k(x)$ are well-defined  for some $k\ge 0$. Set  
\begin{align*}
&	r_k (x,z):= \sup\bigg\{  r \ge 0:  r\le \frac{\delta_0d(x,z)}{10},\, \exists w \in B(z,r) \text{ such that } \frac{\mu( N_k(x)\cap B(w,r))}{V(w,r)} \ge 1- \frac{\eps_0^2\sigma}{2} \bigg\},\\
	&	\sJ_k(x,y,z):=\begin{cases}
		\displaystyle	\frac{J_k(x,z)}{	V(z, r_k(x,z))} \wedge \frac{d(y,z)^{\beta s}J(y,z)}{V(y, d(x,y))  d(x,y)^{\beta s}} &\mbox{if  $r_k(x,z)>0$ and $y \in B(z, 4r_k(x,z))$},\\[10pt]
		\displaystyle 0 &\mbox{otherwise} 
	\end{cases}.
\end{align*}
Here $\eps_0 \in (0,1)$ is the constant in \eqref{e:VD-5}. Then we define
\begin{align*}
		J_{k+1} (x,y)& = \int_{B(y, 2d(x,y))}  \sJ_k(x,y,z) \mu(dz),\nn\\
 N_{k+1}(x)&=\left\{ z \in M  : J_{k+1}(x,z) \ge \frac{\theta \lambda^k(1-s)}{ V(x,d(x,z)) d(x,z)^{\beta s}}  \right\}.
\end{align*}

\begin{lemma}\label{l:coercivity-1}
	For any $k\ge 0$, there exists $C>0$ depending on $k,s_0$ and the constants in {\rm \ref{e:VD}} only such that for any ball $B(x_0,R)$ and a function $f:B(x_0,R)\to \R$,
	\begin{align}\label{e:coercivity-1-1}
	&	\int_{B(x_0,R)\times B(x_0,R)} (f(x)-f(y))^2 J(x,y)\mu(dx)\mu(dy) \nn\\
	&\ge C 	\int_{B(x_0,5^{-k}R)\times B(x_0,5^{-k}R)} (f(x)-f(y))^2 J_k(x,y)\mu(dx)\mu(dy).
	\end{align}
\end{lemma}
\pf  Write $aB:=B(x_0,aR)$ for $a>0$. Clearly, \eqref{e:coercivity-1-1} holds for $k=0$. Suppose that \eqref{e:coercivity-1-1} holds for $k-1$. We have
\begin{align*}
&	\int_{5^{-k}B\times 5^{-k}B} (f(x)-f(y))^2 J_k(x,y)\mu(dx)\mu(dy) \nn\\
&\le 2	\int_{5^{-k}B\times 5^{-k}B}  \int_{B(y, 2d(x,y))} \left( (f(x)-f(z))^2 +  (f(y)-f(z))^2 \right) \sJ_{k-1}(x,y,z) \mu(dz) \mu(dx)\mu(dy)\nn\\
&\le 2	\int_{5^{-k}B\times 5^{-k}B} \int_{B(y,2d(x,y))} (f(x)-f(z))^2\frac{  J_{k-1}(x,z)}{	V(z, r_{k-1}(x,z))} \1_{B(z, 4r_{k-1}(x,z))}(y) \, \mu(dz)\mu(dx)\mu(dy) \nn\\
&\quad + 2 	\int_{5^{-k}B\times 5^{-k}B}  \int_{B(y, 2d(x,y))}  (f(y)-f(z))^2 \frac{d(y,z)^{\beta s}J(y,z)}{V(y, d(x,y))  d(x,y)^{\beta s}} \mu(dz) \mu(dx)\mu(dy) \nn\\
&=:2(I_1+I_2).
\end{align*}
Using  Fubini's theorem, \ref{e:VD} and the induction hypothesis, we obtain 
\begin{align}\label{e:coercivity-1-2}
	I_1 &\le \int_{5^{-k}B\times 5^{1-k}B}  (f(x)-f(z))^2  J_{k-1}(x,z)  \frac{  V(z, 4r_{k-1}(x,z))}{	V(z, r_{k-1}(x,z))}\mu(dx)\mu(dz) \nn\\
	&\le c_1 \int_{5^{1-k}B\times 5^{1-k}B}  (f(x)-f(z))^2  J_{k-1}(x,z)\mu(dx)\mu(dz)  \nn\\
	&\le c_2\int_{B\times B}  (f(x)-f(z))^2  J(x,z)\mu(dx)\mu(dz).
\end{align}
Further, using Fubini's theorem and Lemma \ref{l:robust-jumptail}, we  get
\begin{align*}
I_2&\le	\int_{5^{-k}B\times 5^{1-k}B}  (f(y)-f(z))^2 d(y,z)^{\beta s}J(y,z) \int_{B(y, d(y,z)/2)^c}   \frac{\mu(dx)}{V(y, d(x,y))d(x,y)^{\beta s}}  \mu(dy)\mu(dz)\nn\\
&\le \frac{c_3}{\beta s_0}\int_{5^{-k}B\times 5^{1-k}B}  (f(y)-f(z))^2 J(y,z) \mu(dy)\mu(dz).
\end{align*}
Combining this with \eqref{e:coercivity-1-2}, we arrive at the result by induction.  \qed 

For the proof of the next lemma, we adapt the argument from \cite[Lemma 4.2]{CS20}. However, there was a gap in the last display of their proof: in their notation, the set $\Xi$ may not be contained in $\Omega_j(x,y)$. Using Lemma \ref{l:coercivity-0}, we address this issue.

\begin{lemma}\label{l:coercivity-3}
		There exists $c_0 \in (0,1]$ depending on $\delta_0,\sigma$ and the constants in {\rm \ref{e:VD}}  only such that if $\lambda \le c_0$, then for all   $x \in M$ and $k\ge 0$,
		\begin{align*}
			\left\{ z \in B(x,2R_0): r_k(x,z)>0\right\} \subset    N_{k+1}(x).
		\end{align*}
\end{lemma}
\pf    Let $y \in B(x,2R_0)$ be such that $r_k(x,y)>0$. We prove that there exists $c_0 \in (0,1]$ depending on $\delta_0,\sigma$ and the constants in \ref{e:VD}  only such that 
\begin{align}\label{e:coercivity-3-claim}
	J_{k+1}(x,y) \ge \frac{c_0\theta \lambda^{k} (1-s)}{V(x, d(x,y))  d(x,y)^{\beta s} }.
\end{align}
This yields the desired result.

Set $v_0:=y$. Since $r_k(x,v_0)>0$, there exists $w_0 \in B(v_0, r_k(x,v_0))$ such that $	\mu( N_k(x) \cap  B(w_0, r_k(x,v_0)) ) \ge (1-\eps_0^2\sigma/2) V(w_0, r_k(x,v_0))$. We construct sequences $(v_j)_{j \ge 0}$ and $(w_j)_{j\ge 0}$ as follows:  Set
$$
\eps:= \frac{3\delta_0-\delta_0^2}{5}\quad \text{and} \quad K:=1 + \frac{2}{\eps}.
$$
 If  there exists $v_j \in B(w_{j-1},  r_k(x,v_{j-1}))$ with $r_k(x,v_j)>Kr_k(x,v_{j-1}),$
 select  such $v_j$. Since $r_k(x,v_j)>0$, there is $w_j \in B(v_j, r_k(x,v_j))$ such that   
\begin{align}\label{e:def-w-j}
	\frac{\mu( N_k(x) \cap B(w_j,   r_k(x,v_j)))}{V(w_j, r_k(x,v_j)) } \ge 1-\frac{\eps_0^2\sigma}{2}.
\end{align} 
If $v_j$ is well-defined,  then  
\begin{align}\label{e:y-v-j}
	d(y,v_j) &\le d(y, v_{j-1}) + d(v_{j-1}, w_{j-1})+ d(w_{j-1},v_j)<  d(y, v_{j-1})  + 2r_k(x,v_{j-1}) \nn\\
	&   \le \cdots \le 2\sum_{i=0}^{j-1} r_k(x,v_{i}) \le 2r_k(x, v_{j}) \sum_{i=1}^\infty K^{-i}= \eps r_k(x,v_j).
\end{align}
 In particular, we have $10d(y,v_j) < 10r_k(x,v_j) < d(x,v_j) \le  d(y,v_j) + d(x,y)$ which implies 
 \begin{align}\label{e:x-v-j} 
 	10r_k(x,v_j) < d(x,v_j)  < \frac{10}{9}d(x,y).
 \end{align}
 Since $r_k(x,v_j)>K^j r_k(x,v_0)$ whenever $v_j$ is well-defined,  the iterative construction   stops after $n$  steps for some $n<\infty$.

 By \eqref{e:def-w-j} and Lemma \ref{l:coercivity-0}, there exists $w_n' \in B(w_n, (\delta_0/5) r_k(x,v_n))$ such that 
 \begin{align}\label{e:def-w-n}
 \frac{\mu( N_k(x) \cap B(w_n',  (1-\delta_0/5) r_k(x,v_n)))}{V(w_n', (1-\delta_0/5) r_k(x,v_n)) }\ge 1- \frac{\sigma}{2}.
 \end{align} 
   By \eqref{e:y-v-j},   we have
 \begin{align}\label{e:y-w-j}
 	d(y, w_n') &\le d(y,v_n) + d(v_n,w_n) + d(w_n,w_n')\nn\\
 & < \left(1+ \eps + \delta_0/5\right)r_k(x,v_n)  = (1+\delta_0) \left(1- \delta_0/5\right) r_k(x,v_n).
 \end{align} 
Further, by \eqref{e:x-v-j}, $(1-\delta_0/5)r_k(x,v_n)<R_0$.  Set $A:=  B(w_n',  (1-\delta_0/5) r_k(x,v_n))$. Using \eqref{e:def-w-n}, \eqref{e:y-w-j} and \eqref{e:ass-coercivity}, we get
\begin{align}\label{e:coercivity-3-1}
		\mu \left(   N_k(x) \cap N_0(y)  \cap A \right) 	& \ge 	\mu \left(  N_k(x) \cap A \right)  +  	\mu \left(    N_0(y) \cap A \right)  -  \mu(A) \ge \frac{\sigma}{2} \mu(A).
\end{align} 
For any $z \in A$,  we have  $r_k(x,z) \le Kr_k(x,v_n)$ by the maximality of $n$,
\begin{align*}
	d(x,z) \le d(x,v_n) + d(v_n,w_n) + d(w_n,z) \le d(x,v_n) + 2r_k(x,v_n) \le  \frac{11d(x,y)}{9}
\end{align*}
by \eqref{e:x-v-j} and 
\begin{align}\label{e:dist-y-z}
	d(y,z) \le d(y,w_n') + d(w_n', z) <  (2+\delta_0)\left(1- \frac{\delta_0}{5}\right) r_k(x,v_n) 
\end{align}
by \eqref{e:y-w-j}. 
Thus, for any $z \in A \cap N_k(x) \cap N_0(y)$, using \ref{e:VD2}, we get
\begin{align}\label{e:coercivity-3-2}
	\frac{J_k(x,z)}{	V(z, r_k(x,z))} &\ge 	\frac{\theta \lambda^k (1-s)}{ V(x, d(x,z))d(x,z)^{\beta s}	V(z, Kr_k(x,v_n))}\nn\\
	&\ge 	\frac{c_1\theta \lambda^k (1-s)}{ V(x, d(x,y))d(x,y)^{\beta s}	V(w_n,r_k(x,v_n))}
\end{align}
and
\begin{align}\label{e:coercivity-3-3}
	\frac{d(y,z)^{\beta s}J(y,z)}{V(y, d(x,y))  d(x,y)^{\beta s}} &\ge \frac{\theta (1-s)}{V(y, d(x,y))  d(x,y)^{\beta s} V(y, (2+\delta_0)(1-\delta_0/5)r_k(x,v_n)) }\nn\\
	&\ge \frac{c_2\theta (1-s)}{V(x, d(x,y))  d(x,y)^{\beta s} V(w_n, r_k(x,v_n)) }.
\end{align}
Moreover, for any $z \in A$, since
\begin{align*}
\frac{\delta_0	d(x,z)}{10} \ge \frac{\delta_0}{10}\left( d(x,v_n) - d(v_n,w_n)-d(w_n,z) \right) \ge \left(1 - \frac{\delta_0}{5}\right)r_k(x,v_n),
\end{align*}
we deduce from \eqref{e:def-w-n} and \eqref{e:dist-y-z} that
\begin{align}\label{e:coercivity-3-4}
	r_k(x,z) \ge  \left(1 - \frac{\delta_0}{5}\right)r_k(x,v_n) > \frac{d(y,z)}{2+\delta_0}  \ge \frac{d(y,z)}{4}.
\end{align}
Combining \eqref{e:coercivity-3-1}, \eqref{e:coercivity-3-2}, \eqref{e:coercivity-3-3} and \eqref{e:coercivity-3-4}, and using \ref{e:VD2}, we arrive at
\begin{align*}
	J_{k+1}(x,y)&\ge  \int_{A\cap  N_k(x) \cap N_0(y)  } \sJ_k(x,y,z) \mu(dz) \nn\\
	&\ge  \frac{(c_1 \wedge c_2)\theta \lambda^{k} (1-s) }{V(x, d(x,y))  d(x,y)^{\beta s} V(w_n, r_k(x,v_n)) } \int_{A\cap  N_k(x) \cap N_0(y)  } \mu(dz)\nn\\
	&\ge  \frac{\sigma(c_1 \wedge c_2)\theta \lambda^{k} (1-s) \mu(A)}{2V(x, d(x,y))  d(x,y)^{\beta s} V(w_n, r_k(x,v_n)) } \ge  \frac{c_3\theta \lambda^{k} (1-s) }{V(x, d(x,y))  d(x,y)^{\beta s}},
\end{align*}
proving that \eqref{e:coercivity-3-claim} holds. The proof is complete. \qed

In the remainder of this section, we let $\lambda=c_0$  where $c_0$ is the constant in Lemma \ref{l:coercivity-3}.

\begin{lemma}\label{l:coercivity-4}
For all $x \in M$ and $k\ge 0$, we have $	\mu((N_{k+1}(x) \setminus N_k(x)) \cap B(x,2R_0))=0 $.
\end{lemma}
\pf By Lebesgue's differentiation theorem (see \cite[Theorem 1.8]{He01}), for $\mu$-a.e. $y \in N_k(x)$, we have $R_k(x,y)>0$. The result follows from Lemma \ref{l:coercivity-3}. \qed

The proof below is  based on that of \cite[Proposition 4.3]{CS20}, with some non-trivial modifications needed since $M$ is not geodesic in general.

\begin{lemma}\label{l:coercivity-6}
 There exists $a_0\in (0,1)$  depending on $\delta_0,\sigma$ and the constants in {\rm \ref{e:VD}} only such that the following holds: Let  $x \in M$ and $z \in B(x,R_0)$. Set $R:=d(x,z)/(1+\delta_0)$. For all $k\ge 0$,  we have either
	\begin{align*}
		\text{$B(z, R) \subset N_{k+1}(x)$  \quad or \quad $\frac{\mu((N_{k+1}(x) \setminus N_k(x)) \cap  B(z,2R))}{V(z,2R)} \ge a_0$}.
	\end{align*}
\end{lemma}
\pf By Lemma \ref{l:coercivity-4}, \eqref{e:ass-coercivity} and \ref{e:VD}, we have
\begin{align}\label{e:coercivity-6-0}
	\mu(N_{k}(x) \cap B(z,R)) \ge \mu(N_{0}(x) \cap B(z,R))  \ge  \sigma V(z,R) \ge c_1 V(z,2R).
\end{align} 
 Set  $\sigma':=\eps_0^4\sigma/2$ where $\eps_0$ is the constant in \eqref{e:VD-5}. We deal with two cases separately. 

\medskip

Case 1: $\mu( N_k(x)\cap B(z,R)) < (1-\sigma') V(z,R)$.

 Let $y \in  N_k(x)\cap B(z,R)$ be a Lebesgue point of $\1_{N_k(x)}$. Then there is $l^y \in (0, R-d(z,y))$ such that $ \mu(N_k(x) \cap B(y,l^y)) \ge (1-\sigma') V(y,l^y)$. Define
\begin{align*}
	\sA^y:=\left\{ (w,r) \in \overline{B(z,R-l^y)} \times [l_y,\infty) :  d(y,w)+l^y\le r\le R- d(z,w) \right\}
\end{align*}
and $	\sA_0^y:=\left\{ (w,r) \in \sV^y: \mu(N_k(x) \cap B(w,r)) \ge (1-\sigma') V(w,r)\right\}.$
We have $(y,l^y)\in \sA_0^y$ and $(z,R)\in \sA^y \setminus \sA_0^y$. Note that  $\sA^y$ is compact.
Thus,  by Lemma \ref{l:coercivity-0+} and the intermediate value theorem,  there exists $(w^y,r^y) \in \sA^y_0$ such that
\begin{align}\label{e:coercivity-6-2}
\mu(N_k(x) \cap B(w^y,r^y)) = (1-\sigma') V(w^y,r^y).
\end{align} 
Note that $B(w^y,r^y) \subset B(z,R)$.
We claim that there exists $c_2>0$ depending on $\delta_0,\sigma$ and the constants in \ref{e:VD} only such that 
\begin{align}\label{e:coercivity-6-claim}
	\mu\big((N_{k+1}(x) \setminus N_k(x) ) \cap B(w^y,2r^y)\big) \ge c_2 V(w^y,2r^y).
\end{align}
To prove \eqref{e:coercivity-6-claim}, we consider the following three subcases separately:

\smallskip

(i) Assume $r^y\le \delta_0d(x,w^y)/11$. For any $v \in B(w^y,r^y)$, we have $r^y\le \delta_0d(x,v)/10$. Hence, by \eqref{e:coercivity-6-2}, $r_k(x,v)>0$ for all $v \in B(w^y,r^y)$. By Lemma \ref{l:coercivity-3}, this implies  $B(w^y,r^y) \subset N_{k+1}(x)$. Thus, using  \eqref{e:coercivity-6-2} and \ref{e:VD}, we get
\begin{align}\label{e:coercivity-6-claim-case-1}
		&\mu\big((N_{k+1}(x) \setminus N_k(x) ) \cap B(w^y,2r^y)\big)  \ge \mu\big(  B(w^y,r^y) \setminus N_k(x)\big)  = \sigma' V(w^y,r^y) \ge c_3 V(w^y,2r^y).
\end{align}

(ii) Assume $r^y>\delta_0d(x,w^y)/11$ and there is a covering $\{B(w_i, r_i)\}_{i=1}^\infty$ of $B(w^y,r^y)$ with balls such that for all $i\ge 1$,
$$r_i \le \frac{\delta_0^3}{132}d(x,w_i)\quad \text{and} \quad \mu( N_k(x) \cap B(w_i,r_i) ) \ge (1-\eps_0^2\sigma/2) V(w_i,r_i).$$  
By Lemma \ref{l:coercivity-3}, we see $B(w_i,r_i) \subset N_{k+1}(x)$ for all $i\ge 1$. Hence,  $B(w^y,r^y) \subset \cup_{i=1}^\infty B(w_i,r_i) \subset N_{k+1}(x)$ and \eqref{e:coercivity-6-claim-case-1} remains valid.

(iii) Assume $r^y>\delta_0d(x,w^y)/11$ and there is no covering as in (ii). Then there exists $B(w_0,r_0)$ with $w_0 \in B(w^y,r^y)$ and $r_0:= \delta_0^3d(x,w_0)/132$ such that $\mu( N_k(x)\cap B(w_0,r_0)) < (1-\eps_0^2\sigma/2) V(w_0,r_0)$. Note that
\begin{align}\label{e:r-0-r-y}
	r_0 \le \frac{\delta_0^3}{132}(d(x,w^y) + r^y) < \frac{\delta_0^2}{11}r_y.
\end{align}
By Lemma \ref{l:coercivity-0} and \eqref{e:coercivity-6-2}, there exists $w_0' \in B(w^y,r^y)$ such that $\mu( N_k(x)\cap B(w_0',r_0)) \ge  (1-\eps_0^2\sigma/2) V(w_0,r_0)$. By the continuity of $w \mapsto \mu(N_k(x)\cap B(w,r_0))/V(w,r_0)$, we deduce that there exists $v_0 \in \overline{B(w^y,r^y)}$ such that  $\mu( N_k(x)\cap B(v_0,r_0)) =  (1-\eps_0^2\sigma/2) V(v_0,r_0)$. Since 
\begin{align*}
	\frac{\delta_0}{11}d(x,v_0) \ge \frac{\delta_0}{11} (d(x,z) - d(z,v_0)) \ge \frac{\delta_0^2 }{11}R  \ge \frac{\delta_0^2}{11}r_y >r_0
\end{align*}
by \eqref{e:r-0-r-y}, we have $r_k(x,v)>0$ for all $v \in B(v_0,r_0)$.
Applying Lemma \ref{l:coercivity-3}, we get $B(v_0,r_0) \subset N_{k+1}(x) \cap B(w^y, 2r^y)$. Note that 
\begin{align*}
	r_0 \ge \frac{\delta_0^3}{132}(d(x,z)-d(z,w_0)) \ge \frac{\delta_0^4}{132}R \ge \frac{\delta_0^4}{132} d(w^y,v_0).
\end{align*}
Hence, by \ref{e:VD2}, we get $V(v_0,r_0) \ge c_3V(w^y,2r^y)$. It follows that
\begin{align*}
&	\mu\big((N_{k+1}(x) \setminus N_k(x) ) \cap B(w^y,2r^y)\big)  \ge  \mu(  B(v_0,r_0) \setminus N_k(x))  = \frac{\eps_0^2\sigma}{2} V(v_0,r_0) \ge \frac{c_4\eps_0^2\sigma}{2}V(w^y,2r^y).
\end{align*}

The proof of \eqref{e:coercivity-6-claim} is complete. By Lebesgue's differentiation theorem, almost every point of $N_k(x) \cap B(z,R)$ is a Lebesgue point of $\1_{N_k(x)}$. By the Vitali covering lemma, there exists a collection $\{B(w^{y_i},2r^{y_i})\}_{i=1}^\infty$ of pairwise disjoint open balls with $y_i \in N_k(x) \cap B(z,R)$ such that $N_k(x) \cap B(z,R) \subset \cup_{i=1}^\infty B(w^{y_i},10r^{y_i})$ $\mu$-a.e. Using  \eqref{e:coercivity-6-claim}, \ref{e:VD} and \eqref{e:coercivity-6-0}, we arrive at 
\begin{align*}
&	\mu\big((N_{k+1}(x) \setminus N_k(x))  \cap B(z,2R)   \big)  \ge \sum_{i=1}^\infty 	\mu\big((N_{k+1}(x) \setminus N_k(x))  \cap B(w^{y_i},2r^{y_i})   \big)\nn\\
	&\ge c_1\sum_{i=1}^\infty 	V(w^{y_i},2r^{y_i}) \ge c_5\sum_{i=1}^\infty 	V(w^{y_i},10r^{y_i}) \ge c_5 \mu( N_k(x) \cap B(z,R)) \ge  c_1c_2 V(z,2R).
\end{align*}

Case 2: $\mu( N_k(x)\cap B(z,R)) \ge (1-\sigma') V(z,R)$. If  there exists a covering $\{B(w_i, r_i)\}_{i=1}^\infty$ of $B(z,R)$ with balls such that 
$$r_i \le \frac{\delta_0^3}{132}d(x,w_i)\quad \text{and} \quad \mu( N_k(x) \cap B(w_i,r_i) ) \ge (1-\eps_0^2\sigma/2) V(w_i,r_i) \quad \text{for all $i \ge 1$,}$$  
then,  by Lemma \ref{l:coercivity-3}, we obtain  $B(z,R) \subset \cup_{i=1}^\infty B(w_i,r_i) \subset N_{k+1}(x)$. 
If no such  covering of $B(z,R)$ exists, we can apply the argument from Case 1(iii) to deduce that there exist $v_0 \in \overline{B(z,R)}$ and $\delta_0^4 R/132 \le r_0 < \delta_0 d(x,v_0)/11$ 
such that $\mu(N_k(x)\cap B(v_0,r_0)) = (1-\eps_0^2 \sigma/2)V(v_0,r_0)$. Using \ref{e:VD2},  we obtain
\begin{align*}
	&	\mu\big((N_{k+1}(x) \setminus N_k(x) ) \cap B(z,2R)\big)  \nn\\
	& \ge  \mu(  B(v_0,r_0) \setminus N_k(x))  = \frac{\eps_0^2\sigma}{2} V(v_0,r_0) \ge  \frac{\eps_0^2\sigma}{2} V(v_0,\delta_0^4 R/132) \ge c_5 V(z,2R).
\end{align*}
The proof is complete. \qed

Let $n_0 \ge3$ be the smallest natural number such that 
\begin{align}\label{e:coercivity-K}
	a_0(n_0-1)>1
\end{align}
where $a_0\in (0,1)$ is the constant in Lemma \ref{l:coercivity-6}. 
\begin{cor}\label{c:coercivity-6}
 For all $x \in M$, we have $	B(x, R_0) \subset N_{n_0} (x)$ $\mu$-a.e.
\end{cor}
\pf Let $z \in B(x,R_0)$ and set $R:=d(x,z)/(1+\delta_0)$. If $\mu(N_{n_0}(x)\cap B(z,R))<V(z,R)$, then by Lemmas \ref{l:coercivity-4} and \ref{l:coercivity-6}, and \eqref{e:coercivity-K},    we get
\begin{align*}
V(z,2R) \ge 	\mu(N_{n_0}(x) \cap B(z,2R)) \ge \sum_{i=0}^{n_0-1} \mu( (N_{i+1}(x)-N_{i}(x)) \cap B(z,2R)) \ge \frac{n_0}{n_0-1} V(z,2R),
\end{align*}
which is a contradiction. Thus, we get $B(z,R) \subset N_{n_0}(x)$ $\mu$-a.e., implying the desired result. \qed

Now we  present the proof of Theorem \ref{t:coercivity}.

\medskip

\noindent \textbf{Proof of Theorem \ref{t:coercivity}.} Let $x_0 \in M$, $r \in (0,R_0)$ and $f \in L^2(B(x_0,r))$. Using Lemma \ref{l:coercivity-1} and Corollary \ref{c:coercivity-6}, we obtain
\begin{align*}
	&\int_{B(x_0,r)\times B(x_0,r)} (f(x)-f(y))^2 J(x,y) \mu(dx)\mu(dy)\nn\\
	 &\ge c_1^{n_0} \int_{B(x_0, 5^{-n_0}r)\times B(x_0,5^{-n_0}r)} (f(x)-f(y))^2 J_{n_0}(x,y)\mu(dx)\mu(dy)\\
	& = c_1^{n_0} \int_{B(x_0, 5^{-n_0}r)} \int_{B(x_0, 5^{-n_0}r) \cap B(x,2\cdot 5^{-n_0}r) \cap  N_{n_0}(x)} (f(x)-f(y))^2 J_{n_0}(x,y)\mu(dy) \, \mu(dx) \nn\\
	&\ge c_1^{n_0} \theta \lambda^{n_0}(1-s)  \int_{B(x_0, 5^{-n_0}r)\times B(x_0,5^{-n_0}r)} \frac{(f(x)-f(y))^2}{V(x,d(x,y))d(x,y)^{\beta s}}\mu(dx)\mu(dy),
\end{align*}
proving that \EC  \ holds with $K_0=5^{n_0}$. \qed 

\medskip

\noindent\large{\bf 	
	Declaration of interests}

\medskip

\noindent \small The authors declare that they have no known competing financial interests or personal relationships that could have appeared to influence the work reported in this paper.

\vskip 0.4truein

\noindent {\bf Soobin Cho:} Department of Mathematics,
University of Illinois Urbana-Champaign,
Urbana, IL 61801, U.S.A.
Email: \texttt{soobinc@illinois.edu}

\end{document}